\documentclass[11pt,a4paper]{article}

\usepackage{amsmath,amssymb}
\newcommand{\R}{\mathbb{R}}
\newcommand{\D}{\Delta}

\newcommand{\n}{\nabla}
\newcommand{\ep}{\epsilon}
\newcommand{\p}{\partial}

\newcommand{\de}{\delta}

\newtheorem{lem}{Lemma}

\newcommand{\Om}{\mathbb{R}^{N}}
\newcommand{\N}{\frac{N}{2}}

\newcommand{\e}{\epsilon}
\newcommand{\va}{\varphi}

\newtheorem{definition}{Definition}
\newtheorem{theorem}{Theorem}

\newtheorem{proposition}{Proposition}
\newtheorem{corollaire}{Corollary}
\newtheorem{remarka}{Remark}

\newtheorem{lemme}{Lemma}

\title{Porous media, Fast diffusion equations and the existence of global weak solution for the quasi-solution of compressible Navier-Stokes equations}
\author{Boris Haspot %\thanks{Basque Center of Applied Mathematics, Bizkaia Technology Park, Building 500,
\thanks{Ceremade UMR CNRS 7534
Universit\'e  Paris  Dauphine,
Place du MarŽchal DeLattre De Tassigny
75775 PARIS CEDEX 16 , haspot@ceremade.dauphine.fr }}
\date{}
\begin{document}
%\tableofcontents
\maketitle
\begin{abstract}
In \cite{arxiv,arxiv1,Kor,cras1,cras2}, we have developed a new tool called \textit{quasi solutions} which approximate in some sense the compressible Navier-Stokes equation. In particular it allows us to obtain global strong solution for the compressible Navier-Stokes equations with \textit{large} initial data on the irrotational part of the velocity (\textit{large} in the sense that the smallness assumption is subcritical in terms of scaling, it turns out that in this framework we are able to obtain a family of large initial data in the energy space in dimension $N=2$). In this paper we are interested in proving the result anounced in \cite{cras3} concerning the existence of global weak solution for the quasi-solutions, we also observe that for some choice of initial data (irrotationnal) the quasi solutions verify the porous media, the heat equation or the fast diffusion equations in function of the structure of the viscosity coefficients. In particular it implies that it exists classical quasi-solutions in the sense that they are $C^{\infty}$ on $(0,T)\times\R^{N}$ for any $T>0$. Finally we show the convergence of the global weak solution of compressible Navier-Stokes equations to the quasi solutions in the case of a vanishing pressure limit process.  In particular we show that for highly compressible equations the speed of propagation of the density is quasi finite when the viscosity corresponds to $\mu(\rho)=\rho^{\alpha}$ with $\alpha>1$. Furthermore the density  is not far from converging asymptotically in time to the Barrenblatt solution of mass the initial density $\rho_{0}$.%We are also going to discuss  the notion of scaling of the solution for compressible Navier-Stokes equations which justifies the notion of quasi solution.
\end{abstract}

%The title of your section 1
%Invariance du scaling nous donne ici:
%$$(\rho_{0},u_{0})\in B^{\N+\frac{2}{2\lambda-\gamma-1}}_{2,\infty}\times B^{\N-1+2\frac{\lambda-1}{2\lambda-\gamma-1}}_{2,\infty}.$$
% Text of your Version française abrégée here.
% Note you do not need to repeat here equations that you use in the
% main text - for example 'voir (3)' is quite acceptable.

%\selectlanguage{english}
% main text
\section{Introduction}
The motion of a general barotropic compressible fluid is described by the following system:
\begin{equation}
\begin{cases}
\begin{aligned}
&\p_{t}\rho+{\rm div}(\rho u)=0,\\
&\p_{t}(\rho u)+{\rm div}(\rho u\otimes u)-{\rm div}(2\mu (\rho) D(u))-\n(\lambda(\rho){\rm div}u)+\n P(\rho)=\rho f,\\
&(\rho,u)_{/t=0}=(\rho_{0},u_{0}).
\end{aligned}
\end{cases}
\label{a01}
\end{equation}
Here $u=u(t,x)\in\R^{N}$ stands for the velocity field, $\rho=\rho(t,x)\in\R^{+}$ is the density and $D(u)=\frac{1}{2}(\n u+^{t}\n u)$ the strain tensor.
The pressure $P$ is such that $P(\rho)=a\rho^{\gamma}$ with $\gamma> 1$ and $a>0$.
We denote by $\mu(\rho)$ and $\lambda(\rho)$ the two-Lam\'e viscosity coefficients depending on the density and satisfying:
\begin{equation}
\mu(\rho)>0\;\;\;2\mu(\rho)+N\lambda(\rho)\geq 0.
\label{2coeff}
\end{equation}
Throughout the paper, we assume that the space variable $x\in\R^{N}$. In this article, we are going to investigate the existence of global weak quasi solutions for the system (\ref{a01}), a notion which has been introduced in \cite{cras1,cras2,arxiv,arxiv1} in order to prove the existence of global strong solution with large initial data on the rotational and irrotational part of the velocity for the scaling of the equation (we refer to the remark \ref{Danchin} for the definition of this notion of scaling when we assume that the initial density is far away from the vacuum).
\\
Before entering in the heart of the topic and defining the notion of quasi-solutions, we would like to explain more in details the notion of invariance by scaling for compressible Navier-Stokes equations when we suppose that the density $\rho_{0}$ is in a Besov space. Up our knowledge these results are quite new and shall allow to understand the relation between the quasi solutions and the so called Barrenblatt solution for the porous media equation in terms of invariance by scaling.%we are going to explain one of the major difficulty relied with the notion of scaling.
\subsubsection*{Scaling of the equations}
A natural way to understand the  equations of fluid mechanics correspond to search self similar solution, it means that there is a scaling of the variables after which the system become stationary solutions. Precisely it holds when we set:
\begin{equation}
\begin{aligned}
&\rho(t,x)=t^{-\alpha}F(x t^{-\beta}),\\
&u(t,x)=t^{-\alpha_{1}}G(x t^{-\beta})
\end{aligned}
\label{12scaling}
\end{equation}
The exponents $\alpha, \alpha_{1}$ and $\beta$ are called \textit{similarity exponents}, and functions $F$ and $G$ are the \textit{self similar profiles}. In particular $\alpha$ and $\alpha_{1}$ are the density contraction rate and $\beta$ the space expansion rate. In the sequel we shall assume that $\mu(\rho)=\mu\rho^{\theta}$ and $\lambda(\rho)=\lambda\rho^{\theta}$ with $\theta\geq 0$. Simple calculus give when we set $\eta=x t^{-\beta}$:
%On a alors avec $\eta=x t^{-\beta}$:
$$
\begin{aligned}
&\p_{t}\rho(t,x)=-t^{-\alpha-1}(\alpha F(\eta)+\beta\n F(\eta)\cdot \eta),\\
&{\rm div}(\rho u)=t^{-\alpha-\alpha_{1}-\beta}{\rm div}(GF)(\eta).
\end{aligned}
$$
Next we have:
$$
\begin{aligned}
&\rho\p_{t}u=-t^{-\alpha-\alpha_{1}-1}(\alpha_{1} F(\eta)G(\eta)+\beta\eta\cdot\n G(\eta) F(\eta)),\\
&\rho u\cdot\n u=t^{-\alpha-2\alpha_{1}-\beta}F(\eta)G(\eta)\cdot\n G(\eta),\\
&2\mu{\rm div}(\rho^{\theta}D u)=2\mu t^{-\theta\alpha-\alpha_{1}-2\beta}{\rm div}(F^{\theta}(\eta)D G(\eta)),\\
&\n(\rho^{\lambda}{\rm div} u)=t^{-\theta\alpha-\alpha_{1}-2\beta}\n(F^{\theta}(\eta){\rm div}G(\eta)).
\end{aligned}
$$
and finally:
$$\n\rho^{\gamma}=\gamma t^{-\alpha\gamma-\beta}F^{\gamma-1}(\eta)\n F(\eta)= t^{-\alpha\gamma-\beta}\n(aF(\eta)^{\gamma}).$$
In order to ensure the existence of solution for the system(\ref{a01}) under the specific form introduced in (\ref{12scaling}) we need to assume that:
$$
\begin{cases}
\begin{aligned}
&\alpha+\alpha_{1}+\beta=\alpha+1,\\
&\alpha+\alpha_{1}+1=\alpha+2\alpha_{1}+\beta,\\
&\alpha+\alpha_{1}+1=\theta\alpha+\alpha_{1}+2\beta,\\
&\alpha+\alpha_{1}+1=\alpha\gamma+\beta,
\end{aligned}
\end{cases}
$$
which is equivalent to the following system:
$$
\begin{cases}
\begin{aligned}
&\alpha_{1}+\beta=1,\\
&(\theta-1)\alpha+2\beta=1,\\
&\alpha(\gamma-1)+\beta=\alpha_{1}+1.
\end{aligned}
\end{cases}
$$
The solution of the previous system is:
\begin{equation}
\begin{aligned}
&\alpha=\frac{-1}{\theta-\gamma},\;\alpha_{1}=\frac{1-\gamma}{2(\theta-\gamma)}\;\;\mbox{et}\;\;\beta=\frac{2\theta-\gamma-1}{2(\theta-\gamma)}.
\end{aligned}
\label{scal}
\end{equation}
With this choice on the parameter $\alpha$, $\alpha_{1}$ and $\beta$ we finaly get the profile equation:
$$
\begin{aligned}
\begin{cases}
&\alpha F(\eta)+\beta\n F(\eta)\cdot \eta-{\rm div}(GF)(\eta)=0,\\
&\alpha_{1} F(\eta)G(\eta)+\beta\eta\cdot\n G(\eta) F(\eta)-F(\eta)G(\eta)\cdot\n G(\eta)+2\mu {\rm div}(F^{\theta}(\eta)D G(\eta))\\
&\hspace{6cm}+\n(F^{\theta}(\eta){\rm div}G(\eta))-\n(aF(\eta)^{\gamma})=0.
\end{cases}
\end{aligned}
$$
\begin{remarka}
A first remark consists in observing that there is no scaling invariance when $\theta=\gamma$.
\end{remarka}
An other way to express this scaling invariance corresponds to consider a classical solution $(\rho,u)$ of the system (\ref{a01}) and to check that the family:
$$(\rho_{l},u_{l})(t,x)=(l^{\alpha}\rho(lt,l^{\beta}x), l^{\alpha_{1}}\rho(lt,l^{\beta}x),$$
is a solution of (\ref{a01}) for any $l\in\R$ when $\alpha$, $\alpha_{1}$ and $\beta$ verify (\ref{scal}).
\begin{remarka}
We shall say that a functional space $E$ embedded in ${\cal S}^{'}(\R^{N})\times ({\cal S}^{'}(\R^{N})))^{N}$ is a critical space for (\ref{a01}) if the associated norm is invariant under the transformation
$(\rho_{0},u_{0})\rightarrow ((\rho_{0})_{l},(u_{0})_{l})$ for any $l\in\R$. In particular we observe that:
$$B^{\N+\frac{2}{2\theta-\gamma-1}}_{2,\infty}\times B^{\N-1+2\frac{\theta-1}{2\theta-\gamma-1}}_{2,\infty},$$ verifies a such property (
we refer to \cite{BCD} for the definition of Besov space) .
\end{remarka}
\begin{remarka}
Let us point out that when $\mu(\rho)=\mu\rho$ and $\lambda(\rho)=0$ with $\gamma>1$ (the case of the shallow-water equation) then $B^{\N-\frac{2}{\gamma-1}}_{2,\infty}\times B^{\N-1}_{2,\infty}$ is a critical space invariant by the scaling of the equation. In particular let us mention that it is the same scaling invariance for the initial velocity $u_{0}$ than for the incompressible Navier-Stokes equations.\\
In the case of constant viscosity coefficients we observe that $B^{\N-\frac{2}{\gamma+1}}_{2,\infty}\times B^{\N-1+\frac{2}{\gamma+1}}_{2,\infty}$ is a critical space for the scaling of the system (\ref{a01}).
\end{remarka}
\begin{remarka}
To finish we would like to mention as it has been observed in \cite{Da} that there is no invariance by scaling when we wish to work with density far from the vacuum, typically $\rho_{0}=1+q_{0}$ with $q_{0}$ in a Besov space. However in \cite{Da}  Danchin makes abstraction of the pressure term and define a other notion of critical space for the compressible Navier-Stokes equation when the initial density is far away from the vacuum. More precisely we can remark that (\ref{a01}) is invariant by the transformation:
$$
\begin{aligned}
&(\rho_{0}(x),u_{0}(x))\rightarrow (\rho_{0}(lx),l u_{0}(lx))
&(\rho(t,x),u(t,x))\rightarrow (\rho(l^{2}t,lx),lu(l^{2}t,lx))
\end{aligned}
$$
up to a change of the pressure law $P$ into $l^{2}P$. In particular $B^{\N}_{2,1}\times B^{N-1}_{2,1}$ is norm invariant by the previous transformation and is critical for the initial data $(q_{0},u_{0})$. This notion of critical space seems well adapted to the case of initial density far away from the vacuum and is also  relevant at least for the initial velocity since $B^{\N-1}_{2,1}$ is also critical for the incompressible Navier-Stokes equations. Let us recall that with this notion of critical space we have proved the existence of global strong solution with large initial data on the irrotational and rotational part of the velocity (see \cite{arxiv,arxiv1,Kor}) by involving the notion  of quasi-solutions verifying regularizing effects.
\label{Danchin}
\end{remarka}
Before coming back on this notion of scaling in order to make some link between self similar solutions for the porous media equation and the behavior of the quasi-solutions, we would like now to give precise assumptions on the viscosity coefficients with which we are going to work.
%Invariance du scaling nous donne ici:
%$$(\rho_{0},u_{0})\in B^{\N+\frac{2}{2\lambda-\gamma-1}}_{2,\infty}\times B^{\N-1+2\frac{\lambda-1}{2\lambda-\gamma-1}}_{2,\infty}.$$
%En particulier lorsque:
%$$\gamma\rightarrow (2\lambda-1)_{-}\;\;\mbox{ou}\;\;\gamma\rightarrow (2\lambda-1)_{-}.$$
%En particulier cela implique que $\lambda\geq 1$ puisque:
%$$\lambda=\frac{\gamma}{2}+\frac{1}{2}\geq 1,$$
%car $\gamma\geq 1$.\\
% or to the periodic
%box ${\cal T}^{N}_{a}$ with period $a_{i}$, in the i-th direction. 
\subsubsection*{Condition on $\mu(\rho)$ and $\lambda(\rho)$}
In this paper we are interested in extending the notion of quasi-solution developed in \cite{cras1, cras2,arxiv,arxiv1} (where only shallow water coefficients were considered) for general viscosity coefficients following the algebraic equality discovered by Bresch and Desjardin in \cite{BD,BD1}:
\begin{equation}
\lambda(\rho)=2\rho\mu^{'}(\rho)-2\mu(\rho).
\label{BD}
\end{equation}
In the sequel we shall deal with the function $\va(\rho)$ defining by $\va^{'}(\rho)=\frac{2\mu^{'}(\rho)}{\rho}$ and the function $f(\rho)$ by $f^{'}(\rho)=\sqrt{\rho}\va^{'}(\rho)$. Let us mention that the equality (\ref{BD}) implies that the viscosity coefficients are degenerated inasmuch as it imposes that $\mu(0)=\lambda(0)=0$.  With this choice of viscosity coefficients Bresch and Desjardin have obtained a remarkable new entropy for compressible Navier-Stokes equations (\ref{a01}) providing a $L^{\infty}_{T}(L^{2}(\R^{N}))$ control for any $T>0$ on the gradient of the density (more precisely on $\sqrt{\rho}\n\va(\rho)$). In particular it allows them to prove the existence of global weak solution for a specific choice of pressure, more precisely what they describe as a cold pressure. Compared with the case of viscosity coefficient the pressure term is quite simple to deal with by using Sobolev embedding since we have uniform estimate on $\sqrt{\rho}\n\va(\rho)$ in $L^{\infty}_{T}(L^{2}(\R^{N})$; however a new difficulty appears coming  from the degenerescence of the viscosity coefficient. Indeed we lose the control of $\n u\in L^{2}((0,T)\times\R^{N})$ what makes delicate the compactness study of the term $\rho u\otimes u$ (in particular via the classical energy estimates we have only a convergence in the sense of the measure) due to the existence of vacuum. In order to overcome this difficulty Mellet and Vasseur in \cite{MV} obtained new entropy on the velocity which furnish them a gain of integrability on the velocity. With this new ingredient they are able in \cite{MV}   to prove the stability of the global weak solution for compressible Navier Stokes equations with such viscosity coefficients and with classical $\gamma$ law pressure $P(\rho)=a\rho^{\gamma}$ with $a>0$ and $\gamma> 1$.\\
We are going to detail here the assumptions  on the viscosity coefficients which allow to Mellet and Vasseur  to obtain additional informations on the integrability of the velocity and we are going even to relax their hypothesis (it will be crucial in order to consider in the sequel quasi solution verifying porous media equations). We shall suppose the following inequalities on $\mu$ and $\lambda$, let $\nu_{1}\in (0,1)$ and $\nu_{2}>0$ such that:
%In the sequel we will work only with such viscosity coefficients verifying the relation (\ref{BD}).\\
%The entropy discovered in \cite{BD} shall play a crucial role but also an other one derived by Mellet and Vasseur in \cite{MV} which provides a gain of integrability on the velocity. In particular in order to obtain this entropy we shall assure as in \cite{MV} the following conditions, there exists a positive constant $\nu_{1}\in(0,1)$ and $\nu_{2}>0$ such that:
\begin{equation}
\begin{aligned}
&|\lambda^{'}(\rho)|\leq\frac{1}{\nu_{1}}\mu^{'}(\rho),\\
&\nu_{1}\mu(\rho)\leq 2\mu(\rho)+N\lambda(\rho)\leq \nu_{2}\mu(\rho).%\frac{1}{\nu}\mu(\rho).
\end{aligned}
\label{11}
\end{equation}
%When $\gamma\geq 3$ and $N=3$, we also require that:
%\begin{equation}
%\lim\inf_{\rho\rightarrow +\infty}\frac{\mu(\rho)}{\rho^{\frac{\gamma}{3}+\ep}}>0,
%\label{a11}
%\end{equation}
%for some small $\e>0$.
\begin{remarka}
If we assume that $\mu(\rho)=\mu\rho^{\alpha}$ with $\alpha> 0$ then the relation (\ref{BD}) gives:
\begin{equation}
\lambda(\rho)=2(\alpha-1)\mu\rho^{\alpha},
\end{equation}
and:
\begin{equation}
2\mu(\rho)+N\lambda(\rho)=2(1+N(\alpha-1))\mu(\rho).
\end{equation}
In this situation we have $\nu_{1}=\nu_{2}=2(1+N(\alpha-1))$.
\end{remarka}
Following \cite{MV} let us briefly make some comments  on the conditions (\ref{11}).
\begin{remarka}
The condition (\ref{11}) is crucial in order to obtain the estimates (\ref{entropcle}) and (\ref{gain}). In particular we can observe that the second condition in (\ref{11}) is similar to the classical assumption on the Lam\'e coefficient $2\mu(\rho)+N\lambda(\rho)\geq 0$ when $\mu(\rho)=\mu\rho^{\alpha}$ and $\lambda(\rho)$ verifies (\ref{BD}).
\end{remarka}
\begin{remarka}
The lower estimate in the second inequality in (\ref{11}) is trivial when $\lambda(\rho)\geq 0$, while the upper estimate is trivial when $\lambda(\rho)\leq 0$. Together this provides:
$$|\lambda(\rho)|\leq C\mu(\rho)\;\;\forall \rho>0.$$
This inequality and the first inequality of (\ref{11}) will be crucial for estimating the limit of $\n(\lambda(\rho_{n}){\rm div}u_{n})$.
\label{r1}
\end{remarka}
\begin{remarka}
Condition (\ref{11}) and (\ref{BD}) implies that:
$$\frac{N-1+\frac{\nu_{1}}{2}}{N\rho}\leq\frac{\mu^{'}(\rho)}{\mu(\rho)}\leq \frac{N-1+%\frac{1}{\nu}
\frac{\nu_{2}}{2}}{N\rho}\;\;\forall\rho>0.$$
It yields:
\begin{equation}
\begin{cases}
\begin{aligned}
&C\rho^{1-\frac{1}{N}+\frac{\nu_{1}}{2N}}\leq \mu(\rho)\leq C\rho^{1-\frac{1}{N}+\frac{\nu_{2}}{2N}}\;\;\forall \rho>1,\\
&C\rho^{1-\frac{1}{N}+\frac{\nu_{2}}{2N}}\leq \mu(\rho)\leq C\rho^{1-\frac{1}{N}+\frac{\nu_{1}}{2N}}\;\;\forall \rho\leq1,
\end{aligned}
\end{cases}
\label{13}
\end{equation}
\label{r2.2}
\end{remarka}
%\\
%\subsubsection*{Notion of quasi-solutions for compressible Navier-Stokes equations}
%This paper is devoted to prove the existence of quasi solutions for compressible Navier-Stokes equations with degenerate viscosity coefficients. 
We can now recall briefly the definition of the quasi solutions introduced in \cite{arxiv,arxiv1,cras1,cras2} which doughly spaaking are solutions of the compressible Navier-Stokes equations (\ref{a01}) where we have cancelled out the pressure term (in the sequel we shall give a more accurate definition).
\begin{definition}
We say that $(\rho,u)$ is a quasi solution if $(\rho,u)$ verifies in distribution sense:
\begin{equation}
\begin{cases}
\begin{aligned}
&\frac{\p}{\p t}\rho+{\rm div}(\rho u)=0,\\
&\frac{\p}{\p t}(\rho u)+{\rm div}(\rho
u\otimes u)-\rm div(2\mu(\rho)\, Du)-\n(\lambda(\rho){\rm div}u)=0,\\
&(\rho,u)_{\ t=0}=(\rho_{0},u_{0})
\end{aligned}
\end{cases}
\label{3systeme}
\end{equation}
\end{definition}
As we explained previously this notion of quasi solution is interesting inasmuch as it allows to exhibit \textit{large} initial velocity on the irrotational part in the sense of the scaling of the remark \ref{Danchin} (in particular we assume no vacuum on the initial density) providing global strong solution for compressible Navier-Stokes equation (see \cite{arxiv, arxiv1,Kor} for more details). This result is based on a strange phenomena on the quasi solution since these last one verify regularizing effects allowing to neglect in terms of scaling the pressure term in high frequencies. An other way to express the things is that the quasi solutions preserves a structure of irrotationality of the system when we choose irrotational initial data (it will be the case in the sequel), we shall say that the quasi solution, typically $u=-\n\va(\rho)$ is purely compressible.%by taking into account the non linear term coming from the convection terms.Indeed in this paper we are working around the quasi solutions which exhibit regular effect of type heat equation, and by combining the lack of scaling invariance on the pressure we are able to choose large initial data $q_{0}\in B^{\N}_{2,1}$ with $\rho_{0}=1+q_{0}$. In particular we are able to prove the global existence of strong solution for (\ref{a01}) with large initial data in the energy space when $N=2$. 
For other results on the existence of strong solution with critical initial data for variable viscosity coefficients we refer to \cite{CH,arma,JDE}.\\
\\
%pin \cite{arxiv}  we obtain global strong solution with initial data small in subcritical space for the scaling of the equations. In this sense quasi solutions are good approximate in order to obtain global strong solution with large initial data  in terms of scaling (in particular in dimension $N=2$ we can choose large initial data in energy space).
% Let $\rho(t,x)$ and $u(t, x)$ denote the density and the velocity of a barotropic compressible viscous fluid (as usual, $\rho$ is a non-negative function and $u$ is a vector field in $\R^{N}$ with $N\geq 2$).
%\begin{remarka}
%Here $\lambda$ and $\mu$ verifies the condition (\ref{BD}), in particular we have classical energy estimates by multiplying the momentum equation by $u$ except that we have no information on the density of the type $\rho^{\gamma}\in L^{\infty}(\R^{+},L^{1}(\R^{N}))$ as for compressible Navier-Stokes equation when ($P(\rho)=\rho^{\gamma}$ with $\gamma\geq 1$) because here $P(\rho)=0$. However using the entropy discovered in \cite{BD} we can prove that $\sqrt{\rho}\n \va(\rho)$ belongs in $L^{\infty}(\R^{+},L^{2}(\R^{N}))$ with 
%\begin{equation}
%\va^{'}(\rho)=\frac{2\mu^{'}(\rho)}{\rho}
%\label{BD1}
%\end{equation}
 %and that $\rho$ belongs in $L^{\infty}(\R^{+},L^{1}(\R^{N}))$ by conservation of the mass. It will be sufficient to prove the stability of global weak solution and it explains the assumption of the definition \ref{defMV}.
%\end{remarka}
%\begin{remarka}
%\end{remarka}
We now are going to investigate the existence of such quasi solution for the viscosity coefficients verifying (\ref{BD}) when the initial data is assumed to be close from the vacuum, typically $\rho_{0}\in L^{1}(\R^{N})$. More precisely as in \cite{arxiv,arxiv1,Kor} we are going to search in a first time irrotational solution under the form $u(t,x)=\n c(t,x)$ for the system (\ref{3systeme}). Let us assume now to simplify that:
\begin{equation}
\mu(\rho)=\mu\rho^{\alpha}\;\;\mbox{with}\;\alpha>0\;\;\mbox{and}\;\;\lambda(\rho)=2(\alpha-1)\mu\rho^{\alpha},
\end{equation}
with $\alpha\geq 1-\frac{1}{N}$ in order to ensure the relation $2\mu(\rho)+N\lambda(\rho)>0$.
We observe here that $\mu(\rho)$ and $\lambda(\rho)$ verify the relation (\ref{BD}). In this case we will verify (see the theorem \ref{theo1}) that at least for suitable initial data on $\rho_{0}$ then it exists a explicit solution to the problem (\ref{3systeme}) written under the form $(\rho,-\n\va(\rho))$ with $\rho$ verifying the porous media or the fast diffusion equation when $\alpha\ne 1$:
\begin{equation}
\begin{cases}
\begin{aligned}
&\p_{t}\rho-2\mu\D\rho^{\alpha}=0,\\
&\rho(0,\cdot)=\rho_{0}.
\end{aligned}
\end{cases}
\label{P}
\end{equation}
In a very surprising way it means that the quasi solutions are directly related to the porous or the fast diffusion equations. Moreover we will show that when we work with highly compressible Navier-Stokes equations (which correspond to the case where $a$ goes to $0$ with $P(\rho)=a\rho^{\gamma}$), then the properties of porous media or fast diffusion equations are more or less preserved for the solution of (\ref{a01}) when $a$ is small. Before giving more details on your results, let us give few words on the porous and fast diffusion equations for the reader which are not so familiar with these equations.
\subsubsection*{Porous media and fast diffusion equations}
%
%\\
%Let us now briefly recall the so-called porous medium and fast diffusion equations before explaining the link between the quasi-solutions and these solutions. More precisely the solutions of the following  the nonlinear Cauchy problem:
%\begin{equation}
%\begin{cases}
%\begin{aligned}
%&\p_{t}\rho-\mu\D\rho^{\alpha}=0,\\
%&\rho(0,\cdot)=\rho_{0}.
%\end{aligned}
%\end{cases}
%\label{P}
%\end{equation}
%where $\alpha$ is a positive number which we assume different from one are solutions of the porous media equation or fast diffusion equations,  here we assume that $\rho_{0}\in L^{1}(\R^{N})$ is nonnegative.
Let us consider the equation (\ref{P}) when $2\mu=1$ to simplify the notations; the case $\alpha>1$ (the porous media equations) arises as a model of slow diffusion of a gas inside a porous container. Unlike the heat equation $\alpha = 1$, this equation exhibits finite speed of propagation in the sense that solutions associated to compactly supported initial data remain compactly supported in space variable at all times (see \cite{Vaz} and \cite{Aron}). When $0 < \alpha < 1$, the opposite happens. Infinite speed of propagation occurs and solutions may even vanish in finite time. This problem is usually referred to as the fast diffusion equation.\\ 
Let us recall that there exists a theory of global unique solution for initial data $\rho_{0}\in L^{1}(\R^{N})$ (see the section \ref{porous} for some reminders). %which are not classical it means not $C^{\infty}$ even if the initial data is $C^{\infty}$ (see a example due to Aronson in the problem 5.7 of \cite{Vaz}). However if the initial data is strictly positive then the unique global weak solution is classical (see the proposition 7.21 p 177 in \cite{Vaz}).\\
%\end{remarka}
%\begin{remarka}
M. Pierre in \cite{Pierre} has extend this last one in obtaining the existence of unique global weak solution with bounded Radon measure as initial data. 
Let us mention also that the porous media equations  are invariant by scaling, more precisely we can introduce a notion of self similarity (for more informations we refer to the chapter 16 of \cite{Vaz}). The notion of scaling consists in searching some solutions under the following form $\rho(t,x)=t^{-\gamma}F(\frac{x}{t^{\beta}}),$
with $\gamma$ and $\beta$ to be determined. In our case $\gamma$ and $\beta$ have the form:
$\gamma(\alpha-1)+2\beta=1,$ and $F$ verifies the following equation:
$$\D F^{\alpha}+\beta\eta\dot\n F+\gamma F=0.$$
In this case we said that $\rho$ is a self similar solution of type I or a forward self similar solution. In particular it exists self similar solution such that the initial data is a Dirac mass (as in the theorem of M. Pierre in \cite{Pierre}) and such that for $t>0$ this solution conserves a constant mass. This is the so-called Barrenblatt solutions (here $\alpha>1$) that we can write under the following form:
$$U_{m}(t,x)=t^{-\gamma_{1}}F(\frac{x}{t^{\beta}})\;\;
\mbox{with}\;\;F(x)=(C-\frac{(\alpha-1)\gamma_{1}}{2\alpha}|x|^{2})_{+}^{\frac{1}{\alpha-1}}$$
with $C>0$ and $\gamma_{1}=\frac{N}{N(\alpha-1)+2}$, $\beta=\frac{1}{N(\alpha-1)+2}$. Here we have the conservation of the mass $\int U_{m}(t,x)dx=m$ with $m$ depending on $C$ and the initial data corresponds to the Dirac mass $m\delta_{0}$. Similarly when $m_{c}<\alpha<1$ with $m_{c}=\max(0,\frac{N-2}{N})$ it exists also Barrenblatt solutions defined as follows:
$$U_{m}(t,x)=t^{-\gamma_{1}}F(x t^{-\beta})\;\;\mbox{with}\;\;F(x)=(C+\kappa_{1}|x|^{2})^{\frac{-1}{\alpha-1}}_{+},$$
with $\kappa_{1}=\frac{(1-\alpha)\gamma_{1}}{2N\alpha}$. 
We recall that asymptotically in time all the global weak solution with $L^{1}$ initial data converges to a Barrenblatt solution determined by his mass $\|u_{0}\|_{L^{1}}$ (we refer to Friedman and Kamin \cite{13}, V\' azquez and Kamin \cite{18,19} and Dolbeaut and Del Pino \cite{Jean}). As we mentioned previously in the case of fast diffusion equation $0<\alpha<1$, infinite propagation occurs and solution may even vanish in finite time when $\alpha$ is in the interval $(0,m_{c})$ with $m_{c}=\max(0,\frac{N-2}{N})$. In particular it implies a lost of the initial mass when $\rho_{0}$ is in $L^{1}$ (it implies also a lost of the regularity of the solution). We refer to \cite{Vaz1} theorem 5.7 for a necessary condition of extinction, in particular the initial data belongs in an appropriate Marcinkewitz space $M_{p^{*}}(\R^{N})$. Let us finished this subsection by mentioning that we shall recall more results on the porous and the fast diffusion equations in the section \ref{porous}.
%Let us mention that in the case $\alpha\in(m_{c},1)$ the situation is quite different as the mass is preserved which implies no extinction in finite time. Moreover we have self similar solutions also discovered by Barrenblatt that we can writte under the following form:
%$$U_{m}(t,x)=t^{-\gamma_{1}}F(x t^{-\beta})\;\;\mbox{with}\;\;F(x)=(C+\kappa_{1}|x|^{2})^{\frac{-1}{\alpha-1}}_{+},$$
%with $\kappa_{1}=\frac{(1-\alpha)\gamma_{1}}{2N\alpha}$. 
%Similarly to the case $\alpha>1$ in the situation $\alpha_{c}<\alpha<1$ the global strong solution converges asymptotically to a Barrenblatt solution.
%M. A. Herrero and M. Pierre in \cite{HP} have proved the global existence and the uniqueness of global weak solution for initial data in $L^{1}_{loc}$ (this is in sharp contrast with the case $\alpha\geq 1$ where some growth condition at infinity is required on $u_{0}$ to provide even a local solution in time).\\
\subsection{Main results}
Let us now give your first result describing the link between quasi-solutions and the solutions of (\ref{P}).
%\begin{remarka}
%In the case $0<\alpha<\alpha_{c}=1-\frac{2}{N}$, let us give an example of solution which extinct in finite time $T>0$, here we define $\rho$ as follows:
%$$
%\begin{aligned}
%&\rho(t,x)=c_{\alpha}(\frac{T-t}{|x|^{2}})^{\frac{1}{1-\alpha}}\;\;\mbox{and}\;\;u(t,x)=-\frac{2\mu c_{\alpha}\alpha}{\alpha-1}\n (\frac{|x|^{2}}{T-t}),
%\end{aligned}
%$$
%with $c_{\alpha}^{1-\alpha}=2(N-\frac{2}{1-\alpha})$. In particular we observe a blow-up behavior of $u$ at time $T$. However as we mentioned above this case is not physical because we must have: $\alpha>1-\frac{1}{N}$.
%\end{remarka}
Finally we obtain the following theorems.
\begin{theorem}
%Let $\alpha>1-\frac{1}{N}$. 
Let $\mu(\rho)=\mu\rho^{\alpha}$ with $\alpha\geq 1-\frac{1}{N}$ and $\lambda(\rho)$ verifying the relation (\ref{BD}).\\
Let $\rho_{0}\in L^{1}(\R^{N})$ with $\rho_{0}>0$ and continuous and $u_{0}=-\n\va(\rho_{0})$. Then it exists a global weak solution solution of the system (\ref{3systeme}) of the form $(\rho,u=-\n\va(\rho))$ with $(\rho,u)$ belonging in $C^{\infty}((0,+\infty)\times\R^{N})\cap C([0,+\infty]\times\R^{N})$  and  solving the following system almost everywhere :
\begin{equation}
\begin{cases}
&\p_{t}\rho-2\D\mu(\rho)=0,\\
&\rho(0,\cdot)=\rho_{0}.
\end{cases}
\label{3systeme1}
\end{equation}
Furthermore we have:
\begin{equation}
\lim_{t\rightarrow+\infty}\|\rho(t)-U_{m}(t)\|_{L^{1}(\R^{N})}=0.
\label{convasym}
\end{equation}
Convergence holds also in $L^{\infty}$ norm:
\begin{equation}
\lim_{t\rightarrow+\infty}t^{\beta}\|\rho(t)-U_{m}(t)\|_{L^{\infty}(\R^{N})}=0,
\label{convasym1}
\end{equation}
with $\beta=\frac{N}{N(\alpha-1)+2}$ and $U_{m}$ the Barrenblatt of mass $m=\|\rho_{0}\|_{L^{1}(\R^{N})}$.
For every $p\in(1,+\infty)$ we have the following regularizing effect, $\rho(t,\cdot)$ belongs in $L^{p}(\R^{N})$ and:
$$\|\rho(t)\|_{L^{p}(\R^{N})}\leq Ct^{-\sigma_{p }}\|\rho_{0}\|^{\alpha_{p}}_{L^{1}(\R^{N})},$$
with $\sigma_{p}=\frac{N(\alpha-1)+2p}{(N(\alpha-1)+2)p}$ and $\alpha_{p}=\frac{N(p-1)}{(N(\alpha-1)+2)p}$.
%Furthermore this solution is unique in the class of the function $(\rho,u)$ belonging in $C^{\infty}((0,+\infty),\R^{N})$.\\
%When $\alpha=1$ with $(\rho_{0},u_{0}=-2\mu\n\ln\rho_{0})$,  similarly we have particular global weak unique solution solution of the system (\ref{3systeme}) of the form $(\rho,u=-2\mu\n\ln\rho)$ solving the heat equation:
%\begin{equation}
%\begin{cases}
%&\p_{t}\rho-2\mu\D\rho=0,\\
%&\rho(0,\cdot)=\rho_{0}.
%\end{cases}
%\label{3systeme2}
%\end{equation}
%And if $\alpha>1$ we assume that $u=0$ on $\{\rho=0\}$.
\label{theo1}
\end{theorem}
\begin{remarka}
Let us point out that we could have also global strong solution for more general viscosity coefficients with $\mu$ verifying the same conditions than the subsection (\ref{pormu}). We refer to the section \ref{porous} for more details in this situation.
\end{remarka}
\begin{remarka}
We shall remark in the proof of this theorem that any solution of $(\ref{P})$ such that $\rho$ is in $C^{3}((0,+\infty)\times\R^{N})$ is a solution of (\ref{3systeme}) almost everywhere.% It is in particular the case of some Barrenblatt solutions (see the remark \ref{}). 
%n the case where $\rho=0$ the velocity is not defined when $0<\alpha<1$ that is why we assume that $u=0$ on the vacuum set. In other case we could give sense to $\rho u$ as in \cite{MV}.
\end{remarka}
\begin{remarka}
We can observe as in \cite{friction} that if we consider the compressible Navier-Stokes equation with a friction term $a\rho u$ and a pressure term of the form $2\mu a\rho^{\alpha}$ then the solution of the previous theorem \ref{theo1} verify also a such system.
\end{remarka}
%\begin{remarka}
%We recognize here the so called equation of the porous medium when $\alpha>1$ and of the fast diffusion when $0<\alpha<1$. We refer for more details on the theory to the books of J-L V\' azquez (see \cite{Vaz,Vaz1}).
%\end{remarka}
\begin{remarka}
%We have recognized the so called porous medium equation and the fast diffusion equation. Unlike the heat equation $\alpha=1$, when $\alpha>1$ this equation exhibits finite speed of propagation in the sense that solutions associated to compactly supported initial data remains compactly supported in space variable at all times. The situation is different in the case of fast diffusion equation $0<\alpha<1$, indeed in this case infinite propagation occurs. 
Let us mention that when $\alpha$ is in the interval $(0,m_{c})$ with $m_{c}=\max(0,\frac{N-2}{N})$ (then this contradicts the Lam\'e condition on the viscosity coefficients ) then it can appears a phenomena of extinction of the solution in finite time, in particular it implies a loss of the initial mass when $\rho_{0}$ is in $L^{1}$. A typical example is the solution:
%n the case $0<\alpha<\alpha_{c}=1-\frac{N}{2}$, let us give an example of solution which extinct in finite time $T>0$, here we define $\rho$ as follows:
$$
\begin{aligned}
&\rho(t,x)=c_{\alpha}(\frac{T-t}{|x|^{2}})^{\frac{1}{1-\alpha}}\;\;\mbox{and}\;\;u(t,x)=-\frac{2\mu c_{\alpha}\alpha}{\alpha-1}\n (\frac{|x|^{2}}{T-t}),
\end{aligned}
$$
with $c_{\alpha}^{1-\alpha}=2(N-\frac{2}{1-\alpha})$. In particular we observe a blow-up behavior of $u$ at time $T$.%
% ( it is not the case in our framework because $\alpha1-\frac{1}{N}>m_{c}$). %Let us recall that if the initial data is strictly positive then the unique global weak solution is classical (see the proposition 7.21 p 177 in \cite{Vaz}). 
%We refer to \cite{Vaz,Vaz1} for more details on such equations. In particular it implies a lost of the initial mass when $\rho_{0}$ is in $L^{1}$ (it implies also a lost of the regularity of the solution, it is not the case in our framework because $\alpha>1-\frac{1}{N}>m_{c}$).\\
%Let us recall that when $\alpha>1$ and $\rho_{0}$ belongs in $L^{1}$, it exists global unique weak solution and that the solution converges asymptotically to the so called Barrenblatt solution
%$$U_{m}(t,x)=t^{-\gamma}F(\frac{x}{t^{\beta}})\;\;
%\mbox{with}\;\;F(x)=(C-\frac{\alpha-1}{2\alpha}|x|^{2})_{+}^{\frac{1}{\alpha-1}},$$
%which are self similar. 
%finally that there exists global weak solution which are not classical it means not $C^{\infty}$ even if the initial data is $C^{\infty}$ (see a example due to Aronson in the problem 5.7 of \cite{Vaz}).  However if the initial data is strictly positive then the unique global weak solution is classical (see the proposition 7.21 p 177 in \cite{Vaz}). We refer to \cite{Vaz,Vaz1} for more details on such equations.
\end{remarka}
%\begin{corollaire}
%Barrenblatt avec support compact
%\end{corollaire}
%\begin{remarka}
%SENS DE RANKINE-HUGUENOT
%\end{remarka}
We are going to give a general definition of global weak solution for the quasi solutions in the spirit of \cite{MV} including the case where the initial velocity is not necessary irrotational. %Finally we 
\begin{definition}
We say that $(\rho,u)$ is a global weak quasi solution if $(\rho,u)$ verifies in distribution sense:
\begin{equation}
\begin{cases}
\begin{aligned}
&\frac{\p}{\p t}\rho+{\rm div}(\rho u)=0,\\
&\frac{\p}{\p t}(\rho u)+{\rm div}(\rho
u\otimes u)-\rm div(2\mu(\rho)\, Du)-\n(\lambda(\rho){\rm div}u)=0,\\
&(\rho,u)_{\ t=0}=(\rho_{0},u_{0}).
\end{aligned}
\end{cases}
\label{3systemev}
\end{equation}
More precisely  $(\rho,u)$ is a weak solution of (\ref{3systeme}) on $[0,T]\times\R^{N}$%, which the following initial conditions
\begin{equation}
\rho_{/t=0}=\rho_{0}\geq0,\;\;\rho u_{/t=0}=m_{0}.
\label{7}
\end{equation}
with:
\begin{equation}
\begin{aligned}
&\rho_{0} \in L^{1}(\R^{N}),\;\sqrt{\rho_{0}}\n\va(\rho_{0})\in L^{2}(\R^{N}),\;\rho_{0}\geq 0,\\
&\sqrt{\rho_{0}}|u_{0}|(1+\sqrt{\ln(1+|u_{0}|^{2})})\in L^{2}(\R^{N}).%\;\;\;\mbox{for somme small $\delta$}.
\end{aligned}
\label{8}
\end{equation}
if
\begin{itemize}
\item $\rho\in L^{\infty}_{T}(L^{1}(\R^{N})$, $\sqrt{\rho}\n\va(\rho)\in L_{T}^{\infty}(L^{2}(\R^{N}))$, $\sqrt{\rho}u\in L^{\infty}_{T}(L^{2}(\R^{N}))$,
\item $\sqrt{\mu(\rho)}\, \n u\in L^{2}((0,T)\times \R^{N})$, $\sqrt{\rho}|u|\sqrt{\ln(1+|u|^{2})}\in L^{\infty}_{T}(L^{2}(\R^{N}))$.\\
\end{itemize}
with $\rho\geq 0$ and $(\rho,\sqrt{\rho}u)$ satisfying in distribution sense on $[0,T]\times\R^{N}$:
$$
\begin{cases}
\begin{aligned}
&\p_{t}\rho+{\rm div}(\sqrt{\rho}\sqrt{\rho}u)=0,\\
&\rho(0,x)=\rho_{0}(x).
\end{aligned}
\end{cases}
$$
and if the following equality holds for all $\va(t,x)$ smooth test function with compact support such that $\va(T,\cdot)=0$:
\begin{equation}
\begin{aligned}
&\int_{\R^{N}}(\rho u)_{0}\cdot\va(0,\cdot)dx+\int^{T}_{0}\int_{\R^{N}}\sqrt{\rho}(\sqrt{\rho}u)\p_{t}\va+
\sqrt{\rho}u\otimes\sqrt{\rho}u:\n\va dx\\
&\hspace{4cm}-<2\mu(\rho)\, Du,\n \va>-<\lambda(\rho)\, {\rm div}u,{\rm div} \va>=0,
\end{aligned}
\label{equadistrib}
\end{equation}
where we give sense to the diffusion terms by rewriting him according to $\sqrt{\rho}$ and $\sqrt{\rho}u$:
$$
\begin{aligned}
&<2\mu(\rho)\, D u,\n\va>=-\int\frac{\mu(\rho)}{\sqrt{\rho}}(\sqrt{\rho}u_{j})\p_{ii}\va_{j}dx\,dt-\int
2(\sqrt{\rho}u_{j})\mu^{'}(\rho)\p_{i}\sqrt{\rho}\p_{i}\va_{j}dx\,dt\\
&-\int\frac{\mu(\rho)}{\sqrt{\rho}}(\sqrt{\rho}u_{j})\p_{ji}\va_{j}dx\,dt-\int
2(\sqrt{\rho}u_{i})\mu^{'}(\rho)\p_{j}\sqrt{\rho}\p_{i}\va_{j}dx\,dt\\[2mm]
&<\lambda(\rho)\, {\rm div}u,{\rm div} \va>=-\int\frac{\lambda(\rho)}{\sqrt{\rho}}(\sqrt{\rho}u_{i})\p_{ji}\va_{j}dx\,dt-\int
2(\sqrt{\rho}u_{i})\lambda^{'}(\rho)\p_{i}\sqrt{\rho}\p_{j}\va_{j}dx\,dt\\
%&-\int\sqrt{\rho}(\sqrt{\rho}v_{i})\p_{ji}\va_{j}dx\,dt-\int
%2\sqrt{\rho}v_{i}\p_{j}\sqrt{\rho}\p_{i}\va_{j}dx\,dt.
\end{aligned}
$$
We assume also that $\mu$ and $\lambda$ verify the conditions (\ref{BD}) and (\ref{11}).
%Similarly we have exactly the same type of definition for system (\ref{3systeme2}).%, except that we write the diffusion term on the following form:
%$$
%\begin{aligned}
%&<\rho\n v,\n\va>=-\int\sqrt{\rho}(\sqrt{\rho}v_{j})\p_{ii}\va_{j}dx\,dt-\int
%2\sqrt{\rho}v_{j}\p_{i}\sqrt{\rho}\p_{i}\va_{j}dx\,dt.
%\end{aligned}
%$$
%\label{defMV}
%We supplement the problem with initial condition $(\rho_{0},u_{0})$.
\label{defMV}
\end{definition}
We obtain now a general result concerning the stability of the global weak solution for system (\ref{3systeme}) and a result of existence of global weak solution for general initial data of the form $(\rho_{0},-\n\va(\rho_{0}))$ (in particularly $\rho_{0}$ is not assumed only strictly positive).
%. The first one show the existence of global weak solution for system (\ref{3systeme1}), the second one for system (\ref{3systeme2}) and the last one corresponds to the existence of global weak solution for isentropic compressible Navier-Stokes equation, it means the following system:
%\begin{equation}
%\begin{cases}
%\begin{aligned}
%&\frac{\p}{\p t}\rho+{\rm div}(\rho u)=0,\\
%&\frac{\p}{\p t}(\rho u)+{\rm div}(\rho
%u\otimes u)-\rm div(\mu\rho\,\n u)+\n(a\rho^{\gamma})=0.
%\end{aligned}
%\end{cases}
%\label{3systeme3}
%\end{equation}
%with $\gamma\geq 1$. We show in the following result the existence of global weak solution for system (\ref{3systeme1}).
\begin{theorem}
\label{theo2}
%Let $1<\gamma<p$ with $p=+\infty$ if $N=2$ and $p=3$ if $N=3$. 
Assume that $\mu(\rho)$ and $\lambda(\rho)$ are two regular function of $\rho$ verifying (\ref{BD}) and (\ref{11}). Furthermore we shall set $g(x)=\frac{\mu(x)}{\sqrt{x}}$ and we assume that $g$ is bijective and that $g^{-1}$ is continuous on $(0,+\infty)$. When $2+N \leq \nu_{1} $, we assume in addition that $g$ and  $g^{'}$ are increasing on $(0,+\infty)$.\\
Let $(\rho_{n},u_{n})$  be a sequence of global weak solutions of  system (\ref{3systemev}) satisfying entropy inequalities (\ref{21}), (\ref{22}) and (\ref{gain}) with initial data:
$$(\rho_{n})_{/t=0}=\rho_{0}^{n}(x)\;\;\;\mbox{and}\;\;\;(\rho_{n}u_{n})_{/t=0}=\rho_{0}^{n}u^{n}_{0}(x)$$
with $\rho_{0}^{n}$ and $u^{n}_{0}$ such that:
\begin{equation}
\rho_{0}^{n}\geq0,\;\;\rho_{0}^{n}\rightarrow\rho_{0}\;\;\mbox{in}\;L^{1}(\R^{N}),\;\;\rho_{0}^{n}u^{n}_{0}\rightarrow\rho_{0}u_{0}\;\;\mbox{in}\;L^{1}(\R^{N}),
\label{15a}
\end{equation}
and satisfying the following bounds (with $C$ constant independent on $n$):
\begin{equation}
\int_{\R^{N}}\rho_{0}^{n}\frac{|u_{0}^{n}|^{2}}{2}<C,\;\;\;\;\;
\int_{\R^{N}}\sqrt{\rho_{0}^{n}}|\n \va( \rho_{0}^{n})|^{2}dx<C
\label{16a}
\end{equation}
and:
\begin{equation}
\int_{\R^{N}}\rho_{0}^{n}\frac{1+|u_{0}^{n}|^{2}}{2}\ln(1+|u_{0}^{n}|^{2}) dx<C.
\label{17a}
\end{equation}
%\begin{equation}
%\int_{\Om}\big(\rho_{0}^{n}\frac{|u_{0}^{n}|^{2}}{2}+\frac{a(\rho_{0}^{n})^{\gamma}}{\gamma-1})<C,\;\;
%\int_{\Om}\frac{1}{\rho_{0}^{n}}|\n \ln \rho_{0}^{n}|^{2}dx<C,
%\label{16a}
%\end{equation}
%and:
%\begin{equation}
%\int_{\Om}\rho_{0}^{n}\frac{|v_{0}^{n}|^{2+\delta}}{2}dx<C,
%\label{17a}
%\end{equation}
Then, up to a subsequence, $(\rho_{n},\sqrt{\rho_{n}}u_{n})$ converges strongly to a global weak solution $(\rho,\sqrt{\rho}u)$ of (\ref{3systemev}) satisfying entropy inequalities (\ref{21}), (\ref{22}) and (\ref{gain}).\\
Furthermore the density $\rho_{n}$ converges strongly to $\rho$ in $C([0,T], L^{1+\alpha}_{loc}(\Om))$ with $0<\alpha<\nu_{1}$ when $N=3$ and in $C([0,T], L^{q}_{loc}(\Om))$ for any $q\geq 1$ when $N=2$; $\sqrt{\rho_{n}}u_{n}$ converges strongly in $L^{2}(0,T,L^{2}_{loc})$ to $\sqrt{\rho}u$ and the momentum $m_{n}=\rho_{n}u_{n}$ converges strongly in $L^{1}(0,T,L^{1}_{loc}(\Om))$, for any $T>0$.\\[2mm]
If we assume moreover that $(\rho_{0},u_{0})$ verify the initial condition of the definition \ref{defMV}, that $u_{0}=-\n\va(\rho)$  with $\mu(\rho)=\mu\rho^{\alpha}$ and that $\sqrt{\rho_{0}}|u_{0}|^{1+\frac{1}{p}}$ belongs in $L^{2}(\R^{N})$ for $p$ large enough, then it exists a global weak solution $(\rho,u)$ of the system (\ref{3systemev}) where $\rho$ is also the unique solution of the system (\ref{P}) (see the theorem \ref{theo2.4} for the existence of a unique global strong solution for the equation (\ref{P}) with a $L^{1}$ initial data).% (as limit solution of energy solutions of (\ref{3systeme1}) (see chapter 6 in \cite{Vaz}).
\end{theorem}
\begin{remarka}
Let us mention that the additional technical assumption on the viscosity $\mu$ remains quite natural since they are verified in the standard case $\mu(\rho)=\mu\rho^{\alpha}$ with $\alpha>1-\frac{1}{N}$.
In particular we remark that this result extend the analysis of \cite{MV} to general viscosity coefficient, in particular we do not suppose that $\mu^{'}(\rho)\geq c>0$.
\end{remarka}
\begin{remarka}
Let us emphasize that the condition (\ref{11}) implies that we exclude the case of fast diffusion equation with $0<\alpha<1-\frac{1}{N}$, in particular it forbids any phenomen of extinction and loss of mass $\|\rho(t,\cdot)\|_{L^{1}(\R^{N})}$.
%on the importance of the gain of integrability on $v$. This is due to the specific structure of the diffusion on $v$ which is of the form ${\rm div}(\rho\n v)$ and not under the form  ${\rm div}(\rho\n v)$. Following Mellet and Vasseur in \cite{fMV} we obtain a supplementary entropy on the effective velocity $v$ which is a angular stone of the proof. It also shall play a crucial role in the theorem \ref{theo4}.\\
%Let us emphasize that when $N=3$, we have the restriction $1<\gamma<3$ essentially because in other case we are not able to derive a gain of integrability on $v$.
\end{remarka}
\begin{remarka}
Let us mention that our second result of existence of global weak solution when $u_{0}=-\n \va(\rho_{0})$ can be applied to the Barrenblatt solution by choosing $\rho_{0}=U_{m}(\tau,\cdot)$ with $\tau>0$. In this case we have then a free boundary problem with $(\rho,u)$ $C^{\infty}$ when $\rho>0$ and $(\rho,u)=(0,0)$ when $\rho=0$. In particular when $\mu(\rho)=\mu\rho^{\alpha}$ with $\alpha>1$ we observe that the propagation of the support of the Barrenblatt is finite.\\
More generally when we choose a initial density with compact support the support of the density remains bounded along the time when we deal with the case $\mu(\rho)=\mu\rho^{\alpha}$ with $\alpha>1$. It means that quasi-solutions in this case have the same properties than the porous media solutions. We are related to a free boundary problem with on a side the solution which is $C^{\infty}$ and on the other side the solution which is identically null. In particular it implies that we can not hope uniqueness as the velocity can take any value on the vacuum set, that is why it is natural to consider the momentum unknown $\sqrt{\rho}u$ as it is the case in the previous theorem.
\end{remarka}
%\begin{remarka}
%Let us mention that when we have initial data $\rho_{0}$ with compact support than any $\rho_{0}=U_{m}(\tau,\cdot)$ with $\tau>0$ then we can not hope the uniqueness of the quasi-solution. Indeed on the region $\{(t,x),\;\rho(t,x)=0\}$ the velocity $u$ can take any value. There is a problem to define the uniqueness, an option to overcome this difficulty consists in working with the unknowns $\rho$ and $\sqrt{\rho}u$.
%\end{remarka}
\begin{remarka}
Let us emphasize on the fact that the problem of the existence of global weak solutions remains open in the general case (it means when $u_{0}$ is different from $-\n\va(\rho_{0})$). Indeed in the previous theorem we prove the stability of the global weak solutions, however it seems quite complicated to construct approximate global weak solution which verify uniformly in $n$ all the entropies. We have the same problem in the case of the compressible Navier Stokes problem where the problem remains also open.
\end{remarka}
We are going to prove now the convergence of the global weak solution of the compressible Navier-Stokes equations to the quasi solutions when we consider a vanishing pressure process. More precisely let us consider the highly compressible Navier-Stokes system with $\ep$ going to $0$:
\begin{equation}
\begin{cases}
\begin{aligned}
&\p_{t}\rho_{\ep}+{\rm div}(\rho_{\ep} u_{\ep})=0,\\
&\p_{t}(\rho_{\ep} u_{\ep})+{\rm div}(\rho_{\ep} u_{\ep}\otimes u_{\ep})-{\rm div}(2\mu(\rho_{\ep}) D(u_{\ep}))-\n(\lambda(\rho_{\ep}){\rm div}u_{\ep})+\ep\n \rho_{\ep}^{\gamma}=0,\\
&(\rho_{\ep},u_{\ep})_{/t=0}=(\rho_{0},u_{0}).
\end{aligned}
\end{cases}
\label{0.2}
\end{equation}
In the literature we find many result on the incompressible limit which corresponds to take $\e=\frac{1}{\eta^{2}}$ with $\eta$ the Mach number going to $0$. For such results in the framework of the global weak solutions for the ill-prepared data we refer to the pioneering papers of Desjardins, Grenier, Lions and Masmoudi \cite{D1,D2,L1,L2}. Indeed in this last situation we can observe at least heuristically that the density $\rho_{\e}$ converges to a constant $1$ when we are working with initial density of the form $\rho_{0,\e}=1+\eta q_{0,\e}$ with $q_{0,\e}$ uniformly bounded in appropriate space. By considering the mass equation, roughly speaking we can deduce that the limit solution of $u_{\e}$ is incompressible. The main difficulty compared with the well prepared case corresponds to deal with the acoustic waves, in particular in order to overcome such difficulty the authors use Strichartz estimates.\\
In our case we are going to deal with the opposite situation when the solutions are highly compressible and converge to quasi-solution which are in some sense purely compressible since irrotationnal. In the case of constant viscosity it appears impossible to pass in the limit when $\e$ goes to $0$ since we lose the $L^{\infty}_{T}(L^{\gamma}(\R^{N}))$ estimate on $\rho_{\e}$ coming from the pressure term (we conserve only the $L^{1}$ conservation of the mass which is not sufficient to pass to the limit since it does not provide enough compactness information). In our case with viscosity coefficients verifying the relation (\ref{BD}) we know that we have uniform estimate on $\sqrt{\rho_{\e}}\n\va(\rho_{\e})$ in $L^{\infty}_{T}(L^{2}(\R^{N}))$ providing of the entropy (\ref{22a}) (see \cite{BD}) which will be sufficient in term of compactness to pass to the limit when $\e$ goes to $0$. Before giving our main result in this spirit, let us give a definition of global weak solution for the system (\ref{0.2}) in the spirit of \cite{MV}:
\begin{definition}
We say that $(\rho,u)$ is a global weak solution of (\ref{a01}) if $(\rho,u)$ verifies in distribution sense:
\begin{equation}
\begin{cases}
\begin{aligned}
&\frac{\p}{\p t}\rho+{\rm div}(\rho u)=0,\\
&\frac{\p}{\p t}(\rho u)+{\rm div}(\rho
u\otimes u)-\rm div(2\mu(\rho)\, Du)-\n(\lambda(\rho){\rm div}u)+\n P(\rho)=0,\\
&(\rho,u)_{\ t=0}=(\rho_{0},u_{0}),
\end{aligned}
\end{cases}
\label{3systemev}
\end{equation}
with $P(\rho)=a\rho^{\gamma}$, $\gamma>1$.
More precisely  $(\rho,u)$ is a weak solution of (\ref{a01}) on $[0,T]\times\R^{N}$%, which the following initial conditions
\begin{equation}
\rho_{/t=0}=\rho_{0}\geq0,\;\;\rho u_{/t=0}=m_{0}.
\label{7}
\end{equation}
with:
\begin{equation}
\begin{aligned}
&\rho_{0} \in L^{1}(\R^{N})\cap L^{\gamma}(\R^{N}),\;\sqrt{\rho_{0}}\n\va(\rho_{0})\in L^{2}(\R^{N}),\;\rho_{0}\geq 0,\\
&\sqrt{\rho_{0}}|u_{0}|(1+\sqrt{\ln(1+|u_{0}|^{2})})\in L^{2}(\R^{N}).%\;\;\;\mbox{for somme small $\delta$}.
\end{aligned}
\label{8}
\end{equation}
if
\begin{itemize}
\item $\rho\in L^{\infty}_{T}(L^{1}(\R^{N})$, $\sqrt{\rho}\n\va(\rho)\in L_{T}^{\infty}(L^{2}(\R^{N}))$, $\sqrt{\rho}u\in L^{\infty}_{T}(L^{2}(\R^{N}))$,
\item $\sqrt{\mu(\rho)}\, \n u\in L^{2}((0,T)\times \R^{N})$, $\sqrt{\rho}|u|\sqrt{\ln(1+|u|^{2})}\in L^{\infty}_{T}(L^{2}(\R^{N}))$.\\
\end{itemize}
with $\rho\geq 0$ and $(\rho,\sqrt{\rho}u)$ satisfying in distribution sense on $[0,T]\times\R^{N}$:
$$
\begin{cases}
\begin{aligned}
&\p_{t}\rho+{\rm div}(\sqrt{\rho}\sqrt{\rho}u)=0,\\
&\rho(0,x)=\rho_{0}(x).
\end{aligned}
\end{cases}
$$
and if the following equality holds for all $\va(t,x)$ smooth test function with compact support such that $\va(T,\cdot)=0$:
\begin{equation}
\begin{aligned}
&\int_{\R^{N}}(\rho u)_{0}\cdot\va(0,\cdot)dx+\int^{T}_{0}\int_{\R^{N}}\sqrt{\rho}(\sqrt{\rho}u)\p_{t}\va+
\sqrt{\rho}u\otimes\sqrt{\rho}u:\n\va dx\\
&+a\int^{T}_{0}\int_{\R^{N}}\rho^{\gamma}{\rm div}\va(s,x)dx-<2\mu(\rho)\, Du,\n \va>-<\lambda(\rho)\, {\rm div}u,{\rm div} \va>=0,
\end{aligned}
\label{equadistrib}
\end{equation}
where we give sense to the diffusion terms as in the definition \ref{defMV}. We assume also that $\mu$ and $\lambda$ verify the conditions (\ref{BD}) and (\ref{11}).
%Similarly we have exactly the same type of definition for system (\ref{3systeme2}).%, except that we write the diffusion term on the following form:
%$$
%\begin{aligned}
%&<\rho\n v,\n\va>=-\int\sqrt{\rho}(\sqrt{\rho}v_{j})\p_{ii}\va_{j}dx\,dt-\int
%2\sqrt{\rho}v_{j}\p_{i}\sqrt{\rho}\p_{i}\va_{j}dx\,dt.
%\end{aligned}
%$$
%\label{defMV}
%We supplement the problem with initial condition $(\rho_{0},u_{0})$.
\label{defMV1}
\end{definition}
\begin{theorem}
Let $\gamma>1$% with the additional hypothesis (\ref{a11})  if $\gamma\geq 3$ and $N=3$ 
and $(\rho_{0},u_{0})$ verifies the conditions of definition \ref{defMV1}. Let us take the following assumptions on $\gamma$, $\nu_{1}$ and $\nu_{2}$: 
\begin{itemize}
\item $N=3$\\
\\
$\bullet$ $\nu_{1}\geq 2$:
 $$
 \begin{aligned}
 &\frac{5}{6}+\frac{\nu_{2}}{12}<\gamma<2+\frac{\nu_{1}}{2}\;\;\;\mbox{if}\;\;N=3,\\
 &\frac{5}{6}+\frac{\nu_{2}}{12}<\gamma<\frac{5}{6}+\frac{7}{12}\nu_{1}\;\;\;\mbox{if}\;\;N=3,\\
 \end{aligned}
 $$
 $\bullet$ $0<\nu_{1}< 2$:
$$
 \begin{aligned}
 &\frac{5}{6}+\frac{\nu_{2}}{12}<\gamma<\frac{(4-\nu_{1})(1+\nu_{1})}{2-\nu_{1}}\;\;\;\mbox{if}\;\;N=3,\\
 &\frac{5}{6}+\frac{\nu_{2}}{12}<\gamma<\frac{5}{6}+\frac{7}{12}\nu_{1}\;\;\;\mbox{if}\;\;N=3,\\
% &\frac{1}{4}+\frac{\nu_{2}}{8}<\gamma\;\;\;\mbox{if}\;\;N=2.
 \end{aligned}
 $$
 \item $N=2$
 $$
 \begin{aligned}
 %&\frac{5}{6}+\frac{\nu_{2}}{12}<\gamma<\frac{(4-\nu_{1})(1+\nu_{1})}{2-\nu_{1}}\;\;\;\mbox{if}\;\;N=3,\\
% &\frac{5}{6}+\frac{\nu_{2}}{12}<\gamma<\frac{5}{6}+\frac{7}{12}\nu_{1}\;\;\;\mbox{if}\;\;N=3,\\
 &\frac{1}{4}+\frac{\nu_{2}}{8}<\gamma\;\;\;\mbox{if}\;\;N=2.
 \end{aligned}
 $$
 \end{itemize}
 Then under these conditions we have as in the theorem \ref{theo2} the stability of the global weak solutions for the system (\ref{a01}).\\
Assume that there exists global weak solution $(\rho_{\ep},u_{\ep})$ verifying the definition of \cite{MV} with the conditions (\ref{BD}) and (\ref{11}) on $\mu(\rho)$ and $\lambda(\rho)$. Then  $(\rho_{\ep},u_{\ep})$ converges in distribution sense to a global weak quasi-solution $(\rho,u)$ of the system (\ref{3systemev}) in the sense of the definition \ref{defMV}. Furthermore the density $\rho_{\e}$ converges strongly to $\rho$ in $C([0,T], L^{1+\alpha}_{loc}(\Om))$ with $0<\alpha<\nu_{1}$ when $N=3$ and in $C([0,T], L^{q}_{loc}(\Om))$ for any $q\geq 1$ when $N=2$; $\sqrt{\rho_{\e}}u_{\e}$ converges strongly in $L^{2}(0,T,L^{2}_{loc})$ to $\sqrt{\rho}u$ and the momentum $m_{\e}=\rho_{\e}u_{\e}$ converges strongly in $L^{1}(0,T,L^{1}_{loc}(\Om))$, for any $T>0$
%Furthermore $\rho_{\e}$ converges strongly to $\rho$ in $C([0,T],L^{1}_{loc}(\R^{N})$ for any $T>0$. VERIFIER PROPREMENT LES CONVERGENCES!!!!
\label{theo3}
\end{theorem}
\begin{remarka}
Let us mention that the technical assumption on $\nu_{1}$, $\nu_{2}$ and $\gamma$ are important in order to ensure a uniform gain of integrability on the velocity $u_{\e}$ (as in the paper of Mellet and Vasseur \cite{MV}), more precisely we will see that we have a competition between the pressure and the viscosity. 
\end{remarka}
\begin{remarka}
We shall also emphasize on an important question which remains open; indeed when $\mu(\rho)=\mu\rho^{\alpha}$ and $u_{0}=-\n\va(\rho_{0})$  the limit solution $(\rho,u)$ of $(\rho_{\ep},u_{\ep})$ when $\ep$ goes to $0$ is a quasi solution of (\ref{3systeme1}). However it is not clear how to prove the uniqueness of the quasi solutions in the class of the solutions giving by  the definition \ref{defMV} and in particular to show that this quasi solution is solution of the porous media equation for the density $\rho$ when $\alpha>1$. In the following corollary, we shall give properties of the solutions of (\ref{0.2}) when we assume the uniqueness of the quasi solutions.\\ 
Let us mention that in some case if we assume that the solution of the porous media is enough regular (typically the case of some solutions with initial data $\rho_{0}=U_{m}(\tau,\cdot)$ with $\tau>0$ and $U_{m}$ a suitable Barrenblatt solution)  the uniqueness consists in proving a weak-strong uniqueness theorem since we can assume that the solution $\big(U_{m}(T+\tau,\cdot),-\n\va(U_{m}(T+\tau,\cdot))\big)$ is strong.
\end{remarka}
\begin{remarka}
The second important remark consists to point out  the fact that the question of global weak solution for the system (\ref{0.2}) remains open, indeed Mellet and Vasseur have proved the stability of the global weak solution in \cite{MV}. However it seems not so easy to construct a regular sequel $(\rho_{n},u_{n})$ approximating (\ref{0.2}) (typically by a Friedrich process) and verifying uniformly all the entropies of \cite{MV} (\ref{21a}), (\ref{22a}) and (\ref{gaina}).
\end{remarka}
\begin{corollaire}
Let $\gamma>1$, $\nu_{1}$ and $\nu_{2}$  with the hypothesis of theorem \ref{theo3}. Let $\rho_{0}$ and $u_{0}$ verifying the hypothesis of the theorem \ref{theo3} with $u_{0}=-\n\va(\rho_{0})$ and $\mu(\rho)=\mu\rho^{\alpha}$, $\lambda(\rho)$ verifying also the hypothesis of theorem \ref{theo3}. We assume here that there exists a unique quasi solution of system (\ref{3systemev}) with such initial data, in particular the density of this quasi solution verifies (\ref{P}). Then as in theorem \ref{theo3}  $(\rho_{\ep},u_{\ep})$ converges in distribution sense to a global weak quasi-solution $(\rho,u)$ of the system (\ref{3systemev}) such that $\rho$ is solution of (\ref{P}). Moreover:
\begin{itemize}
\item If $\alpha>1$, $\rho_{0}\in L^{\infty}(\R^{N})$ and the support of $\rho_{0}$ is compact then we have for $T>0$ $\rho_{\e}$ converges strongly to $\rho$ in $C([0,T], L^{1}(\R^{N})\cap L^{p}(\R^{N}))$ with $1<p<1+\nu_{1}$ if $N=3$ and in $C([0,T], L^{p}(\R^{N}))$ for any $p>1$.\\
For all $\eta>0$ it exists $\ep_{0}>0$ such that for all $0<\ep\leq\ep_{0}$ $\forall t\in[0,T]$ we have:
$$\rho_{\ep}(t,\cdot)=\rho(t,\cdot)+f_{\ep}(t,\cdot),$$
with $\rho(t,\cdot)$ with compact support for any $t\in[0,T]$ (it exists $C>0$ such that ${\mbox supp}\, \rho(t,\cdot)\leq C t^{\beta}$ with $\beta=\frac{1}{N(\alpha-1)+2}$) and we have $\|f_{\e}(t,\cdot)\|_{L^{1}(\R^{N})}\leq\eta$ $\forall t\in[0,T]$.\\
\\
In particular for all $\eta>0$ it exists $\ep_{0}>0$ such that for all $0<\ep\leq\ep_{0}$ we have:
$$\rho_{\ep}(t,\cdot)=\rho(t,\cdot)+f_{\ep}(t,\cdot),$$
with $\rho(t,\cdot)$ verifying for $p$ as above:
$$\|\rho(t)\|_{L^{p}(\R^{N})}\leq Ct^{-\sigma_{p }}\|\rho_{0}\|^{\alpha_{p}}_{L^{1}(\R^{N})},$$
with $\sigma_{p}=\frac{N(\alpha-1)+2p}{(N(\alpha-1)+2)p}$ and $\alpha_{p}=\frac{N(p-1)}{(N(\alpha-1)+2)p}$,
and we have $\|f_{\e}(t,\cdot)\|_{L^{p}(\R^{N})}\leq\eta$ for all $t\in[0,T]$.
\item If $\alpha\geq 1-\frac{1}{N}$, for all $\eta>0$ for all compact $K$ it exists $T>0$, it exists $\ep_{0}>0$ such that for $0<\ep\leq \ep_{0}$ we have:
$$\|\rho_{\ep}(t,\cdot)-U_{m}(t,\cdot)\|_{L^{1}(K)}\leq\eta\;\;\;\forall t\in[T,2T].$$
\end{itemize}
\label{cor}
\end{corollaire}
\begin{remarka}
It is very surprising to observe that for $\e$ small enough $\rho_{\e}$ is subjected to a type of damping effect in $L^{p}(\R^{N})$ with $p$ choose as above up to a small remainder term in $L^{p}(\R^{N})$. Let us point out that this effect is similar to the dispersion property for the Schro\"dingier or the waves equations. In \cite{Da,arma} we observe a damping effect on the density due to the role of the pressure, but in our case the pressure tends to disappear. As for the porous media equation this effect seems purely non linear and is exhibited because the particular structure of the viscosity coefficients.
\end{remarka}
\begin{remarka}
Under the hypothesis of uniqueness of the quasi solution when $u_{0}=-\n\va(\rho_{0})$ we show that the solution of the highly compressible Navier-Stokes equation are not so far to have a finite speed of propagation when we take a initial density with compact support. Indeed this is the case modulo a perturbation $f_{\e}$ of small $L^{1}$ norm. Similarly modulo this hypothesis of uniqueness we expect a asymptotic convergence of $\rho_{\e}$ to the Barrenblatt solution of $L^{1}$ norm $\|\rho_{0}\|_{L^{1}(\R^{N})}$ modulo a small perturbation.
\end{remarka}
The paper is structured in the following way: in section \ref{porous} we recall some important results on the porous media and the fast diffusions equations. In section \ref{section3} we adapted the entropy of \cite{BD} and \cite{MV} to the case of the quasi-solutions. In section \ref{section4}, we give a few notation, some compactness results and briefly introduce the basic Fourier analysis techniques needed to prove our result. In section \ref{section4} we prove theorem \ref{theo1} and in section \ref{section5} we show the theorem \ref{theo2}. In section \ref{section6} we conclude with the proof of the theorem \ref{theo3} and the forollary \ref{cor}. An appendix is postponed in order to prove rigorously some technical lemmas.
\section{Important results on the porous media equation}
\label{porous}
for the sake of completeness for the reader which are not familiar with the porous and the fast diffusion equations, we are going to recall in this section some essential results on the porous media and the fast diffusion equations. The majority of them are directly issue from the
excellent book \cite{Vaz}, \cite{Vaz1} from V\'azquez. In this part in order to simplify the problem we shall only consider the following equation with $\alpha>1-\frac{1}{N}$:
\begin{equation}
\begin{cases}
\begin{aligned}
&\p_{t}\rho-2\mu\D\rho^{\alpha}=0,\\
&\rho(0,\cdot)=\rho_{0}\geq 0.
\end{aligned}
\end{cases}
\label{porprim}
\end{equation}
Let us start with the case where $\alpha>1$
\subsection{Porous Media, $\alpha>1$}
In the sequel we shall set $Q=(0,+\infty)\times\R^{N}$. Let us recall the notion of global strong solution for the equation (\ref{P}) of the porous medium equation ($\alpha>1$) %and of the fast diffusion equation ($0<\alpha<1$) 
(see \cite{Vaz} chapter 9 for more details).% and \cite{Vaz1}). 
\begin{definition}
We say that a function $\rho\in C([0,+\infty), L^{1}(\R^{N}))$ positive is a strong $L^{1}$ solution of problem (\ref{porprim})  if:
\begin{itemize}
\item $\rho^{\alpha}\in L^{1}_{loc}(0,+\infty, L^{1}(\R^{N}))$ and $\rho_{t}, \D \rho^{\alpha}\in L^{1}_{loc}((0,+\infty)\times\R^{N})$
\item $\rho_{t}=\mu\D\rho^{\alpha}$ in distribution sense.
\item $u(t)\rightarrow \rho_{0}$ as $t\rightarrow0$ in $L^{1}(\R^{N})$.
\end{itemize}
%A locally integrable function $\rho$ defined in $(0,T)\times\R^{N}$
\end{definition}
Let us mention (see \cite{Vaz} p197) that we have the following theorem of global strong solution, we are going to give a sketch of the proof of the existence and uniqueness of the $L^{1}$ solution (which is perfectly detailed in \cite{Vaz} Chapter 6 and 9). Indeed it will be important to understand this point for the proof of the last part of theorem \ref{theo2}. 
%\begin{theorem}
%Let $\alpha>0$. For every $\rho_{0}\in L^{1}(\R^{N})$ positive there exists a unique global strong solution $\rho$ positive of problem (\ref{P}) such that $\rho\in C([0,+\infty), L^{1}(\R^{N}))\cap L^{\infty}((\tau,+\infty)\times\R^{N})$ for every $\tau>0$. %The solution satisfies estimate \ref
%\end{theorem} 
%\begin{remarka}
\begin{theorem}
\label{theo2.4}
Let $\alpha>1$ For every non-negative function $\rho_{0}\in L^{1}(\R^{N})\cap L^{\infty}(\R^{N})$ there exists a unique global  strong solution $\rho\geq 0$ of (\ref{porprim}). Moreover, $\p_{t}\rho\in L^{p}_{loc}(Q)$ for $1\leq p<\frac{(\alpha+1)}{\alpha}$ and:
$$
\begin{aligned}
&\p_{t}\rho\geq-\frac{\rho}{(\alpha-1)t}\;\;\;\mbox{in}\;\;{\cal D}^{'}(Q),\\
&\|\p_{t}\rho(t,\cdot)\|_{L^{1}(\R^{N})}\leq \frac{2\|\rho_{0}\|_{L^{1}(\R^{N})}}{(\alpha-1)t}.
\end{aligned}
$$
Let $\rho_{1}$ and $\rho_{2}$ be two strong solutions of (\ref{porprim}) in $(0,T)\times\R^{N}$ then for every $0\leq \tau<t$ 
\begin{equation}
\|\big(\rho_{1}(t,\cdot)-\rho_{2}(t,\cdot)\big)_{+}\|_{L^{1}(\R^{N})}\leq \|\big(\rho_{1}(\tau,\cdot)-\rho_{2}(\tau,\cdot)\big)_{+}\|_{L^{1}(\R^{N})}.
\label{contraction}
\end{equation}
If  $\rho_{1}$ and $\rho_{2}$ are two strong solution with initial data $\rho_{01}$ and $\rho_{02}$ with $\rho_{01}\leq \rho_{02}$ in $\R^{N}$, then $\rho_{1}\leq\rho_{2}$ almost everywhere in $(0,+\infty)\times\R^{N}$.
\end{theorem}
{\bf Proof:} The proof is decomposed in three steps. Let us recall that in the sequel the initial density is always positive.
\subsubsection*{Energy global weak solution on a bounded domain $\Omega$}
Let us mention that a global weak solution in our setting is given by the definition 5.4 of \cite{Vaz} and verify on the initial data $\rho_{0}\in L^{1}(\Omega)\cap L^{1+\alpha}(\Omega)$. By considering $\rho_{0n}=\max(\rho_{0}+\frac{1}{n},n)$  which is strictly positive with boundary condition such that $\rho_{n}(t,\cdot)=\frac{1}{n}$ on the boundary of $\Omega$, by the theory of the quasilinear equations (see Ladyzhenskaya et al \cite{La} or Friedman in \cite{F}) it exists global classical solutions which verify the energy estimates:
$$\int^{T}_{0}\int_{\Omega}|\n\rho_{n}^{\alpha}|^{2} dxdt+\frac{1}{\alpha+1}\int_{\Omega}\rho_{n}^{\alpha+1}dxdt\leq \frac{1}{\alpha+1}\int_{\Omega}\rho_{0n}^{\alpha+1}.$$
Furthermore by applying the maximum principle we know that $(\rho_{n})_{n\in\mathbb{N}}$ is a decreasing sequel so to converges everywhere to a limit $\rho$. Since $\rho_{n}$ is uniform bounded in $L^{\infty}_{T}(L^{1+\alpha}(\Omega))$ for any $T>0$, it implies that up to a subsequence $\rho_{n}$ converge weakly to $\rho$ in $L^{\infty}_{T}(L^{1+\alpha}(\Omega))$ . Furthermore $\n \rho_{n}^{\alpha}$ is uniformly bounded in $L^{2}_{T}(L^{2}(\Omega))$ for any $T>0$ it implies that up to a subsequence $\n\rho_{n}^{\alpha}$ converges to a limit $\psi$.\\
By applying the lemma \ref{lemmeimp} to $\rho_{n}^{\alpha}$ and the fact that $\rho_{n}^{\alpha}$ is uniformly bounded in $L^{\infty}_{T}(L^{1+\e}(\Omega))$ with $\e>0$ and converges almost everywhere to $\rho^{\alpha}$ it shows that $\rho_{n}^{\alpha}$ converges strongly to $\rho^{\alpha}$ in $L^{1}_{loc}((0,T)\times\Omega)$. In particular we deduce that $\psi=\n\rho^{\alpha}$.
\subsubsection*{$L^{1}$  global weak solution on a bounded domain $\Omega$}
By the fundamental $L^{1}$ contraction principle which ensures that for any solution $\rho_{1}$, $\rho_{2}$ of the porous media in $L^{1}$ we have:
\begin{equation}
\|\rho_{1}(t)-\rho_{2}(t)\|_{L^{1}(\Omega)}\leq \|\rho_{01}-\rho_{02}\|_{L^{1}(\Omega)}.
\label{maxcru}
\end{equation}
The $L^{1}$ limit solution $\rho$ of the equation (\ref{P}) consists in considering the $L^{1}$ limit of an approximate energy solution $\rho_{n}$ of (\ref{P}) with $\rho_{0n}$ converging to $\rho_{0}$ in $L^{1}(\Omega)$. By the estimate (\ref{maxcru}) we check easily that the limit $\rho$ does not depend on the choice of the regularizing sequel $\rho_{0n}$. In addition we can verify that $\rho$ is in $C([0,+\infty),L^{1}(\Omega))$ (see p129 in \cite{Vaz}). However an important question remains, is the limit solution $\rho$ is a weak solution according the definition 5.2 of \cite{Vaz}. The answer is positive. Indeed for any approximative solution we can check that:
\begin{equation}
\begin{aligned}
&\int_{\Omega}|\rho_{n}-\rho_{m}|(t,x)\xi(x)dx+\int^{t}_{0}\int_{\Omega}\|\rho_{n}^{\alpha}(s,x)-\rho_{m}^{\alpha}(s,x)|dxds\leq \int_{\Omega}|\rho_{0n}-\rho_{0m}|(x)\xi(x)dx,
\end{aligned}
\label{6.5}
\end{equation}
where $\xi$ is the unique solution of the problem:
$$\D\xi=-1\;\;\mbox{in}\;\Omega,\;\;\xi=0\;\;\mbox{on}\;\;\p\Omega.$$
It implies that $\rho_{n}^{\alpha}$ is a Cauchy sequence in $L^{1}((0,T)\times\Omega)$ for any $T>0$ and then converge strongly to $\psi$ but $\rho_{n}$ converge strongly in $L^{\infty}_{T}(L^{1}(\Omega)$ for any $T>0$ then up to a subsequence $rho_{n}$ converges almost everywhere to $\rho$ and $\rho_{n}^{\alpha}$ up to a subsequence converges to $\rho^{\alpha}$ and $\psi$ then $\psi=\rho^{\alpha}$ which implies that $\rho$ is a very weak solution of (\ref{P}) in the sense of the definition 5.2 of \cite{Vaz}.
\subsubsection*{$L^{1}$ solution on $\R^{N}$}
We start by approximating the initial data by setting:
$$\rho_{0n}(x)=\max(-n,\min(\rho_{0}(x),n)\chi (nx),$$
whit $\chi\in C^{\infty}_{0}(\R^{N})$ such that the support of $\chi$ is embedded in the ball $B(0,2)$ and $\chi=1$ on $B(0,1)$. We have seen that it exists global $L^{1}$ solutions for the homogeneous Cauchy-Dirichlet problem in bounded domain $\Omega_{n}=B(0,2n)$. It suffices to pass to the limit when $n$ goes to infinity by using the same type of compactness argument.\\
Let us mention to finish that the uniqueness of the $L^{1}$ solution is a direct consequence of the $L^{1}$ contraction principle. \hfill{$\Box$}
\begin{remarka}
Let us recall that there exists global weak solution which are not classical it means not $C^{\infty}$ even if the initial data is $C^{\infty}$ (see a example due to Aronson in the problem 5.7 of \cite{Vaz}). 
\end{remarka}
We are now to recall the so called $L^{1}-L^{\infty}$ smoothing effect (as for the dispersive equations), we refer to \cite{Vaz} p 202.
\begin{theorem}
For every $t>0$ we have:
$$\rho(t,x)\leq C\|\rho_{0}\|_{L^{1}(\R^{N})}^{\sigma} t^{-\beta},$$
with $\sigma=\frac{2}{N(\alpha-1)+2}$, $\beta=\frac{N}{N(\alpha-1)+2}$ and $C>0$ depends only on $\alpha$ and $N$. The exponents are sharp.
\end{theorem}
Let us finnish by giving a more general theorem of existence of global strong solution for (\ref{porprim}) with some properties on the solutions (see \cite{Vaz} p 204-205).
\begin{theorem}
\label{theo2.6}
For every $\rho_{0}\in L^{1}(\R^{N})$ there exists a unique global strong solution of (\ref{porprim}) such that $\rho\in C([0,+\infty),L^{1}(\R^{N}))\cap L^{\infty}((\tau,+\infty)\times\R^{N})$ with $\tau>0$. Furthermore we have the following $L^{\infty}$ estimate:
$$|\rho(t,x)|\leq C\|\rho_{0}\|_{L^{1}(\R^{N})}^{\sigma} t^{-\beta},$$
with $\sigma=\frac{2}{N(\alpha-1)+2}$, $\beta=\frac{N}{N(\alpha-1)+2}$ and $C>0$ depends only on $\alpha$ and $N$.
Moreover we have the following properties:
\begin{enumerate}
\item The solutions are continuous functions of $(t,x)$ in $Q$ with a uniform modulus of continuity for $t\geq \tau>0$.
\item The maximum principles holds.
\item if $\rho_{0}$ is strictly positive and continuous, then $\rho\in C^{\infty}(Q)\cap C(\bar{Q})$ and is a classical solution of (\ref{porprim}).
\item For every $p\in(1,+\infty)$ we have the following regularizing effect, $\rho(t,\cdot)$ belongs in $L^{p}(\R^{N})$ and:
$$\|\rho(t)\|_{L^{p}(\R^{N})}\leq Ct^{-\sigma_{p }}\|\rho_{0}\|^{\alpha_{p}}_{L^{1}(\R^{N})},$$
with $\sigma_{p}=\frac{N(\alpha-1)+2p}{(N(\alpha-1)+2)p}$ and $\alpha_{p}=\frac{N(p-1)}{(N(\alpha-1)+2)p}$.
\end{enumerate}
\end{theorem} 
Let us conclude this section with two important theorem showing the finite speed of propagation for the porous media equation (see \cite{Vaz} p 210) and the time asymptotic behavior of the solution which converges to Barrenblatt solutions (see \cite{Vaz2} p 69).
\begin{proposition}
Let $\rho$ be the global strong solution of (\ref{porprim}) with initial data $\rho_{0}\in L^{1}(\R^{N})\cap L^{\infty}(\R^{N})$ and assume that $\rho_{0}$ has a compact support then for every $t>0$ the support of $\rho(t,\cdot)$ is a bounded set.
\label{speed}
\end{proposition}
\begin{theorem}
\label{theo2.7}
Let $\rho(t,x)$ be the unique global strong solution of (\ref{porprim}) with initial data $\rho_{0}\in L^{1}(\R^{N})$, $\rho_{0}\geq 0$. Let $U_{m}$ be the Barrenblatt with the same mass as $\rho_{0}$. Then we have:
\begin{equation}
\lim_{t\rightarrow+\infty}\|\rho(t)-F_{m}(t)\|_{L^{1}(\R^{N})}=0.
\label{convasym}
\end{equation}
Convergence holds also in $L^{\infty}$ norm:
\begin{equation}
\lim_{t\rightarrow+\infty}t^{\beta}\|\rho(t)-F_{m}(t)\|_{L^{\infty}(\R^{N})}=0,
\label{convasym1}
\end{equation}
with $\beta=\frac{N}{N(\alpha-1)+2}$.
\end{theorem}
\begin{remarka}
For more results in this direction we refer also to \cite{Jean}.
\end{remarka}
Let us conclude this section by giving general results (essentially extracted form the chapter 9 from \cite{Vaz}) on porous media equation of the form:
\begin{equation}
\begin{cases}
\begin{aligned}
&\p_{t}\rho-2\D\mu(\rho)=0,\\
&\rho(0,\cdot)=\rho_{0}.
\end{aligned}
\end{cases}
\end{equation}
\subsubsection{General viscosity coefficients}
\label{pormu}
We assume here that $\mu$ verifies the following assumptions:
\begin{itemize}
\item $\mu$ is a continuous and increasing function: $\R\rightarrow \R$ with $\mu(0)=0$
\item $\mu$ has at least linear growth at infinity in the sense that it exists $c>$ such that for large $|s|$ we have:
$$|\mu(s)|\geq c |s|>0.$$
\end{itemize}
\begin{definition}
A locally integrable function $\rho$ defined in $Q_{T}$ is said to be a weak solution of the problem if:
\begin{enumerate}
\item $\mu(\rho)\in L^{2}(0,T,H^{1}(\R^{N}))$
\item $\rho$ satisfies the identity:
\begin{equation}
\int\int_{Q_{T}}(\n\mu(\rho)\cdot\n \va-\rho\p_{t}\va)dxdt=\int_{\R^{N}}\rho_{0}(x)\va(0,x)dx,
\label{9.45}
\end{equation}
\end{enumerate}
for any function $\va\in C^{1}(\bar{Q}_{T})$ which vanishes for $t=T$ and has uniformly bounded support in the space variable.
\end{definition}
We define $L_{\mu}(\R^{N})$ by the set of measurable function $\rho_{0}$ such that $\mu(\rho_{0})\in L^{1}(\R^{N})$. We shall consider $\psi$ the primitive of $\mu$:
$$\psi(s)=\int^{s}_{0}\mu(\tau)d\tau.$$
Let $X_{T}=L_{\mu}(\R^{N})\cap L^{1}(\R^{N})$ and $Y=L^{\infty}(Q_{T})\cap L^{1}(Q_{T})$ for $0<T\leq+\infty$.
\begin{theorem}
Let $\rho_{0}\in X$. Then it exists a unique global weak solution defined in $(0,+\infty)$ and $\n\mu(\rho)\in L^{2}(Q_{T})$. Moreover we also have $\rho\in L^{\infty}((0,T),X)$.
\end{theorem}
\subsection{Fast diffusion equations, $0<\alpha<1$}
The situation is different in the case of fast diffusion equation $0<\alpha<1$, indeed in this case infinite propagation occurs and solution may even vanish in finite time. Let us mention that when $\alpha$ is in the interval $(0,m_{c})$ with $m_{c}=\max(0,\frac{N-2}{N})$ then it can appears a phenomena of extinction of the solution in finite time. In particular it implies a lost of the initial mass when $\rho_{0}$ is in $L^{1}$ (it implies also a lost of the regularity of the solution). We refer to \cite{Vaz1} theorem 5.7 for a necessary condition of extinction, in particular the solution belongs in an appropriate Marcinkewitz space $M_{p^{*}}(\R^{N})$.\\
Let us mention that in the case $\alpha\in(m_{c},1)$ the situation is quite similar to the case of the porous media equation (except the infinite propagation speed) as the mass is preserved which implies no extinction in finite time. Moreover we have self similar solutions also discovered by Barrenblatt that we can write under the following form:
$$U_{m}(t,x)=t^{-\gamma_{1}}F(x t^{-\beta})\;\;\mbox{with}\;\;F(x)=(C+\kappa_{1}|x|^{2})^{\frac{-1}{\alpha-1}}_{+},$$
with $\kappa_{1}=\frac{(1-\alpha)\gamma_{1}}{2N\alpha}$. In particular the proof of the most of the previous result in the last section are based on the existence of Barrenblatt solutions and on the maximum principle or in other word the $L^{1}$ contraction principle. This two fundamental point arise also in the case $\alpha\in(m_{c},1)$ which implies that we have the most of the result of the case $\alpha>1$ exists also in this case (using essentially the same proof).
In particular similarly to the case $\alpha>1$ in the situation $\alpha_{c}<\alpha<1$ the global strong solution converges asymptotically to a Barrenblatt solution and we have also regularizing effect $L^{1}-L^{\infty}$. For more details in this situation we refer the reader to the excellent books of V\'azquez \cite{Vaz,Vaz1}.
\section{Entropy inequality for the quasi-solutions and basic tools}
\label{section3}
\subsection{Entropy for the quasi-solution of the system (\ref{3systeme}) }
We now want to establish new entropy inequalities for system (\ref{3systeme}) by applying the entropy inequalities discovered in \cite{BD,MV}. More precisely if we assume that $(\rho,u)$ are classical solutions of system (\ref{3systeme}), we obtain the following proposition.
\begin{proposition}
%\frac{1}{\gamma-1}\rho^{\gamma}
Assume that $(\rho,u)$ are classical solutions of system (\ref{3systeme}) then for all $t>0$ we have the two following entropy:
\begin{equation}
\begin{aligned}
&\int_{\R^{N}}\frac{1}{2}\rho|u|^{2}(t,x)\,dx+\int^{t}_{0}\int_{\Om}2\mu(\rho)|D u|^{2}dxdt+\int^{t}_{0}\int_{\Om}\lambda(\rho)|{\rm div}u|^{2}dxdt \\
&\hspace{9cm}=\int_{\R^{N}}\rho_{0}|u_{0}|^{2}(x)\,dx.
\end{aligned}
\label{21}
\end{equation}
\begin{equation}
\begin{aligned}
&\int_{\R^{N}}\frac{1}{2}\rho|u+\n\va(\rho)|^{2}(t,x)\,dx+\frac{1}{2}\int^{t}_{0}\int_{\R^{N}}\mu(\rho)|\n u-^{t}\n u|^{2}dxdt\\
&\hspace{7cm}=\int_{\R^{N}}\frac{1}{2}\rho_{0}|u_{0}+\n\va(\rho_{0})|^{2}(x)\,dx.
\end{aligned}
\label{22}
\end{equation}
\label{propBD}
\end{proposition}
{\bf Proof:}
In order to obtain (\ref{21}) it suffices to multiplying the momentum equation by $u$ and integrating over $(0,t)\times\R^{N}$.\\
Let us briefly recall the proof of the second entropy (\ref{22}) introduced by Bresch and Desjardins in \cite{BD}. To this purpose, we have to study:
$$\frac{d}{dt}\int\big[\frac{1}{2}\rho|u|^{2}+\rho u\cdot\n\va(\rho)+\frac{1}{2}\rho|\n\va(\rho)|^{2}\big]dx.$$
\subsubsection*{Step 1:}
First by the mass equation, we have:
$$\begin{aligned}
&\frac{1}{2}\int\frac{d}{dt}(\rho|\n\va(\rho)|^{2})dx=\int\rho\frac{d}{dt}\frac{|\n\va(\rho)|^{2}}{2}dx-\int\frac{|\n\va(\rho)|^{2}}{2}{\rm div}(\rho u)dx,\\
&\hspace{0,5cm}=-\int\rho\n u\n\va(\rho)\otimes\n \va(\rho)dx+\int\rho^{2}\va^{'}(\rho)\D\va(\rho){\rm div}udx+\int\rho|\n\va(\rho)|^{2}{\rm div}udx.
\end{aligned}
$$
\subsubsection*{Step 2}
It remains to estimate the derivative of the cross product:
$$
\begin{aligned}
\frac{d}{dt}\int\rho u\cdot\n\va(\rho)dx&=\int\n\va(\rho)\cdot\frac{d}{dt}(\rho u)+
\int\rho u\cdot\frac{d}{dt}\n\va(\rho)dx\\
&=\int\n\va(\rho)\cdot\frac{d}{dt}(\rho u)-\int{\rm div}(\rho u)\va^{'}(\rho)\frac{d}{dt}\rho dx\\
&=\int\n\va(\rho)\cdot\frac{d}{dt}(\rho u)+\int{\rm div}(\rho u)^{2}\va^{'}(\rho)dx.
\end{aligned}
$$
Multiplying the momentum equation by $\n\va(\rho)$, we get:
$$
\begin{aligned}
&\int\n\va(\rho)\cdot\frac{d}{dt}(\rho u)=-\int(2\mu(\rho)+\lambda(\rho))\D\va(\rho){\rm div}udx+2\int\n u:\n\va(\rho)\otimes\n \mu(\rho)dx\\
&\hspace{4cm}-2\int\n\va(\rho)\cdot\n \mu(\rho){\rm div}u dx-\int\n\va(\rho){\rm div}(\rho u\otimes u)dx,
\end{aligned}
$$
where we use the fact that:
$$\int\n(\lambda(\rho){\rm div}u)\cdot\n\va(\rho)dx=-\int \lambda(\rho)\D\va(\rho){\rm div}udx,$$
and:
$$
\begin{aligned}
&\int{\rm div}(2 \mu(\rho)D(u))\cdot\n\va(\rho)dx=\int\p_{j}(\mu(\rho)\p_{j}u_{i})\p_{i}\va(\rho)dx+\int\p_{j}(\mu(\rho)\p_{i}u_{j})\p_{i}\va(\rho)dx,\\
&\hspace{5cm}=\int\p_{i}(\mu(\rho)\p_{j}u_{i})\p_{j}\va(\rho)dx+\int\p_{j}(\mu(\rho)\p_{i}u_{j})\p_{i}\va(\rho)dx,\\
&\hspace{5cm}=\int\p_{i}\mu(\rho)\p_{j}u_{i})\p_{j}\va(\rho)dx-\int\p_{j}\mu(\rho)\p_{i}u_{i}\p_{j}\va(\rho)dx\\
&\hspace{1,5cm}=-\int \mu(\rho)\p_{i}u_{i}\p_{jj}\va(\rho)dx+\int\p_{j}\mu(\rho)\p_{i}u_{j}\p_{i}\va(\rho)dx-\int\p_{i}\mu(\rho)\p_{j}u_{j}\p_{i}\va(\rho)dx\\
&\hspace{10cm}-\int \mu(\rho)\p_{j}u_{j}\p_{ii}\va(\rho)dx\\
&=2\int\n u:\n \mu(\rho)\otimes\n \va(\rho)dx-2\int\n \mu(\rho)\cdot\n\va(\rho){\rm div}u dx-2\int \mu(\rho)\D\va(\rho){\rm div}udx.
\end{aligned}
$$
\subsubsection*{Step 4}
Since$\va$, $\mu$ and $\lambda$ satisfies (\ref{BD}) and (\ref{11}), then we obtain:
$$
\begin{aligned}
&\frac{d}{dt}\big(\int\rho u\cdot\n\va(\rho)+\rho\frac{|\n\va(\rho)|^{2}}{2}dx\big)=-\int\n\va(\rho){\rm div}(\rho u\otimes\rho u)dx+\int\va^{'}(\rho)({\rm div}(\rho u))^{2}dx.
\end{aligned}
$$
Finally we have:
$$
\begin{aligned}
&-\int\n\va(\rho){\rm div}(\rho u\otimes\rho u)dx+\int\va^{'}(\rho)({\rm div}(\rho u))^{2}dx\\
&\hspace{1cm}=-\int\va^{'}(\rho)u\cdot\n\rho{\rm div}(\rho u)-\va^{'}(\rho)\n\rho(\rho u\cdot\n u)+\va^{'}(\rho)({\rm div}(\rho u))^{2}dx\\
&\hspace{1cm}=\int\rho\va^{'}(\rho){\rm div}u\,{\rm div}(\rho u)-\rho\va^{'}(\rho)\n\rho(u\cdot\n u)dx\\
&\hspace{1cm}=\int\rho^{2}\va^{'}(\rho)({\rm div}u)^{2}+\rho\va^{'}(\rho)u\cdot\n\rho{\rm div}u-\rho\va^{'}(\rho)\n\rho(u\cdot\n u)dx,
\end{aligned}
$$
then by (\ref{BD}) and (\ref{11}), we get:
$$
\begin{aligned}
&-\int\n\va(\rho){\rm div}(\rho u\otimes\rho u)dx+\int\va^{'}(\rho)({\rm div}(\rho u))^{2}dx,\\
&\hspace{4cm}=2\int\rho h^{'}(\rho)({\rm div}u)^{2}+\n h(\rho)\cdot u{\rm div}u-\n(h(\rho))(u\cdot\n u)dx,\\
&=2\int\rho \mu^{'}(\rho)({\rm div}u)^{2}-\mu(\rho)({\rm div}u)^{2}-\mu(\rho)u\cdot\n{\rm div}udx+2\int \mu(\rho)
\p_{i}u_{j}\p_{j}u_{i}\\
&\hspace{10cm}+\mu(\rho)u\cdot\n{\rm div}u dx,\\
&\hspace{5cm}=\int(2\rho \mu^{'}(\rho)-2\mu(\rho))({\rm div}u)^{2}+2\mu(\rho)\p_{i}u_{j}\p_{j}u_{i}dx\\
&\hspace{5cm}=\int \lambda(\rho)({\rm div}u)^{2}+\int 2\mu(\rho)\p_{i}u_{j}\p_{j}u_{i}dx,
\end{aligned}
$$
whichgives:
$$
\begin{aligned}
&\frac{d}{dt}\big(\int\rho u\cdot\n\va(\rho)+\rho\frac{|\n\va(\rho)|^{2}}{2}dx\big)\\
&\hspace{4cm}=\int \lambda(\rho)({\rm div}u)^{2}dx+\int2\mu(\rho)\p_{i}u_{j}\p_{j}u_{i}dx,\\
\end{aligned}
$$
Adding this equality and (\ref{21}), and using the fact that:
$$\int2\mu(\rho)|(u)|^{2}-\int2\mu(\rho)\p_{i}u_{j}\p_{j}u_{i}dx=\int2\mu(\rho)(\frac{\p_{i}u_{j}-\p_{j}u_{i}}{2})^{2},$$
we easily get (\ref{22}).\hfill{$\Box$}
As in \cite{MV} we are also interested in getting a gain of integrability on the velocity. We have then the following proposition.
\begin{proposition}
Assume that:
$$2\mu(\rho)+N\lambda(\rho)\geq\nu\lambda(\rho)$$
for some $\nu\in (0,1)$ (which is a part of (\ref{11})). Then it exists $C>0$ such that smooth solutions of (\ref{3systeme})
satisfy the following inequality:
\begin{equation}
\begin{aligned}
&\frac{d}{dt}\int\rho\frac{1+|u|^{2}}{2}{\rm ln}(1+|u|^{2})dx+\frac{\nu}{2}\int \mu(\rho)[1+{\rm ln}(1+|u|^{2})|D u|^{2}dx\leq C\int \mu(\rho)|\n u|^{2}dx.
\end{aligned}
\label{gain}
\end{equation}
for any $\de\in(0,2)$, and with $|\n u|^{2}=\sum_{i}\sum_{j}|\p_{i}u_{j}|^{2}$.
\end{proposition}
{\bf Proof:} Multiplying the momentum equation by $\big(1+{\rm ln}(1+|u|^{2})\big)u$, we get:
$$
\begin{aligned}
&\int\rho\frac{d}{dt}\big[\frac{1+|u|^{2}}{2}{\rm ln}(1+|u|^{2})\big]dx+\int\rho u\cdot\n\big(\frac{1+|u|^{2}}{2}{\rm ln}(1+|u|^{2})\big)dx\\
&\hspace{0,5cm}+\int 2\mu(\rho)(1+{\rm ln}(1+|u|^{2})|D(u)|^{2}dx+\int 2\mu(\rho)\frac{2u_{i}u_{k}}{1+|u|^{2}}\p_{j}u_{k}D_{ij}(u)dx\\
&\hspace{1cm}+\int \lambda(\rho)(1+{\rm ln}(1+|u|^{2}){\rm div}u|^{2}dx+\int \lambda(\rho)\frac{2u_{i}u_{k}}{1+|u|^{2}}\p_{i}u_{k}{\rm div}udx=0.\\
\end{aligned}
$$
Since:
$$|{\rm div}u|^{2}\leq N|\n u|^{2}\;\;\mbox{and}\;\;\nu\mu(\rho)\leq2\mu(\rho)+N\lambda(\rho),$$
we obtain:
\begin{equation}
\begin{aligned}
&\int\rho\frac{d}{dt}\big[\frac{1+|u|^{2}}{2}{\rm ln}(1+|u|^{2})\big]dx+\int\rho u\cdot\n\big(\frac{1+|u|^{2}}{2}{\rm ln}(1+|u|^{2})\big)dx\\
&\hspace{3cm}+\nu\int \mu(\rho)(1+{\rm ln}(1+|u|^{2})|D(u)|^{2}dx\leq C\int \mu(\rho)|\n u|^{2}dx.\\
\end{aligned}
\label{entropcle}
\end{equation}
Moreover multiplying the mass equation by $\frac{1+|u|^{2}}{2}{\rm ln}(1+|u|^{2})$ and integrating by parts, we have:
$$\int\frac{1+|u|^{2}}{2}{\rm ln}(1+|u|^{2})\frac{d}{dt}\rho dx-\int\rho u\cdot\n\big(\frac{1+|u|^{2}}{2}{\rm ln}(1+|u|^{2})\big) dx=0$$
We deduce that:
$$
\begin{aligned}
&\frac{d}{dt}\int\rho\frac{1+|u|^{2}}{2}{\rm ln}(1+|u|^{2})dx+\frac{\nu}{2}\int \mu(\rho)[1+{\rm ln}(1+|u|^{2})]|D(u)|^{2}dx\\
&\hspace{3cm}\leq C\int \mu(\rho)|\n u|^{2}dx.
\end{aligned}
$$
It concludes the proof. \hfill{$\Box$}
\subsection{Entropy for the compressible Navier-Stokes equations}
We are going now to consider the following system with $\mu(\rho)$ and $\lambda(\rho)$ verifying (\ref{BD}):
\begin{equation}
\begin{cases}
\begin{aligned}
&\p_{t}\rho+{\rm div}(\rho u)=0,\\
&\p_{t}(\rho u)+{\rm div}(\rho u\otimes u)-{\rm div}(2\mu(\rho) D u)-\n (\lambda(\rho)\n u)+\e\n P(\rho)=0,\\
&(\rho,u)(0,\cdot)=(\rho_{0},u_{0}).
\end{aligned}
\end{cases}
\label{syscompress}
\end{equation}
Here $P$ corresponds to the pressure and we shall consider a $\gamma$ law $P(\rho)=a\rho^{\gamma}$ with $\gamma>1$ and $a>0$.
\begin{proposition}
%\frac{1}{\gamma-1}\rho^{\gamma}
Assume that $(\rho,u)$ are classical solutions of system (\ref{syscompress}) then for all $t>0$ we have the two following entropy:
\begin{equation}
\begin{aligned}
&\int_{\R^{N}}\big[\frac{1}{2}\rho|u|^{2}(t,x)+\frac{\e a}{\gamma-1}\rho^{\gamma}(t,x)\big]\,dx+\int^{t}_{0}\int_{\Om}2\mu(\rho)|D u|^{2}dxdt\\
&\hspace{2cm}+\int^{t}_{0}\int_{\Om}\lambda(\rho)|{\rm div}u|^{2}dxdt=\int_{\R^{N}}\big[\rho_{0}|u_{0}|^{2}(x)+\frac{\e a}{\gamma-1}\rho_{0}^{\gamma}(x)\big]\,dx.
\end{aligned}
\label{21a}
\end{equation}
\begin{equation}
\begin{aligned}
&\int_{\R^{N}}\big[\frac{1}{2}\rho|u+\n\va(\rho)|^{2}(t,x)+\frac{\e a}{\gamma-1}\rho^{\gamma}(t,x)\big]\,dx+\frac{1}{2}\int^{t}_{0}\int_{\R^{N}}\mu(\rho)|\n u-^{t}\n u|^{2}dxdt\\
&+\e a\int^{t}_{0}\int_{\R^{N}}\n\rho^{\gamma}\cdot\n\va(\rho)(s,x)ds\,dx=\int_{\R^{N}}\big[\frac{1}{2}\rho_{0}|u_{0}+\n\va(\rho_{0})|^{2}(x)\\
&\hspace{9cm}+\frac{\e a}{\gamma-1}\rho_{0}^{\gamma}(t,x)\big]\,dx.
\end{aligned}
\label{22a}
\end{equation}
\label{propBDa}
\end{proposition}
{\bf Proof:} We refer to \cite{BD} for the proof.
 \hfill{$\Box$}
 \begin{proposition}
Assume that:
$$2\mu(\rho)+N\lambda (\rho)\geq\nu_{1}\lambda(\rho)$$
for some $\nu\in (0,1)$ (which is a part of (\ref{11})). Then it exists $C>0$ such that smooth solutions of (\ref{3systeme})
satisfy the following inequality:
\begin{equation}
\begin{aligned}
&\frac{d}{dt}\int\rho\frac{1+|u|^{2}}{2}{\rm ln}(1+|u|^{2})dx+\frac{\nu}{2}\int \mu(\rho)[1+{\rm ln}(1+|u|^{2})|D u|^{2}dx\leq \\
&\hspace{1cm}C\int \mu(\rho)|\n u|^{2}dx+C\e^{2}\biggl(\int\big(\frac{\rho^{2\gamma-\frac{\delta}{2}}}{\mu(\rho)}\big)^{\frac{2}{2-\de}}dx\biggl)^{\frac{2}{2-\de}}\big(\int(\rho|u|^{2}+\rho)dx\big)^{\frac{\de}{2}}.
\end{aligned}
\label{gaina}
\end{equation}
for any $\de\in(0,2)$.
\end{proposition}
{\bf Proof:} The proof follows the same lines than the lemma 3.2 in \cite{MV}, for the sake of completeness let us adapt this proof to our situation.  Multiplying the momentum equation by $\big(1+{\rm ln}(1+|u|^{2})\big)u$, we get as in (\ref{entropcle}):
\begin{equation}
\begin{aligned}
&\int\rho\frac{d}{dt}\big[\frac{1+|u|^{2}}{2}{\rm ln}(1+|u|^{2})\big]dx+\int\rho u\cdot\n\big(\frac{1+|u|^{2}}{2}{\rm ln}(1+|u|^{2})\big)dx\\
&+\nu_{1}\int \mu(\rho)(1+{\rm ln}(1+|u|^{2})|D(u)|^{2}dx\leq -a\e\int  [1+{\rm ln}(1+|u|^{2})]u\cdot\n \rho^{\gamma}dx\\
&\hspace{9cm}+C\int \mu(\rho)|\n u|^{2}dx.\\
\end{aligned}
\label{entropcle1}
\end{equation}
%FAIRE PROPREMENT
It remains to bound the right hand side. We have:
$$
\begin{aligned}
&|\e \int  [1+{\rm ln}(1+|u|^{2})]u\cdot\n \rho^{\gamma}dx|\\
&\leq \e|\int\frac{2 u_{i}u_{k}}{1+|u|^{2}}\p_{i}u_{k}\rho^{\gamma}dx|+\e|\int [1+{\rm ln}(1+|u|^{2})]{\rm div}u\rho^{\gamma}dx|,\\
&\leq 2\e (\int \mu(\rho)|\n u|^{2}dx)^{\frac{1}{2}}(\int\frac{\rho^{2\gamma}}{\mu(\rho)}dx)^{\frac{1}{2}}+\e|\int [1+{\rm ln}(1+|u|^{2})]{\rm div}u\rho^{\gamma}dx|.
\end{aligned}
$$
Let us deal with the last term on the right hand side:
$$
\begin{aligned}
&\e|\int [1+{\rm ln}(1+|u|^{2})]{\rm div}u\rho^{\gamma}dx|\leq\\
&\leq \e(\int [1+{\rm ln}(1+|u|^{2})]\mu(\rho)({\rm div}u)^{2}dx)^{\frac{1}{2}}(\int [1+{\rm ln}(1+|u|^{2})]\frac{\rho^{2\gamma}}{\mu(\rho)}dx)^{\frac{1}{2}},\\
%&\hspace{8cm}+\frac{\nu_{1}}{2}(\int [1+{\rm ln}(1+|u|^{2})]\mu(\rho)({\rm div}u)^{2}dx)\\
&\hspace{8cm}+\frac{\e^{2}}{2\nu_{1}}(\int [1+{\rm ln}(1+|u|^{2})]\frac{\rho^{2\gamma}}{\mu(\rho)}dx).
\end{aligned}
$$
We deduce that it exists $C>0$ such that:
$$
\begin{aligned}
&|\e \int  [1+{\rm ln}(1+|u|^{2})]u\cdot\n \rho^{\gamma}dx|\leq \int\mu(\rho)|\n u|^{2}dx+
\frac{\nu_{1}}{2}(\int [1+{\rm ln}(1+|u|^{2})]\mu(\rho)({\rm div}u)^{2}dx)\\
&\hspace{8cm}+\frac{C\e^{2}}{2\nu_{1}}(\int [2+{\rm ln}(1+|u|^{2})]\frac{\rho^{2\gamma}}{\mu(\rho)}dx)
\end{aligned}
$$
where the last term satisfies (if $\de\in(0,2)$) for a $C>0$:
$$
\begin{aligned}
\e^{2}\int\frac{\rho^{2\gamma}}{h(\rho)}|u|^{\de}dx&\leq \e^{2}\biggl(\int\big(\frac{\rho^{2\gamma-\frac{\delta}{2}}}{\mu(\rho)}\big)^{\frac{2}{2-\de}}dx\biggl)^{\frac{2}{2-\de}}
\big(\int\rho[2+\ln(1+|u|^{2})]^{\frac{2}{\delta}}dx)^{\frac{\delta}{2}},\\
&\leq C\e^{2}\biggl(\int\big(\frac{\rho^{2\gamma-\frac{\delta}{2}}}{\mu(\rho)}\big)^{\frac{2}{2-\de}}dx\biggl)^{\frac{2}{2-\de}}\big(\int(\rho|u|^{2}+\rho)dx\big)^{\frac{\de}{2}}
\end{aligned}
$$
and the proposition follows.  \hfill{$\Box$}
%\begin{remarka}
%Let us mention that from proposition \ref{BD}, we can easily show by interpolation that $\rho^{\gamma}$ is bounded in $L^{\frac{5}{3}}((0,T)\times \Om)$ for $N=3$ (we refer to the lemma \ref{pression} for more details). In order to prove that $\rho^{\frac{1}{2+\de}}v$ belongs in $L^{\infty}((0,T),L^{2+\de}(\Om))$ for $\de$ small enough, it is necessary to control the integral $\int^{T}_{0}\int_{\Om}\big(\rho^{2\gamma-1-\frac{\delta}{2}}\big)^{\frac{2}{2-\de}}dxdt$. Since $\rho$ is bounded in $L^{\frac{5}{3}}((0,T)\times \Om)$, a necessary condition is:
%$$2\gamma-1<\frac{5}{3}\gamma \Leftrightarrow \gamma<3.$$
%It explains in particular why in the theorem \ref{theo1}, we assume that $\gamma<3$ for $N=3$.\\
%For $N=2$, we show that  $\rho$ is bounded in $L^{r}((0,T)\times \Om)$ for any $1\leq r<2$. In particular we have always $2\gamma-1\leq 2\gamma$, that is why we do not need any assumption on $\gamma$ for $N=2$ in the theorem \ref{theo1}.
%\label{remimp}
%\end{remarka}
\subsection{Basic results of compactness}
We would like to finish this section by giving  very useful propositions of compactness that we shall often apply. We are going to recall the so-called Aubin-Lions theorem.
\begin{proposition}
\label{Aubin}
Let $X\hookrightarrow \hookrightarrow B\hookrightarrow Y$ be Banach spaces (with $X$ which is compactly imbedded in $B$) and $(f_{n})_{n\in\mathbb{N}}$ a sequence bounded in $L^{q}((0,T),B)\cap L^{1}((0,T),X)$ (with $1< q\leq +\infty$) and $(\frac{d}{dt}f_{n})_{n\in\mathbb{N}}$ bounded in $L^{1}((0,T),Y)$.Then $(f_{n})_{n\in\mathbb{N}}$ is relatively compact in $L^{p}((0,T),B)$ for any $1\leq p<q$.
\end{proposition}
Let us recall now the theorem of Arz\`ela-Ascoli.
\begin{proposition}
\label{Ascoli}
Let $B$ and $X$ Banach spaces such that  $B\hookrightarrow \hookrightarrow X$ is compact. Let $f_{N}$ be a sequence of functions $\bar{I}\rightarrow B$ (with $I$ an interval) uniformly bounded in $B$ and uniformly continuous in $X$. Then there exists $f\in C^{0}(\bar{I}, B)$ such that $f_{n}\rightarrow f$ strongly in $f\in C^{0}(\bar{I}, X)$ up to a subsequence.
\end{proposition}
\begin{lem}
Let $K$ a compact subset of $\mathbb{R}^{N}$ (with $N\geq 1$) and $v^{\e}$ a sequel such that:
\begin{itemize}
\item $v^{\e}$ is uniformly bounded in $L^{1+\alpha}(K)$ with $\alpha>0$,
\item $v^{\e}$ converge almost everywhere to $v$,
\end{itemize}
then $v^{\e}$ converges strongly to $v$ in $L^{1}(K)$ with $v\in L^{1+\alpha}(K)$.
\label{lemmeimp}
\end{lem}
{\bf Proof}: First by the Fatou lemma $v$ is in $L^{1+\alpha}(K)$. Next we have for any $M>0$:
\begin{equation}
\int_{K}|v^{\e}-v|dx\leq \int_{K\cap\{|v^{\e}-v|\leq M\}}|v^{\e}-v|dx+ \int_{K\cap\{|v^{\e}-v|\geq M\}}|v^{\e}-v|dx.
\label{lemme1}
\end{equation}
We are dealing with the second member of the right hand side, by H\"older inequality and Tchebychev lemma we have for a $C>0$:
\begin{equation}
\begin{aligned}
\int_{K\cap\{|v^{\e}-v|\geq M\}}|v^{\e}-v|dx&\leq (\int_{K}|v^{\e}-v|^{1+\alpha}dx)^{\frac{1}{1+\alpha}}(\{|v^{\e}-v|\geq M\})^{\frac{\alpha}{1+\alpha}},\\
&\leq \frac{C}{M^{\frac{\alpha}{1+\alpha}}}.
\end{aligned}
\label{lemme2}
\end{equation}
In particular we have shown the strong convergence of $v^{\e}$ to $v$, indeed from the inequality (\ref{lemme1}) it suffices to use the Lebesgue theorem for the first term on the right hand side and the estimate (\ref{lemme2}) with $M$ going to $+\infty$. {\hfill $\Box$}
\begin{lem}
Let $f\in\dot{H}^{s}$ with $s>0$ and $f\in L^{p}+L^{2}$ with $1\leq
p<2$. Then $f\in L^{2}$. \label{lemutile}
\end{lem}
{\bf Proof:}
Indeed we have as $f\in\dot{H}^{s}$:
$$\int_{\R^{N}}|\xi|^{2s}|\widehat{f}|^{2}d\xi<+\infty,$$
so $\widehat{f}1_{\{|\widehat{f}|\geq1\}}\in L^{2}(\R^{N})$. And as
$f=f_{1}+f_{2}$ with $f_{1}\in L^{p}(\R^{N})$ and $f_{2}\in L^{2}$.
By using the Riesz-Thorin theorem, we know that
$\widehat{f_{1}}\in L^{q}(\R^{N})$ with $\frac{1}{p}+\frac{1}{q}=1$.
As $q\geq2$ we then have $\widehat{f}1_{\{|\widehat{f}|\leq1\}}\in
L^{2}(\R^{N})$. This concludes the proof.
\hfill {$\Box$}
\section{Proof of theorem \ref{theo1}}
\label{section4}
Let us assume in a first time that the solution $(\rho,u)$ of system (\ref{3systeme}) are classical, we are going to search solution under the form:
$(\rho,-\n\va(\rho))$.
The mass equation give us:
\begin{equation}
\p_{t}\rho-{\rm div}(\rho\n\va(\rho))=0
\label{aporeux1}
\end{equation}
Since $\va^{'}(\rho)=\frac{2\mu^{'}(\rho)}{\rho}$ we get:
\begin{equation}
\p_{t}\rho-2\D\mu(\rho)=0
\label{poreux1}
\end{equation}
Let us check that the second equation is compatible with the first and keep an irrotational structure. First we have:
\begin{equation}
\begin{aligned}
\p_{t}(\rho u)&=-\p_{t}(\rho\n\va(\rho))=-2\n\p_{t}(\mu(\rho)).\\[3mm]
-2{\rm div}(\mu(\rho)Du)&=2{\rm div}(\mu(\rho)\n\n\va(\rho)),\\
&=2\mu(\rho)\n\D\va(\rho)+2\n\mu(\rho)\cdot\n\n\va(\rho).\\[3mm]
-\n(\lambda(\rho){\rm div}u)&=-\n(\lambda(\rho)\D\va(\rho)).
\end{aligned}
\label{1}
\end{equation}
Next we have:
$$
\begin{aligned}
{\rm div}(\rho u\otimes u)_{j}&=\sum_{i}\p_{i}(\rho u_{i}u_{j})=\sum_{i}\p_{i}(\rho \p_{i}\va(\rho)\p_{j}\va(\rho)),\\
&=\sum_{i}\p_{i}(\rho \va^{'}(\rho)\p_{j}\rho \p_{i}\va(\rho))=2\sum_{i}\p_{i}(\p_{j}\mu(\rho) \p_{i}\va(\rho)),\\
&=2\D\va(\rho)\p_{j}\mu(\rho)+2(\n\mu(\rho)\cdot\n\n\va(\rho))_{j}.
\end{aligned}
$$
We have then:
\begin{equation}
{\rm div}(\rho u\otimes u)=2\D\va(\rho)\n\mu(\rho)+2\n\mu(\rho)\cdot\n\n\va(\rho).
\label{2}
\end{equation}
Combining (\ref{1}) and (\ref{2}) we obtain:
\begin{equation}
\begin{aligned}
&{\rm div}(\rho u\otimes u)-2{\rm div}(\mu(\rho)Du)=2\D\va(\rho)\n\mu(\rho)+2\n\mu(\rho)\cdot\n\n\va(\rho)\\
&\hspace{6cm}+2\mu(\rho)\n\D\va(\rho)+2\n\mu(\rho)\cdot\n\n\va(\rho),\\
&=2\n(\mu(\rho)\D\va(\rho))+2\n(\n\mu(\rho)\cdot\n\va(\rho)).
\end{aligned}
\label{3}
\end{equation}
Finally using (\ref{1}), (\ref{3}) and the fact that $\lambda(\rho)+2\mu(\rho)=2\mu^{'}(\rho)$, we obtain:
\begin{equation}
\begin{aligned}
&\p_{t}(\rho u)+{\rm div}(\rho u\otimes u)-2{\rm div}(\mu(\rho) Du)-\n(\lambda(\rho){\rm div}u)=\\[2mm]
&-\n\big(2\p_{t}\mu(\rho)-2\mu(\rho)\D\va(\rho)-2\n\mu(\rho)\cdot\n\va(\rho)-\lambda(\rho)\D\va(\rho)\big),\\
&-\n\big(2\p_{t}\mu(\rho)-2\rho\mu^{'}(\rho)\D\va(\rho)-2\n\mu(\rho)\cdot\n\va(\rho)\big).
\end{aligned}
\label{4}
\end{equation}
Next we have since $\mu^{'}(\rho)=\frac{1}{2}\rho\va^{'}(\rho)$:
$$
\begin{aligned}
\D\mu(\rho)&=\sum_{i}\p_{ii}\mu(\rho)=\sum_{i}\frac{1}{2}\p_{i}(\rho\va^{'}(\rho)\p_{i}\rho),\\
&=\frac{1}{2}\sum_{i}\p_{i}(\rho\p_{i}\va(\rho)),\\
&=\frac{1}{2}\rho\D\va(\rho)+\frac{1}{2}\n\rho\cdot\n\va(\rho).
\end{aligned}
$$
and:
\begin{equation}
\begin{aligned}
\mu^{'}\rho)\D\mu(\rho)&=\frac{1}{2}\rho\D\va(\rho)+\frac{1}{2}\mu^{'}(\rho)\n\rho\cdot\n\va(\rho),\\
&=\frac{1}{2}\rho\D\va(\rho)+\frac{1}{2}\n\mu(\rho)\cdot\n\va(\rho).
\end{aligned}
\label{5}
\end{equation}
In particular from (\ref{5}) we have:
\begin{equation}
4\mu^{'}(\rho)=2\rho\D\va(\rho)+2\n\mu(\rho)\cdot\n\va(\rho)
\label{6}
\end{equation}
Combining (\ref{4}) and (\ref{6}) we have:
\begin{equation}
\begin{aligned}
&\p_{t}(\rho u)+{\rm div}(\rho u\otimes u)-2{\rm div}(\mu(\rho) Du)-\n(\lambda(\rho){\rm div}u)\\
&\hspace{6cm}=-\n\big(2\mu^{'}(\rho)(\p_{t}\rho)-2\D\mu(\rho))\big).
\end{aligned}
\label{4}
\end{equation}
This concludes the proof inasmuch as via the above equation the momentum equation is compatible to the mass equation and must verify the equation (\ref{poreux1}).\\
But when we take initial density in $L^{1}$ non negative and continuous, we know via the theorem \ref{theo2.6} that the unique global solution of (\ref{poreux1}) is classical and non negative. It justify in particular all the previous formal calculus and prove that $(\rho,u=-\n\va(\rho))$ is a classical solution of (\ref{3systeme}) with $\rho$ verifying (\ref{poreux1}).  It concludes the proof.\\
Furthermore the different properties on $\rho$ are a direct consequence of theorem \ref{theo2.4}, \ref{theo2.6} and \ref{theo2.7}.
\section{Proof of  the theorems \ref{theo2}}
\label{section5}
We now present the proof of theorem \ref{theo1} extending to general viscosity coefficient the results of \cite{MV} in the case of the quasi solutions. Let us  begin with recalling the assumptions on the initial data. Indeed we assume that it exits a sequence $(\rho_{n},u_{n})$ of regular global weak solution verifying the system (\ref{3systemev}) (or at least an approximated system, typically by using Friedrich approximations).
\subsubsection*{Initial data:}
In particular the initial data $\rho_{0}^{n},u_{0}^{n})$ must uniformly in $n$ satisfy (\ref{7}), and (\ref{8}) in order to verify the entropy inequalities from section \ref{section3}, more precisely we shall have:
\begin{itemize}
\item $\rho_{0}^{n}$ is bounded in $L^{1}(\Om)$, $\rho_{0}^{n}\geq0$ a.e in $\Om$,
\item $\rho_{0}^{n}|u_{0}^{n}|^{2}$ is bounded in $L^{1}(\Om)$,
\item $\sqrt{\rho^{n}}\n \va(\rho_{0}^{n})=\n f(\rho^{n})$ is bounded in $L^{2}(\Om)$,
\item $\rho_{0}^{n}|u_{0}^{n}|^{2}\ln(1+|u_{0}^{n}|^{2})$ is bounded in $L^{1}(\Om)$.
\end{itemize}
With those assumptions, and using the entropy inequalities (\ref{21}), (\ref{22}) and the mass equation, we have the following bounds:
\begin{equation}
\begin{aligned}
&\|\sqrt{\rho_{n}}\|_{L^{\infty}((0,T),L^{2}(\Om))}\leq C,\\
%&\|\rho_{n}\|_{L^{\infty}((0,T),L^{\gamma}(\Om))}\leq C,\\
&\|\sqrt{\rho_{n}}u_{n}\|_{L^{\infty}((0,T),L^{2}(\Om))}\leq C,\\
&\|\n f(\rho_{n})\|_{L^{\infty}((0,T),L^{2}(\Om))}\leq C,\\
%&\|\sqrt{\rho_{n}}\p_{ij}\ln\rho_{n}\|_{L^{2}((0,T)\times \Om)}\leq C,\\
&\|\sqrt{\rho_{n}}\n u_{n}\|_{L^{2}((0,T)\times\Om)}\leq C,
\end{aligned}
\label{27}
\end{equation}
and %for $\delta$ small enough:
\begin{equation}
\begin{aligned}
%&\|\rho_{n}^{\frac{\gamma}{2}-1}\n\rho_{n}\|_{L^{2}((0,T)\times\Om)}\leq C,\\
&\|\, \rho_{n}|u_{n}|^{2}\ln(1+|u_{n}|^{2})\,\|_{L^{\infty}((0,T),L^{1}(\Om))} \leq C.
\end{aligned}
\label{28}
\end{equation}
\begin{remarka}
Let us point out that the gain of integrability on $u_{n}$ in (\ref{28}) will be a direct consequence of a gain of integrability on the pressure with some restriction on $\gamma$, $\nu_{1}$ and $\nu_{2}$.
\end{remarka}
The proof of theorem \ref{theo1} will be derived in three steps and follows some arguments developed in \cite{MV}. In the first step, we deal with the strong convergence of the density $(\rho_{n})_{n\in\mathbb{N}}$ which enables us to show the convergence almost everywhere of $(\rho_{n})_{n\in\mathbb{N}}$ us to a subsequence. We shall also prove the strong convergence of a momentum sequel of the form $\sqrt{\rho_{n}}h(\rho_{n})u_{n}$ with a function $h$ to precise to a function $\sqrt{\rho}h(\rho)u$. In the second step we derive the strong convergence of $\sqrt{\rho_{n}}u_{n}$ to $\sqrt{\rho}u$  in $L^{2}_{loc}((0,T)\times \R)$ (it allows us to give sense to the momentum product $\rho_{n}u_{n}\otimes u_{n}$) by taking advantage of the uniform gain of integrability on $u{n}$ via the entropy inequality (\ref{gain}). Indeed it will suffice to use the lemma \ref{lemmeimp} after proving almost everywhere convergence via Sobolev injection. In this part, we also shall deal with the strong convergence in the distribution sense of the product $\sqrt{\rho_{n}} \sqrt{\rho_{n}}u_{n}$. In the last step we will treat the diffusion term which will achieve the proof of the theorem \ref{theo2}.
\subsubsection*{Step 1: Convergence almost everywhere on  $\rho_{n}$ and $\rho_{n}u_{n}$}
We are going to begin with proving a technical lemma giving uniforms estimates on $\rho_{n}$ via the entropy (\ref{21}) and (\ref{22}).
\begin{lem}
When $N=2, 3$ $\n(\frac{\mu(\rho_{n})}{\sqrt{\rho_{n}}})$ is uniformly bounded in $L^{\infty}_{T}(L^{2}(\R^{N}))$ for any $T>0$.\\
\\
$\bullet$ When $N=3$ it implies that $\frac{\mu(\rho_{n})}{\sqrt{\rho_{n}}}$, $\frac{\lambda(\rho_{n})}{\sqrt{\rho_{n}}}$ are uniformly bounded in $L^{\infty}((0,T),L^{6}(\R^{N}))$ for any $T>0$ which gives:
\begin{equation}
\; \;\rho_{n}^{\frac{1}{6}+\frac{\nu_{1}}{2N}}\;\;\mbox{is uniformly bounded in}\;\;L^{\infty}(0,T;L^{6}(\R^{N})).\hspace{2cm}
\label{estimimp}
\end{equation}
$\bullet$ When $N=2$ we distinguish two cases:
\begin{itemize}
\item $\nu_{2}\geq 2$, $\rho_{n}$ is uniformly bounded in $L^{\infty}_{T}(L^{q}(\R^{N}))$ for any $q\in[1,+\infty[$ and any $T>0$. It implies that $\frac{\mu(\rho_{n})}{\sqrt{\rho_{n}}}$ and $\frac{\lambda(\rho_{n})}{\sqrt{\rho_{n}}}$ are uniformly bounded in $L^{\infty}((0,T),L^{q}(K))$ for any compact $K$. 
\item $0<\nu_{2}<2$, $\rho_{n}$ is bounded in $L^{\infty}_{T}(L^{q}(\R^{N}))$  for any $T>0$ and any $q\in[1,+\infty[$. It implies that $\frac{\mu(\rho_{n})}{\sqrt{\rho_{n}}}$ and $\frac{\lambda(\rho_{n})}{\sqrt{\rho_{n}}}$ are uniformly bounded in $L^{\infty}((0,T),L^{q}(K))$ for any compact $K$. 
\end{itemize}
\label{tech}
\end{lem}
{\bf Proof:} When $N=3$, we observe that:
$$
\begin{aligned}
\n(\frac{\mu(\rho)}{\sqrt{\rho}})&=2\mu^{'}(\rho)\n\sqrt{\rho}-\frac{\mu(\rho)}{2\rho^{\frac{3}{2}}}\n\rho,\\
&=\frac{1}{2}\n f(\rho)-\frac{\mu(\rho)}{2\rho^{\frac{3}{2}}}\n\rho.
\end{aligned}
$$
and from conditions (\ref{BD}), (\ref{11}) and the fact that $\mu(\rho)\geq 0$ we obtain:
\begin{equation}
2\mu^{'}(\rho)\rho=\lambda(\rho)+2\mu(\rho)=\frac{3\lambda(\rho)+2\mu(\rho)}{3}\geq\frac{\nu}{3}\mu(\rho).
\label{controlmu}
\end{equation}
We deduce that:
$$
\begin{aligned}
|\frac{\mu(\rho)}{2\rho^{\frac{3}{2}}}\n\rho|&\leq \frac{3}{\nu}|\frac{\mu^{'}(\rho)}{\sqrt{\rho}}||\n \rho|,\\
&\leq  \frac{3}{2\nu}|\va^{'}(\rho)\sqrt{\rho}||\n \rho|,\\
&\leq \frac{3}{2\nu}|\n f(\rho)|.
\end{aligned}
$$
It provides that
$$|\n(\frac{\mu(\rho)}{\sqrt{\rho}})|\leq C|\n f(\rho)|.$$
It implies by energy estimates that:
$$\|\n(\frac{\mu(\rho_{n})}{\sqrt{\rho_{n}}})\|_{L^{\infty}(0,T;L^{2}(\R^{N}))}\leq C.$$
Sobolev's embedding ensured that $\frac{\mu(\rho_{n})}{\sqrt{\rho_{n}}}$ is bounded in
$L^{\infty}(0,T;L^{6}(\R^{N}))$. Next by using (\ref{13}) we have:
$$
\begin{aligned}
&C\rho_{n}^{\frac{1}{2}-\frac{1}{N}+\frac{\nu_{1}}{2N}}\leq \frac{\mu(\rho_{n})}{\sqrt{\rho_{n}}}\;\;\;\mbox{when}\;\;\rho_{n}>1,\\
&C\rho_{n}^{\frac{1}{2}-\frac{1}{N}+\frac{\nu_{2}}{2N}}\leq \frac{\mu(\rho_{n})}{\sqrt{\rho_{n}}}\;\;\;\mbox{when}\;\;\rho_{n}\leq1.
\end{aligned}
$$
It implies in particular $\rho_{n}^{\frac{1}{2}-\frac{1}{N}+\frac{\nu_{1}}{2N}}$ is uniformly bounded in $L^{\infty}(0,T;L^{6}(\R^{N}))$ (this is due to the fact that $\rho_{n}1_{\{\rho_{n}\leq 1\}}$ is uniformly bounded in $L^{\infty}(0,T;L^{1}(\R^{N}))\cap L^{\infty}(0,T;L^{\infty}(\R^{N}))$ and that $3-\frac{6}{N}+\frac{2\nu_{1}}{N}\geq 1$ when $N=3$).\\
\\
When $N=2$ by following the same proof than in the case $N=3$ we obtained that $\n(\frac{\mu(\rho)}{\sqrt{\rho}})$ is uniformly bounded in $L^{\infty}(0,T;L^{2}(\R^{N}))$ for any $T>0$. By using in a similar way (\ref{controlmu}) when $N=2$ and (\ref{13}) we have:
\begin{equation}
\begin{aligned}
&C\rho_{n}^{-1+\frac{\nu_{1}}{4}}   |\n\rho_{n}|\leq 2|\frac{\mu^{'}(\rho_{n})}{\sqrt{\rho_{n}}}||\n\rho_{n}|=|\n f(\rho_{n})|\;\;\;\mbox{when}\;\;\rho_{n}>1,\\
&C\rho_{n}^{-1+\frac{\nu_{2}}{4}}  |\n\rho_{n}|\leq 2|\frac{\mu^{'}(\rho_{n})}{\sqrt{\rho_{n}}}||\n\rho_{n}|=|\n f(\rho_{n})|\;\;\;\mbox{when}\;\;\rho_{n}\leq1.
\end{aligned}
\label{techi}
\end{equation}
When $\nu_{1}\geq 2$, choosing $\psi\in C^{\infty}_{0}(\R^{N})$ with $\psi=1$ on $B(0,1)$ and $\mbox{supp}\,\psi$ included in $B(0,2)$ we have:
$(1-\psi(\rho_{n}))\n\sqrt{\rho_{n}}$ is uniformly bounded in $L^{\infty}_{T}(L^{2}(\R^{N}))$ for any $T>0$. Since $\sqrt{\rho_{n}}$ is uniformly bounded in $L^{\infty}_{T}(L^{2}(\R^{N}))$ for any $T>0$ we deduce that $(1-\psi(\rho_{n}))\sqrt{\rho_{n}}$ is uniformly bounded in $L^{\infty}_{T}(H^{1}(\R^{N}))$ for any $T>0$. It implies that $(1-\psi(\rho_{n}))\rho_{n}$ is bounded in $L^{\infty}_{T}(L^{q}(\R^{N}))$ for any $q\in[1,+\infty[$ by Sobolev embedding. Let us deal now with the term $\psi(\rho_{n})\rho_{n}$ which is bounded in $L^{\infty}_{T}(L^{1}(\R^{N}))\cap L^{\infty}_{T}(L^{\infty}(\R^{N}))$, it proves that $\rho_{n}$ is bounded in $L^{\infty}_{T}(L^{q}(\R^{N}))$ for any $q\in[1,+\infty[$ . Via (\ref{13}) It implies that $\frac{\mu(\rho_{n})}{\sqrt{\rho_{n}}}$ and $\frac{\lambda(\rho_{n})}{\sqrt{\rho_{n}}}$ are uniformly bounded in $L^{\infty}((0,T),L^{q}(\R^{N}))$ for any compact $K$. \\
Let us deal with the case $0<\nu_{2}<2$. By using (\ref{techi}) we show that $\n\big( (1-\psi(\rho_{n}))\rho_{n}^{\frac{\nu_{1}}{4}}\big)$ is bounded in $L^{\infty}_{T}(L^{2}(\R^{N}))$ and $(1-\psi(\rho_{n}))\rho_{n}^{\frac{\nu_{1}}{4}}$ is bounded in $L^{\infty}_{T}(L^{\frac{4}{\nu_{1}}}(\R^{N}))$. Since by Tchebytchev lemma $(1-\psi(\rho_{n}))\rho_{n}^{\frac{\nu_{1}}{4}}$ is strictly positive only on a set of finite measure it implies that  $(1-\psi(\rho_{n}))\rho_{n}^{\frac{\nu_{1}}{4}}$ is also bounded in $L^{\infty}_{T}(L^{2}(\R^{N}))$. We deduce that 
$(1-\psi(\rho_{n}))\rho_{n}^{\frac{\nu_{1}}{4}}$ is  bounded in $L^{\infty}_{T}(H^{1}(\R^{N}))$. By Sobolev embedding it yields that $ (1-\psi(\rho_{n}))\rho_{n}^{\frac{\nu_{1}}{4}}$ is bounded in $L^{\infty}_{T}(L^{q}(\R^{N}))$ for any  $T>0$ and any $q\in[1,+\infty[$. Since $\psi(\rho_{n})\rho_{n}$ is bounded in  $L^{\infty}_{T}(L^{\infty}(\R^{N})\cap L^{1}(\R^{N}) )$ for  any $T>0$ we conclude that $\rho_{n}$ is bounded in $L^{\infty}_{T}(L^{q}(\R^{N}))$  for any $T>0$ and any $q\in[1,+\infty[$.
\\
\\
The proof is similar for $\frac{\lambda(\rho_{n})}{\sqrt{\rho_{n}}}$ by using the remarks \ref{r1}. \null{\hfill $\Box$}
\begin{lemme}
If $\mu(\rho)$, $\lambda(\rho)$ satisfies (\ref{BD}), (\ref{11})  and in addition we assume that $g(x)=\frac{\mu(x)}{\sqrt{x}}$ is a bijective function on $[0,+\infty)$ and that $g^{-1}$ is continuous, then when we distinguish the two following cases, we have:
\begin{itemize}
\item $\nu_{1}\geq 2$
\begin{enumerate}
\item $\frac{\mu(\rho_{n})}{\sqrt{\rho_{n}}}$ is uniformly bounded in $L^{\infty}(0,T;H^{1}_{loc}(\R^{N}))$
\item $\p_{t}\frac{\mu(\rho_{n})}{\sqrt{\rho_{n}}}$ is bounded in $L^{2}(0,T;W_{loc}^{-1,2}(\R^{N}))$.
\end{enumerate}
As a consequence up to a subsequence (via the Cantor's diagonal process) $\frac{\mu(\rho_{n})}{\sqrt{\rho_{n}}}$ converges almost everywhere and strongly in $C([0,T],L^{2}_{loc}(\R^{N})))$ to $v$.  In particular it implies that $\frac{\mu(\rho_{n})}{\sqrt{\rho_{n}}}$ converges up to a subsequence almost everywhere to $v$ and we define $\rho$ as follows:
$$\rho=g^{-1}(v).$$
It implies that $\rho_{n}$ converge up to a subsequence almost everywhere to $\rho$.
% since $g(x)=\frac{\mu(x)}{\sqrt{x}}$ is a bijective function on $[0,+\infty)$.
\item $0<\nu_{1}<2$\\
Let us consider:
$$\beta(\rho_{n})=\psi(\rho_{n})\sqrt{\rho_{n}}+(1-\psi(\rho_{n}))\rho_{n}^{\alpha_{1}},$$
with $0<\alpha_{1}<\min(\frac{1}{3},\frac{\nu_{1}}{4N})$.
%Let us now assume that $\mu^{'}(\rho)\geq\nu_{3}$ with $\nu_{3}\in(0,1)$
\begin{enumerate}
\item $\beta(\rho_{n})$ is bounded in $L^{\infty}(0,T;H^{1}(\R^{N}))$
\item $\p_{t}\beta(\rho_{n})$ is bounded in $L^{2}(0,T;W^{-1,1}(\R^{N}))$.
\end{enumerate}
As a consequence up to a subsequence, $\beta(\rho_{n})$ converge almost everywhere and strongly in $C([0,T],L^{2}_{loc}(\R^{N})))$ to $v$. We define $\rho$ by:
$$\rho=\beta^{-1}(\rho).$$
In particular we have that $\rho_{n}$ converge up to a subsequence almost everywhere to $\rho$.  
\end{itemize}
%As a consequence up to a subsequence, $\sqrt{\rho_{n}}$ converge almost everywhere and strongly in $L^{2}(0,T;L^{2}_{loc})$. We write:
%$$\sqrt{\rho_{n}}\rightarrow\sqrt{\rho}\;\;\;\mbox{a.e and}\;\;L^{2}_{loc}\;\;\mbox{strong}.$$
Furthermore $\rho_{n}$ converge strongly to $\rho$ in $C([0,T],L^{1+\alpha}_{loc}(\R^{N}))$ if $N=3$ with $\alpha>0$ small enough and in $C([0,T],L^{p}_{loc}(\R^{N}))$ for any $p\geq 1$ if $N=2$. This last result is under the following hypothesis:
\begin{itemize}
\item When $2+N \leq \nu_{1} $, we assume that $g$ and $g^{'}$ are increasing on $(0,+\infty)$.
\end{itemize}
\label{rho}
\end{lemme}
\begin{remarka}
Let us point out that we could weaken the last assumption on $g$ when $\nu_{1}\geq 2$ by assuming that $g$ and $g^{'}$ are only increasing in a neighbor
of $0$ and of $+\infty$. As mentioned above, let us point out that this last condition is quite natural since this is true when we set $\mu(\rho)=\mu\rho^{\alpha}$ with $\alpha\geq\frac{3}{2}$.
\end{remarka}
{\bf Proof:}
Let $T>0$. We are going to distinguish two case when $\nu_{1}\geq 2$ and when $0<\nu_{1}<2$.\\
\\
$\bullet \nu_{1}\geq 2$\\
\\
The first estimate is a direct consequence of the lemma \ref{tech} in the appendix. Indeed we know that $\n \frac{\mu(\rho_{n})}{\sqrt{\rho_{n}}}$ is bounded in $L^{\infty}_{T}(L^{2}(\R^{N}))$. Easily we deduce by lemma \ref{tech} that $\frac{\mu(\rho_{n})}{\sqrt{\rho_{n}}}$ is uniformly bounded in $L^{\infty}_{T}(L^{2}(K))$ for any compact $K$.
%Furthermore via the lemma \ref{tech} we have $\frac{\mu(\rho_{n})}{\sqrt{\rho_{n}}}$ is uniformly bounded in $L^{\infty}_{T}(L^{6}(\R^{N})$ which implies that $\frac{\mu(\rho_{n})}{\sqrt{\rho_{n}}}$
%is uniformly bounded in $L^{\infty}_{T}(L^{2}_{loc}(\R^{N})$. It provides the %\\
%The second estimate in \ref{}, together with the conservation of mass $\|\rho_{n}\|_{L^{1}}=\|\rho_{n,0}\|_{L^{1}}$ gives the 
%in $L^{\infty}(0,T;H^{1}_{loc})$
%bound. 
Next we observe that:
\begin{equation}
\begin{aligned}
\p_{t}(\frac{\mu(\rho_{n})}{\sqrt{\rho_{n}}})&=-{\rm div}(\frac{\mu(\rho_{n})}{\sqrt{\rho_{n}}}u_{n})+\big(\mu^{'}(\rho_{n})\sqrt{\rho_{n}}-\frac{3}{2}\frac{\mu(\rho_{n})}{\sqrt{\rho_{n}}}\big){\rm div}u_{n},\\
&=-{\rm div}(\frac{\mu(\rho_{n})}{\rho_{n}}\sqrt{\rho_{n}}u_{n})+\big(\frac{\mu^{'}(\rho_{n})\sqrt{\rho_{n}}}{\sqrt{\mu(\rho_{n})}} -\frac{3}{2}\sqrt{ \frac{\mu(\rho_{n})}{\rho_{n}}}\big)\sqrt{\mu(\rho_{n})}{\rm div}u_{n}.
\end{aligned}
\label{visco}
\end{equation}
Let us start with estimating the first term on the right hand side of (\ref{visco}).
According to (\ref{13}) we know that:
$$\frac{\psi(\rho_{n})\mu(\rho_{n})}{\rho_{n}}\leq C\psi(\rho_{n})\rho_{n}^{\frac{\nu_{1}}{2N}-\frac{1}{N}}\;\;\mbox{when}\;\;\rho_{n}\leq 1.$$
Since $\frac{\nu_{1}}{2N}-\frac{1}{N}\geq 0$ it implies in particular that $\frac{\psi(\rho_{n})\mu(\rho_{n})}{\rho_{n}}$ is bounded in $L^{\infty}_{T}(L^{\infty}(\R^{N}))$. In order to have local estimates, it suffices to deal with the region when $\rho_{n}>1$, it means $\frac{(1-\psi(\rho_{n}))\mu(\rho_{n})}{\rho_{n}}$. \\
It suffices to use the lemma \ref{tech} which insures that when $N=3$ $\frac{\mu(\rho_{n})}{\sqrt{\rho_{n}}}$ belongs in $L^{\infty}_{T}(L^{6}(\R^{N}))$ and we deduce that
$\frac{(1-\psi(\rho_{n}))\mu(\rho_{n})}{\rho_{n}}$ is bounded in $L^{\infty}_{T}(L^{6}(\R^{N}))$. Finally we have obtain that  $\frac{\mu(\rho_{n})}{\rho_{n}}$ is bounded in $L^{\infty}_{T}(L^{2}_{loc}(\R^{N}))$ when $N=3$. It yields the uniform boundness of $\frac{\mu(\rho_{n})}{\rho_{n}}\sqrt{\rho_{n}}u_{n}$ in  $L^{\infty}_{T}(L^{1}_{loc}(\R^{N}))$ when $N=3$. A similar argument insure the same result when $N=2$.\\
By proceeding similarly we also prove that  $\sqrt{ \frac{\mu(\rho_{n})}{\rho_{n}}}\sqrt{\mu(\rho_{n})}{\rm div}u_{n}$ is uniformly bounded in  $L^{\infty}_{T}(L^{1}_{loc}(\R^{N}))$. Indeed following the same argument than for the previous term, we know that for $N=2,3$ via the lemma \ref{tech} $\sqrt{ \frac{\mu(\rho_{n})}{\rho_{n}}}$ is bounded in 
$L^{\infty}_{T}(L^{12}(K))$ for any compact $K$ and since $\sqrt{\mu(\rho_{n})}{\rm div}u_{n}$ is bounded in $L^{2}_{T}(L^{2}(\R^{N}))$, H\"older(s inequality give the desired result. \\
Let us now deal with the term $\frac{\mu^{'}(\rho_{n})\sqrt{\rho_{n}}}{\sqrt{\mu(\rho_{n})}}\sqrt{\mu(\rho_{n})}{\rm div}u_{n}$. A simple calculus using (\ref{BD}) give us:
$$\frac{\mu^{'}(\rho_{n})\sqrt{\rho_{n}}}{\sqrt{\mu(\rho_{n})}}=\frac{1}{2}(\frac{\lambda(\rho_{n})}{\sqrt{\rho_{n}\mu(\rho_{n})}}+2\sqrt{ \frac{\mu(\rho_{n})}{\rho_{n}}}).$$
We have only to deal with the term $\frac{\lambda(\rho_{n})}{\sqrt{\rho_{n}\mu(\rho_{n})}}$, the other one has been estimated. By the remarks \ref{r2.2} we know that it exists $C>0$ such that:
$$|\lambda(\rho)|\leq C\mu(\rho)\;\;\forall\rho>0.$$
It implies that:
$$|\frac{\lambda(\rho_{n})}{\sqrt{\rho_{n}\mu(\rho_{n})}}|\leq C\sqrt{ \frac{\mu(\rho_{n})}{\rho_{n}}}.$$
And as we know that $ \frac{\mu(\rho_{n})}{\sqrt{\rho_{n}}}$ is uniformly bounded in $L^{\infty}_{T}(L^{6}_{loc}(\R^{N}))$ for $N=2,3$ via the lemma \ref{tech}, it achieves the proof of the second estimates.\\ 
By the Ascoli's theorem, the fact that the application $u\rightarrow \phi u$ with $\phi\in C^{\infty}_{0}(\R^{N})$ is compact from $H^{1}(\R^{N})$ to $L^{2}(\R^{N})$ and the Cantor's diagonal process it entails that  $ \frac{\mu(\rho_{n})}{\sqrt{\rho_{n}}}$ converges strongly up to a subsequence in $C([0,T],L^{2}_{loc}(\R^{N}))$) to $v_{1}= \frac{\mu(\rho)}{\sqrt{\rho}}$ (and in particular  in $L^{2}_{loc}((0,T)\times \R^{N})$ ). We shall define $\rho$ in the sequel by:
$$\rho=g^{-1}(\frac{\mu(\rho)}{\sqrt{\rho}}).$$
Furthermore an immediate  consequence is that up to a subsequence $ \frac{\mu(\rho_{n})}{\sqrt{\rho_{n}}}$ converge almost everywhere to $\frac{\mu(\rho)}{\sqrt{\rho}}$.  And since $g^{-1}$ is continuous, it implies that up to a subsequence $\rho_{n}$ converges almost everywhere to $\rho$. \\
\\
$\bullet$ $0<\nu_{1}<2$\\
\\
In this case we are going to study:
$$\beta(\rho_{n})=\psi(\rho_{n})\sqrt{\rho_{n}}+(1-\psi(\rho_{n}))\rho_{n}^{\alpha_{1}},$$
with $\alpha_{1}$ to choose suitably.  Let us start with dealing with $\psi(\rho_{n})\sqrt{\rho_{n}}$, we have then:
$$\n(\psi(\rho_{n})\sqrt{\rho_{n}})=\psi^{'}(\rho_{n})\sqrt{\rho_{n}}\n\rho_{n}+\frac{\psi(\rho_{n})}{\mu^{'}(\rho_{n})}\frac{\mu^{'}(\rho_{n})}{\sqrt{\rho_{n}}}\n\rho_{n}$$
The first term is easily bounded in $L^{\infty}_{T}(L^{2}(\R^{N}))$ via the entropy (\ref{22}) since the support of $\psi^{'}$ is included in the shell $C(0,1,2)$.The second term is also bounded in $L^{\infty}_{T}(L^{2}(\R^{N}))$ by observing that $\frac{\psi(\rho_{n})}{\mu^{'}(\rho_{n})}$ is bounded in $L^{\infty}_{T}(L^{\infty}(\R^{N}))$ via (\ref{13}), (\ref{BD}), (\ref{11}) and $\nu_{1}\leq 2$. Indeed we have:
$$|\frac{\psi(\rho_{n})}{\mu^{'}(\rho_{n})}|\leq \psi(\rho_{n})\rho_{n}^{\frac{1}{N}-\frac{\nu_{1}}{2N}}.$$
The conservation of mass provides that  $\|\rho_{n}\|_{L^{1}}=\|\rho^{n}_{0}\|_{L^{1}}$which implies the $L^{\infty}(0,T;H^{1})$
bound on $\psi(\rho_{n})\sqrt{\rho_{n}}$. Next we observe that:
\begin{equation}
\begin{aligned}
\p_{t}(\psi(\rho_{n})\sqrt{\rho_{n}})&=-{\rm div}(\psi(\rho_{n})\sqrt{\rho_{n}}u_{n})+(\frac{1}{2}\psi(\rho_{n})\sqrt{\rho_{n}}-\psi^{'}(\rho_{n})\rho^{\frac{3}{2}}_{n}){\rm div}u_{n}.
%\frac{1}{2}\psi(\rho_{n})\sqrt{\rho_{n}}{\rm div}u_{n}-\psi(\rho_{n})u_{n}\cdot\n\sqrt{\rho_{n}}+\psi^{'}(\rho_{n})\sqrt{\rho_{n}}\p_{t}\rho_{n}.\\
%&=\frac{1}{2}\sqrt{\rho_{n}}{\rm div}u_{n}-{\rm div}(u_{n}\sqrt{\rho_{n}}),
\end{aligned}
\label{conv22}
\end{equation}
The first term on the right hand side is obviously bounded in $L^{\infty}_{T}(W^{-1,2}(\R^{N}))$. The two last terms in (\ref{conv22}) are bounded in $L^{2}_{T}(W^{-1,2}(\R^{N}))$. It implies finally that $\psi(\rho_{n})\sqrt{\rho_{n}}$ are uniformly bounded in $L^{\infty}_{T}(H^{1}(\R^{N}))$ and $\p_{t}(\psi(\rho_{n})\sqrt{\rho_{n}})$ in $L^{2}_{T}(W^{-1,2}(\R^{N}))$.\\
Let us deal now with the term $(1-\psi(\rho_{n}))\rho_{n}^{\alpha_{1}}$, we have then:
\begin{equation}
\n((1-\psi(\rho_{n}))\rho_{n}^{\alpha_{1}})=-=\psi^{'}(\rho_{n})\rho_{n}^{\alpha_{1}}\n \rho_{n}+(1-\psi(\rho_{n}))\frac{\rho_{n}^{\alpha_{1}-\frac{1}{2}}}{\mu^{'}(\rho_{n})}
\,\frac{\mu^{'}(\rho_{n})}{\sqrt{\rho_{n}}}\n \rho_{n}.
\label{conv33}
\end{equation}
By (\ref{BD}), (\ref{11}) and (\ref{13}) we show that:
$$
|(1-\psi(\rho_{n}))\frac{\rho_{n}^{\alpha_{1}-\frac{1}{2}}}{\mu^{'}(\rho_{n})}|\leq (1-\psi(\rho_{n}))\rho_{n}^{\alpha_{1}-\frac{1}{2}+\frac{1}{N}-\frac{\nu_{1}}{2N}}.
$$
Since $0<\alpha_{1}\leq\frac{\nu_{1}}{2N}$ it implies that $(1-\psi(\rho_{n}))\frac{\rho_{n}^{\alpha_{1}-\frac{1}{2}}}{\mu^{'}(\rho_{n})}$ is bounded in $L^{\infty}_{T}(L^{\infty}(\R^{N}))$ and $(1-\psi(\rho_{n}))\frac{\rho_{n}^{\alpha_{1}-\frac{1}{2}}}{\mu^{'}(\rho_{n})}
\,\frac{\mu^{'}(\rho_{n})}{\sqrt{\rho_{n}}}\n \rho_{n}$ is bounded in $L^{\infty}_{T}(L^{2}(\R^{N}))$. The first term on (\ref{conv33}) is easy to deal with. Now we have:
\begin{equation}
\begin{aligned}
\p_{t}((1-\psi(\rho_{n}))\rho_{n}^{\alpha_{1}})&=-{\rm div}((1-\psi(\rho_{n}))\rho_{n}^{\alpha_{1}}u_{n})-(-\psi^{'}(\rho_{n})\frac{\rho_{n}^{\alpha_{1}+1}}{\sqrt{\mu(\rho_{n})}} \\
&\hspace{2cm}+(\alpha_{1}-1)(1-\psi(\rho_{n}))\frac{\rho_{n}^{\alpha_{1}}}{\sqrt{\mu(\rho_{n})}}    )\sqrt{\mu(\rho_{n})}{\rm div}u_{n}.
%\frac{1}{2}\psi(\rho_{n})\sqrt{\rho_{n}}{\rm div}u_{n}-\psi(\rho_{n})u_{n}\cdot\n\sqrt{\rho_{n}}+\psi^{'}(\rho_{n})\sqrt{\rho_{n}}\p_{t}\rho_{n}.\\
%&=\frac{1}{2}\sqrt{\rho_{n}}{\rm div}u_{n}-{\rm div}(u_{n}\sqrt{\rho_{n}}),
\end{aligned}
\label{conv223}
\end{equation}
The first term on the right hand side of (\ref{conv223}) is bounded in $L^{\infty}_{T}(W^{-1,1}(\R^{N}))$ since  $(1-\psi(\rho_{n}))\rho_{n}^{\alpha_{1}}$ is  bounded in $L^{\infty}(L^{2}(\R^{N}))$ because $\alpha_{1}\leq\frac{1}{2}$. We have now according to (\ref{13})
$$|(1-\psi(\rho_{n}))\frac{\rho_{n}^{\alpha_{1}}}{\sqrt{\mu(\rho_{n})}}   |\leq (1-\psi(\rho_{n}))\rho_{n}^{\alpha_{1}-\frac{1}{2}+\frac{1}{2N}-\frac{\nu_{1}}{4N}}.$$
Since $\alpha_{1}\leq\frac{1}{3}$ it implies that this term is bounded in $L^{\infty}_{T}(L^{\infty}(\R^{N}))$ which shows the bound of $(1-\psi(\rho_{n}))\frac{\rho_{n}^{\alpha_{1}}}{\sqrt{\mu(\rho_{n})}}    \sqrt{\mu(\rho_{n})}{\rm div}u_{n}$ in $L^{2}_{T}(L^{2}(\R^{N}))$. The second term of \ref{conv223} is easy to treat.\\
Finally we have proved that  $(1-\psi(\rho_{n}))\rho_{n}^{\alpha_{1}}$ is uniformly bounded in $L^{\infty}_{T}(H^{1}(\R^{N}))$ and $\p_{t}((1-\psi(\rho_{n}))\rho_{n}^{\alpha_{1}})$ is also in $L^{2}_{T}(W^{-1,2}(\R^{N}))$. In conclusion it shows that $\beta(\rho_{n})$ is bounded in $L^{\infty}_{T}(H^{1}(\R^{N}))$ and $\p_{t}(\beta(\rho_{n}))$ is bounded $L^{2}_{T}(W^{-1,1}(\R^{N}))$.
\\
%which gives the second estimate. 
Thanks to the Ascoli's theorem and the Cantor's diagonal process it gives the strong convergence in $$C([0,T],L^{2}_{loc}(\R^{N}))$$ of $\beta(\rho_{n})$ to $v_{2}$. It implies in particular up a subsequence of the convergence almost everywhere of $\beta(\rho_{n})$ to $v_{2}$. We shall define $\rho$ in this situation by:
$$\rho=\beta^{-1}(\rho),$$
indeed $\beta$ is inversible and we verify by continuity of $\beta^{-1}$ that $\rho_{n}$ converge up to a subsequence almost everywhere to $\rho$.% (in fact to be more accurate from $C([0,T],L^{2}_{loc}(\R^{N}))$). The last terms in 
%
%which yields the second estimate and, thanks to Aubin's Lions Lemma, gives the strong convergence in $L^{2}_{loc}$.\\
\subsubsection*{Strong convergence of $\rho_{n}$}
We are now interested in proving the strong convergence of $\rho_{n}$ to $\rho$. Let us deal with the first case $\nu_{1}\geq 2$.\\
\\
\textbf{First case: $\nu_{1}\geq 2$}
%$\bullet$ $\nu_{1}\geq 2$.\\
\\
\\
For the moment we have only obtained strong convergence on $\frac{\mu(\rho_{n})}{\sqrt{\rho_{n}}}$ to $\frac{\mu(\rho)}{\sqrt{\rho}}$, let us translate this strong convergence on $\rho_{n}$. We are going to distinguish two different cases when $\nu_{1}\geq 2$, let us start with the first one.\\
\\
$\bullet$ We assume here that $(\frac{\mu(x)}{\sqrt{x}})^{'}$, $\frac{\mu(x)}{\sqrt{x}}$ are increasing on $(0,+\infty)$. We have then the following lemma.
%\begin{lem}
%When $x,y\geq 0$ and $\alpha\geq 1$, we have:
%$$|x-y|^{\alpha}\leq |x^{\alpha}-y^{\alpha}|.$$
%\end{lem}
\begin{lem}
Let $g_{1}$ a regular function with $g_{1}(0)=0$. When $x,y\geq 0$ and $g_{1}$, $g_{1}^{'}$ are increasing, we have:
\begin{equation}
g_{1}(|x-y|)\leq |g_{1}(x)-g_{1}(y)|.
\label{astuce}
\end{equation}
\label{ruse}
\end{lem}
{\bf Proof:} %It suffices to study the function:
%$$f(x)=(x-y)^{\alpha}-(x^{\alpha}-y^{\alpha}),$$
%when $x\geq y$. We observe that for all $x\geq y\geq 0$:
%$$f^{'}(x)=\alpha((x-y)^{\alpha-1}-x^{\alpha-1})\leq 0.$$
%It implies that when $x\geq y$,
%$$|x-y|^{\alpha}\leq |x^{\alpha}-y^{\alpha}|.$$
%We are doing similarly when $0\leq x\leq y$.
It suffices to study the function:
$$p(x)=g_{1}(x-y)-(g_{1}(x)-g_{1}(y)),$$
when $x\geq y$. We observe that for all $x\geq y\geq 0$:
$$p^{'}(x)=g_{1}^{'}(x-y)-g_{1}^{'}(x)\leq 0.$$
It implies that when $x\geq y$, $p^{'}$ is negative and $p$ is decreasing on $[y,+\infty)$ with is equivalent to say that:
$$g_{1}(x-y)-(g_{1}(x)-g_{1}(y))\leq p(y)=g_{1}(0)\;\;\forall x\geq y.$$
It implies that for all $x\geq y$ we have:
$$g_{1}(|x-y|)=g_{1}(x-y)\leq (g_{1}(x)-g_{1}(y))=|g_{1}(x)-g_{1}(y)|.$$
By proceeding similarly when $0\leq x\leq y$ we obtain (\ref{astuce}).
{\hfill $\Box$}\\
\\
In particular  since we assume that $(\frac{\mu(x)}{\sqrt{x}})^{'}$ and $(\frac{\mu(x)}{\sqrt{x}})$ are increasing on $(0,+\infty)$ we deduce from the lemma \ref{ruse}:
$$\frac{\mu(|\rho-\rho_{n}|)}{\sqrt{|\rho-\rho_{n}|}}\leq |\frac{\mu(\rho_{n})}{\sqrt{\rho_{n}}}-\frac{\mu(\rho)}{\sqrt{\rho}}|$$
Using (\ref{13}) we obtain:
\begin{equation}
\begin{aligned}
&|\rho-\rho_{n}|^{\frac{1}{2}-\frac{1}{N}+\frac{\nu_{1}}{2N}}\leq \frac{1}{C} |\frac{\mu(\rho_{n})}{\sqrt{\rho_{n}}}-\frac{\mu(\rho)}{\sqrt{\rho}}|\;\;\forall\;|\rho-\rho_{n}|>1\\
&|\rho-\rho_{n}|^{\frac{1}{2}-\frac{1}{N}+\frac{\nu_{2}}{2N}}\leq \frac{1}{C} |\frac{\mu(\rho_{n})}{\sqrt{\rho_{n}}}-\frac{\mu(\rho)}{\sqrt{\rho}}|\;\;\forall\;|\rho-\rho_{n}|\leq1.
\end{aligned}
\label{utile}
\end{equation}
Since $\frac{\mu(\rho_{n})}{\sqrt{\rho_{n}}}$ converges strongly to $\frac{\mu(\rho)}{\sqrt{\rho}}$ in $C([0,T],L^{2}_{loc})$ and $\frac{\mu(\rho_{n})}{\sqrt{\rho_{n}}}$, $\frac{\mu(\rho)}{\sqrt{\rho}}$ are bounded via lemma \ref{tech} in $L_{T}^{\infty}(L^{6}(\R^{N}))$  when $N=3$ (it suffices to apply Fatou's lemma), we deduce by interpolation that $\frac{\mu(\rho_{n})}{\sqrt{\rho_{n}}}$ converges strongly to $\frac{\mu(\rho)}{\sqrt{\rho}}$ in $C([0,T],L^{6-\alpha}_{loc})$ for any small $\alpha>0$. Similarly when $N=2$ we obtain by interpolation and via lemma \ref{tech} that $\frac{\mu(\rho_{n})}{\sqrt{\rho_{n}}}$ converges strongly to $\frac{\mu(\rho)}{\sqrt{\rho}}$ in $C([0,T],L^{p}_{loc})$ for any $p>1$.\\
By using (\ref{utile}) we deduce that $(1-\psi(\rho-\rho_{n}))|\rho-\rho_{n}|^{\frac{1}{2}-\frac{1}{N}+\frac{\nu_{1}}{2N}}$ and $\psi(\rho-\rho_{n})|\rho-\rho_{n}|^{\frac{1}{2}-\frac{1}{N}+\frac{\nu_{2}}{2N}}$ converge strongly to $0$ in $C([0,T],L^{6-\alpha}_{loc})$ for any small $\alpha>0$ when $N=3$ and in $C([0,T],L^{p}_{loc})$ for any $p>1$ when $N=2$. Since $\nu_{1}>0$  it yields that:
\begin{equation}
\sup_{t\in[0,T]}\||\rho-\rho_{n}|^{p(N)(\frac{1}{2}-\frac{1}{N})}1_{\{|\rho-\rho_{n}|>1\}}(t,\cdot)\|_{L^{1}(K)}\rightarrow_{n\rightarrow+\infty}0,
\label{bout1}
\end{equation}
for any compact $K$ with $p(N)=6-\alpha$ for any small $\alpha>0$ and $p(N)=p$ for any $p\in[1,+\infty[$ when $N=2$. And similarly we have:
\begin{equation}
\sup_{t\in[0,T]}\||\rho-\rho_{n}|^{p(N)(\frac{1}{2}-\frac{1}{N}+\frac{\nu_{2}}{2N})}1_{\{|\rho-\rho_{n}|>1\}}(t,\cdot)\|_{L^{1}(K)}\rightarrow_{n\rightarrow+\infty}0.
\label{bout1}
\end{equation}
It implies since $\nu_{2}>0$ that $\rho_{n}$ converges strongly to $\rho$ in $C([0,T],L^{1+\alpha}_{loc}(\R^{N}))$ for $\alpha>0$ small enough when $N=3$ and $\rho_{n}$ converges strongly to $\rho$ in $C([0,T],L^{p}_{loc}(\R^{N}))$ for any $p\geq 1$.\\
\\
$\bullet$ $\nu_{2}\leq N+2$\\
\\
%that it implies that $(1-\psi(\rho-\rho_{n}))|\rho-\rho_{n}|^{\frac{1}{2}-\frac{1}{N}}$ converges strongly in $C([0,T],L^{2}_{loc}(R^{N}))$.
%Since $\nu_{1}>0$ and $\frac{\mu(\rho_{n})}{\sqrt{\rho_{n}}}$ converges strongly to $\frac{\mu(\rho)}{\sqrt{\rho}}$ in $C([0,T],L^{2}_{loc})$ it implies that $(1-\psi(\rho-\rho_{n}))|\rho-\rho_{n}|^{\frac{1}{2}-\frac{1}{N}}$ converges strongly in $C([0,T],L^{2}_{loc}(R^{N}))$.
%$\frac{1}{2}-\frac{1}{N}+\frac{\nu_{2}}{2N}>\frac{1}{2}-\frac{1}{N}$ a
% (and by interpolation in $C([0,T],L^{p}_{loc})$ for any $p>1$ when $N=2$ and in $C([0,T],L^{6-\e}_{loc})$
%when $N=3$) it implies also the strong convergence of $\rho_{n}$ to $\rho$ in $C([0,T],L^{1+\e}_{loc})$ with $\e>0$ when $N=3$ and in $C([0,T],L^{p}_{loc})$ for any $p>1$ when $N=2$.\\[2mm]
% since via the lemma \ref{tech} $\frac{\mu(\rho_{n})}{\sqrt{\rho_{n}}}$ is uniformly bounded in $L^{\infty}_{T}(L^{6}(\R^{N}))$ , it implies that $\rho_{n}$ converges strongly to $\rho$ in $C([0,T],L^{1}_{loc}(\R^{N}))$.\\
%MORALEMENT CE CAS EST VRAI LORSQUE $\alpha\geq 3$.\\
Let us deal now with the case when $(\frac{\mu(x)}{\sqrt{x}})^{'}$ and $(\frac{\mu(x)}{\sqrt{x}})$ are not increasing on $(0,+\infty)$.
 By calculus and using (\ref{BD})  we have:
$$
\begin{aligned}
(\frac{\mu(\rho)}{\sqrt{\rho}})^{'}&=\frac{\mu^{'}(\rho)}{\sqrt{\rho}}-\frac{1}{2}\frac{\mu(\rho)}{\rho^{\frac{5}{2}}},\\
&=\frac{1}{2}\frac{\lambda(\rho)+\mu(\rho)}{\rho^{\frac{3}{2}}}.
\end{aligned}
$$
Next by (\ref{13}) we obtain:
\begin{equation}
\begin{aligned}
&C\rho^{\frac{\nu_{2}}{2N}-\frac{1}{2}-\frac{1}{N}}\leq |(\frac{\mu(\rho)}{\sqrt{\rho}})^{'}|\leq C\rho^{\frac{\nu_{1}}{2N}-\frac{1}{2}-\frac{1}{N}}\;\;\forall \rho\leq 1,\\
&C\rho^{\frac{\nu_{1}}{2N}-\frac{1}{2}-\frac{1}{N}}\leq |(\frac{\mu(\rho)}{\sqrt{\rho}})^{'}|\leq C\rho^{\frac{\nu_{2}}{2N}-\frac{1}{2}-\frac{1}{N}}\;\;\forall \rho>1.
\end{aligned}
\label{incru}
\end{equation}
Assume that $\frac{\nu_{2}}{2N}-\frac{1}{2}-\frac{1}{N}\leq 0$ which is equivalent to $\nu_{2}\leq N+2$. In this case we obtain:
$$ |(\frac{\mu(\rho)}{\sqrt{\rho}})^{'}|\geq C\;\;\forall \rho\leq 1.$$
Let us recall that the derivative of the inverse function of $g(\rho)=\frac{\mu(\rho)}{\sqrt{\rho}}$ is:
$$(g^{-1})^{'}(\rho)=\frac{1}{g^{'}(g^{-1}(\rho))}.$$
In particular it means that the inverse function  $g^{-1}(\rho)$ is Lipschitz on the region $\rho\leq 1$ and more generally on the region $\rho\leq M$, it provides then the following inequality for $C_{M}>0$ depending on $M$ (since when $\rho\leq M$ we have that $\frac{\mu(\rho)}{\sqrt{\rho}}$ is also bounded via (\ref{13}) and the hypothesis $\frac{1}{2}+\frac{\nu_{1}}{2N}-\frac{1}{N}\geq 0$):
\begin{equation}
\begin{aligned}
|\frac{\mu(\rho)}{\sqrt{\rho}}-\frac{\mu(\rho_{n})}{\sqrt{\rho_{n}}}|1_{\{\rho\leq 1\}\cup \{\rho_{n}\leq 1\}}C_{M}&\geq |g^{-1}(\frac{\mu(\rho)}{\sqrt{\rho}})-g^{-1}(\frac{\mu(\rho_{n})}{\sqrt{\rho_{n}}})|1_{\{\rho\leq M\}\cup \{\rho_{n}\leq M\}}\\
&\geq |\rho-\rho_{n}|1_{\{\rho\leq 1\}\cup \{\rho_{n}\leq 1\}}.
\end{aligned}
\label{butile}
\end{equation}
We deduce the following estimate for any compact $K$:
\begin{equation}
\begin{aligned}
\int_{K}|\rho_{n}(t,x)dx-\rho(t,x)|dx=&\int_{K}|\rho_{n}(t,x)dx-\rho(t,x)|1_{\{\rho>M\}\cup \{\rho_{n}>M\}}(t,x) dx\\
&+\int_{K}|\rho_{n}(t,x)-\rho(t,x)|1_{\{\rho\leq M\}\cup \{\rho_{n}\leq M\}}(t,x) dx
\end{aligned}
\end{equation}
The second term converges uniformly on $(0,T)$ to $0$ when $n$ goes to infinity via (\ref{butile}) applied to $M>0$ and the strong convergence in $C([0,T],L^{2}_{loc}(\R^{N}))$ of $\frac{\mu(\rho_{n})}{\sqrt{\rho_{n}}}$ to $\frac{\mu(\rho)}{\sqrt{\rho}}$. Let us deal with the first term, since we know via the lemma \ref{tech} that $\rho_{n}$ is uniformly bounded in $L^{\infty}_{T}(L^{1+\e}_{loc}(\R^{N}))$ with $\e>0$ we have by H\"older's inequality and Tchebytchev lemma for $C>0$:
$$
\begin{aligned}
&\int_{K}|\rho_{n}(t,x)dx-\rho(t,x)|1_{\{\rho>M\}\cup \{\rho_{n}>M\}}(t,x) dx\leq \\
&(\int_{K}|\rho_{n}(t,x)dx-\rho(t,x)|^{1+\e}1_{\{\rho>M\}\cup \{\rho_{n}>M\}}(t,x) dx)^{\frac{1}{1+\e}}
|\{\rho>1\}\cup \{\rho_{n}>1\}|^{\frac{\e}{1+\e}},\\
&\leq 2C\frac{ \|\rho_{0}\|_{L^{1}(\R^{N})}}{M}.
\end{aligned}
$$
This last term goes uniformly on $(0,T)$ to $0$ in $n$ when $M$ goes to infinity. It show the desired result.
\\
\\
\textbf{Second case $0<\nu_{1}<2$}\\
\\
In this case it suffices to apply exactly the same argument than in the case $\nu_{1}\geq N+2$ except at the place we consider not $\frac{\mu(\rho)}{\sqrt{\rho}}$ but $\beta(\rho)$. It concludes the proof of the lemma. {\hfill $\Box$}
\begin{lem}
Let $\psi\in C^{\infty}_{0}(\R^{N})$ with $\psi=1$ on $B(0,1)$ and $\mbox{supp}\,\psi$ is included in $B(0,2)$. We are going distinguish two cases:\\
\\
$\bullet$ When $\nu_{1}\geq 2$ we set:
$$v_{n}=\psi(\rho_{n})\mu(\rho_{n})u_{n}+(1-\psi(\rho_{n}))\rho_{n}u_{n},$$
we have:
\begin{itemize}
\item $v_{n}$ is uniformly bounded in $L^{2}(0,T,W^{1,1}(K)$ for any compact $K$.
\item $\p_{t}v_{n}$ is uniformly bounded in $L^{2}(0,T,W^{-2,1}(K)$ for any compact $K$.
\end{itemize}
Up to a subsequence, the sequel $v_{n}$ converges strongly in $L^{2}(0,T;L^{p}_{loc}(\R^{N}))$ to some $v(t,x)$ for
all $p\in[1,\frac{3}{2})$. In particular:
$$(\psi(\rho_{n})\mu(\rho_{n})+(1-\psi(\rho_{n})))u_{n}\rightarrow v\;\;\;\mbox{almost everywhere}\;\;(x,t)\in(0,T)\times\R^{N}.$$
Note that we can already define $u(t,x)=\frac{v(t,x)}{ \psi(\rho)\mu(\rho)+(1-\psi(\rho))\rho        }$ outside the vacuum set $\{\rho(t,x)=0\}$,
but we do not know yet whether $v(t,x)$ is zero on the vacuum set (in particular if there is no concentration phenomena for $v$ on $\{\rho(t,x)=0\}$).\\
\\
$\bullet$ When $0<\nu_{1}<2$ we consider:
$$v_{n}=\psi(\rho_{n})\rho_{n}^{\beta_{1}+1}u_{n}+(1-\psi(\rho_{n}))\rho_{n}^{\beta+1}u_{n},$$
with $\beta\leq \frac{-1}{N}$ and $\beta_{1}$ verifying the following assumptions:
\begin{equation}
\begin{aligned}
&\beta_{1}\geq 1,\\
&\beta_{1}+\frac{1}{N}-\frac{\nu_{2}}{2N}\geq 0,\\
&\beta_{1}-\frac{1}{2}+\frac{1}{2N}-\frac{\nu_{2}}{4N}\geq 0
\end{aligned}
\label{trestech}
\end{equation}
we have:
\begin{itemize}
\item $v_{n}$ is uniformly bounded in $L^{2}(0,T,W^{1,1}(K)$ for any compact $K$.
\item $\p_{t}v_{n}$ is uniformly bounded in $L^{2}(0,T,W^{-2,1}(K)$ for any compact $K$.
\end{itemize}
Up to a subsequence, the sequel $v_{n}$ converges strongly in $L^{2}(0,T;L^{p}_{loc}(\R^{N}))$ to some $v(t,x)$ for
all $p\in[1,\frac{3}{2})$. In particular:
$$(\psi(\rho_{n})\rho_{n}^{\beta_{1}+1}+(1-\psi(\rho_{n}))\rho_{n}^{\beta+1})u_{n}\rightarrow v\;\;\;\mbox{almost everywhere}\;\;(x,t)\in(0,T)\times\R^{N}.$$
Note that we can already define $u(t,x)=\frac{v(t,x)}{\psi(\rho)\rho^{\beta_{1}+1}+(1-\psi(\rho))\rho^{\beta+1}      }$ outside the vacuum set $\{\rho(t,x)=0\}$,
but we do not know yet whether $v(t,x)$ is zero on the vacuum set (in particular if there is no concentration phenomena for $v$ on $\{\rho(t,x)=0\}$).
\label{tech2}
\end{lem}
{\bf Proof:} Let us start with proving that $v_{n}$ is uniformly bounded in $L^{2}(0,T,W^{1,1}(\R^{N})$.\\
\\ %In the sequel we shall work with:
%$$v_{n}^{1}=\psi(\rho_{n})\rho_{n}^{\alpha-1}\rho_{n}u_{n}\;\;\mbox{and}\;\;v_{n}^{2}=(1-\psi(\rho_{n}))\rho_{n}u_{n}.$$
$\bullet$ We are going to start with the case $\nu_{1}\geq 2$.\\%\sqrt{\frac{\mu(\rho_{n})}{\rho_{n}}}\rho_{n}u_{n}.$$
Let us prove that $v_{n}^{2}$ is bounded in $L^{2}(0,T,W^{1,1}(K)$ for any compact $K$. We have then:
\begin{equation}
\begin{aligned}
&\p_{i}((1-\psi(\rho_{n}))\rho_{n}u_{nj})=-\frac{\psi^{'}(\rho_{n})}{\mu^{'}(\rho_{n})}   
\rho_{n}\; \frac{\mu^{'}(\rho_{n})}{\sqrt{\rho_{n}}}\p_{i}\rho_{n}\sqrt{\rho_{n}}u_{nj}\\
&+\frac{(1-\psi(\rho_{n}))}{ \mu^{'}(\rho_{n})} \; \frac{\mu^{'}(\rho_{n})}{\sqrt{\rho_{n}}} \p_{i}\rho_{n}\sqrt{\rho_{n}}u_{nj}+\frac{(1-\psi(\rho_{n})) \rho_{n}}{\sqrt{\mu(\rho_{n})}}\sqrt{\mu(\rho_{n})}\p_{i}u_{nj}.
\end{aligned}
\label{t0}
\end{equation}
We observe that by using (\ref{13}), (\ref{BD}) and $\nu_{1}\geq 2$ that $\frac{(1-\psi(\rho_{n}))}{ \mu^{'}(\rho_{n})}$ is uniformly bounded in $L^{\infty}_{T}(L^{\infty}(\R^{N}))$. It implies that the second term $\frac{(1-\psi(\rho_{n}))}{ \mu^{'}(\rho_{n})} \; \frac{\mu^{'}(\rho_{n})}{\sqrt{\rho_{n}}} \p_{i}\rho_{n}\sqrt{\rho_{n}}u_{nj}$ is bounded by H\'older's inequality in $L^{2}_{T}(L^{1}(\R^{N}))$. The first term on the right hand side of (\ref{t0}) is bounded in $L^{2}_{T}(L^{1}(\R^{N}))$ by using the fact that the support of $\psi^{'}$ is included in the shell $C(0,1,2)$. The last term is also bounded in $L^{2}_{T}(L^{1}(\R^{N}))$. Indeed by (\ref{13}) and $\nu_{2}\geq 2$ it suffices to observe that:
$$\frac{(1-\psi(\rho_{n})) \rho_{n}}{\sqrt{\mu(\rho_{n})}}\leq C\sqrt{\rho_{n}},$$
which implies that $\frac{(1-\psi(\rho_{n})) \rho_{n}}{\sqrt{\mu(\rho_{n})}}$ is bounded in $L^{\infty}_{T}(L^{2}(\R^{N}))$.\\
Finally we have seen that $\n v_{n}^{2}$ is bounded in $L^{2}_{T}(L^{1}(\R^{N}))$. And since we see easily that $v^{2}_{n}$ is also bounded in $L^{\infty}_{T}(L^{1}(\R^{N}))$, it implies that $v^{2}_{n}$ is bounded in $L^{2}_{T}(W^{1,1}(\R^{N}))$.\\
Let us deal now with estimating $\p_{t} v_{n}^{2}$. A simple calculus gives for any regular function $g$:
\begin{equation}
\begin{aligned}
&\p_{t}(g(\rho_{n})\rho_{n}u_{n})=-g(\rho_{n})\rho_{n}u_{n}\cdot\n u_{n}-g^{'}(\rho_{n})\rho_{n}^{2}u_{n}{\rm div}u_{n}-u_{n}{\rm div}(g(\rho_{n})\rho_{n}u_{n})\\
&+2{\rm div}(g(\rho_{n})\mu(\rho_{n})D u_{n})-2\n g(\rho_{n})\cdot \mu(\rho_{n}) D u_{n}+\n(g(\rho_{n})\lambda(\rho_{n}){\rm div} u_{n})\\
&\hspace{8cm}-\n g(\rho_{n})\,\lambda(\rho_{n}){\rm div} u_{n},\\
&=-{\rm div}(g(\rho_{n})\rho_{n}u_{n}\otimes  u_{n})-g^{'}(\rho_{n})\rho_{n}^{2}u_{n}{\rm div}u_{n}+2{\rm div}(g(\rho_{n})\mu(\rho_{n})D u_{n})\\
&-2\n g(\rho_{n})\cdot \mu(\rho_{n}) D u_{n}+\n(g(\rho_{n})\lambda(\rho_{n}){\rm div} u_{n})-\n g(\rho_{n})\,\lambda(\rho_{n}){\rm div} u_{n}.
\end{aligned}
\label{moment1}
\end{equation}
When we apply the previous formula to $g(\rho_{n})=(1-\psi(\rho_{n}))\sqrt{\frac{\mu(\rho_{n})}{\rho_{n}}}$, we have to estimate all the terms on the right hand side of (\ref{moment1}), it comes for $T>0$:
\begin{equation}
\begin{aligned}
&g(\rho_{n})\rho_{n}u_{n}\otimes  u_{n}=(1-\psi(\rho_{n}))\sqrt{\rho_{n}}u_{n}\otimes  \sqrt{\rho_{n}}u_{n}.
\end{aligned}
\label{t1}
\end{equation}
%\begin{equation}
%(1-\psi(\rho_{n}))\sqrt{\frac{\mu(\rho_{n})}{\rho_{n}}}\rho_{n}u_{n}\cdot\n u_{n}=\sqrt{\rho_{n}}u_{n}\cdot(\sqrt{\mu(\rho_{n})}\n u_{n})(1-\psi(\rho_{n})).
%\label{t1}
%\end{equation}
By H\"older's inequality we obtain that  $(1-\psi(\rho_{n}))\sqrt{\rho_{n}}u_{n}\otimes  \sqrt{\rho_{n}}u_{n}$ is bounded in $L^{\infty}_{T}(L^{1}(\R^{N}))$. Next we have:
\begin{equation}
g^{'}(\rho_{n})\rho_{n}^{2}u_{n}{\rm div}u_{n}=-\psi^{'}(\rho_{n})\frac{\rho_{n}^{\frac{3}{2}} }{\sqrt{\mu(\rho_{n})}}    \sqrt{\rho_{n}}u_{n}\sqrt{\mu(\rho_{n})}{\rm div}u_{n}
\label{t2}
\end{equation}
%\begin{equation}
%\begin{aligned}
%&( (1-\psi(\rho_{n}))\sqrt{\frac{\mu(\rho_{n})}{\rho_{n}}} )^{'}\rho_{n}^{2}u_{n}{\rm div}u_{n}=\\
%&=\frac{1}{2}(1-\psi(\rho_{n}))\big(\frac{\mu(\rho_{n})}{\rho_{n}} \big)^{-\frac{1}{2}}(\mu^{'}(\rho_{n})\rho_{n}-\mu(\rho_{n}))u_{n}\,{\rm div}u_{n}-\psi^{'}(\rho_{n})\sqrt{\frac{\mu(\rho_{n})}{\rho_{n}}}\rho_{n}^{2}u_{n}\cdot\n u_{n},\\
%&=\frac{1}{2}(1-\psi(\rho_{n}))\big(\frac{\mu^{'}(\rho_{n})\rho_{n}}{\mu(\rho_{n})}-1)\sqrt{\mu(\rho_{n})}{\rm div}u_{n}\sqrt{\rho_{n}}u_{n}-\psi^{'}(\rho_{n})\rho_{n}\sqrt{\rho_{n}}u_{n}\cdot \sqrt{\mu(\rho_{n})}\n u_{n},
%\end{aligned}
%\label{t2}
%\end{equation}
This term is bounded in $L^{2}_{T}(L^{1}(\R^{N}))$ by H\"older's inequality since $\psi^{'}$ is supported in the shell $C(0,1,2)$ which implies that $\psi^{'}(\rho_{n})\frac{\rho_{n}^{\frac{3}{2}} }{\sqrt{\mu(\rho_{n})}}$ is bounded in $L^{\infty}_{T}(L^{\infty}(\R^{N}))$. Similarly we have:
\begin{equation}
\begin{aligned}
&g(\rho_{n})\mu(\rho_{n})D u_{n}=(1-\psi(\rho_{n}))\sqrt{\mu(\rho_{n})}\sqrt{\mu(\rho_{n})}D u_{n},\\
&=(1-\psi(\rho_{n}))\sqrt{\frac{\mu(\rho_{n})}{\sqrt{\rho_{n}}}} \rho_{n}^{\frac{1}{4}}\sqrt{\mu(\rho_{n})}D u_{n}
\end{aligned}
\label{t3}
\end{equation}
According lemma \ref{tech} we deduce easily that $(1-\psi(\rho_{n}))\sqrt{\mu(\rho_{n})}\sqrt{\mu(\rho_{n})}D u_{n}$ is bounded in $L^{2}_{T}(L^{1}(K))$ for any compact $K$. Next we have:
\begin{equation}
2\n g(\rho_{n})\cdot \mu(\rho_{n}) D u_{n}=-2\psi^{'}(\rho_{n})\n\rho_{n}\cdot \mu(\rho_{n}) D u_{n}.
\label{t4}
\end{equation}
This term is bounded in $L^{2}_{T}(L^{1}(\R^{N}))$ by H\"older's inequality since $\psi^{'}$ is supported in the shell $C(0,1,2)$. The two last term in (\ref{moment1}) are similar to treat.
(\ref{t1}), (\ref{t2}), (\ref{t3}) and (\ref{t4}) implies that $\p_{t}v^{2}_{n}$ is uniformly bounded in $L^{2}_{T}(W^{-1,1}(K))$ for any compact $K$.\\
\\
$\bullet$ Case $0<\nu_{1}<2$.\\
\\
In this case we are going to consider:
$$v^{2}_{n}=(1-\psi(\rho_{n}))\rho_{n}^{\beta}\rho_{n}u_{n}.$$
Let us start with proving that $\n v_{n}$ belongs in $L^{2}_{T}(L^{1}(K))$ for any compact $K$. We have then:
\begin{equation}
\begin{aligned}
&\p_{i}((1-\psi(\rho_{n}))\rho_{n}^{\beta}\rho_{n}u_{nj})=-\frac{\psi^{'}(\rho^{\beta+1}_{n})}{\mu^{'}(\rho_{n})}   
\rho_{n}\; \frac{\mu^{'}(\rho_{n})}{\sqrt{\rho_{n}}}\p_{i}\rho_{n}\sqrt{\rho_{n}}u_{nj}\\
&+\frac{(1-\psi(\rho_{n}))\rho_{n}^{\beta}}{ \mu^{'}(\rho_{n})} \; \frac{\mu^{'}(\rho_{n})}{\sqrt{\rho_{n}}} \p_{i}\rho_{n}\sqrt{\rho_{n}}u_{nj}+\frac{(1-\psi(\rho_{n})) \rho^{\beta+1}_{n}}{\sqrt{\mu(\rho_{n})}}\sqrt{\mu(\rho_{n})}\p_{i}u_{nj}.
\end{aligned}
\label{t0c}
\end{equation}
By using (\ref{11}) and (\ref{BD}), we obtain that:
$$
\begin{aligned}
\mu^{'}(\rho)&=\frac{1}{2\rho}\frac{2\mu(\rho)+N\lambda(\rho)}{N}+(1-\frac{1}{N})\frac{\mu(\rho)}{\rho},\\
&\geq (1-\frac{1}{N}+\frac{\nu_{1}}{2N})\frac{\mu(\rho)}{\rho}.
\end{aligned}
$$
We deduce that according to (\ref{13}) that:
$$
\begin{aligned}
|\frac{(1-\psi(\rho_{n}))\rho_{n}^{\beta}}{ \mu^{'}(\rho_{n})}|&\leq C|\frac{(1-\psi(\rho_{n}))\rho_{n}^{\beta+1}}{ \mu(\rho_{n})}|,\\
&\leq C |\frac{(1-\psi(\rho_{n}))\rho^{\beta+1}_{n}}{ \rho^{1-\frac{1}{N}+\frac{\nu_{1}}{2N}}_{n}}|,\\
&\leq C |(1-\psi(\rho_{n}))\rho^{\beta+\frac{1}{N}-\frac{\nu_{1}}{2N}}_{n}|.
\end{aligned}
$$
Since $\beta\leq \frac{-1}{N}$ we deduce that $\frac{(1-\psi(\rho_{n}))\rho_{n}^{\beta}}{ \mu^{'}(\rho_{n})}$ is uniformly bounded in $L^{\infty}_{T}(L^{\infty}(\R^{N}))$. In particular it implies that $\frac{(1-\psi(\rho_{n}))\rho_{n}^{\beta}}{ \mu^{'}(\rho_{n})} \; \frac{\mu^{'}(\rho_{n})}{\sqrt{\rho_{n}}} \p_{i}\rho_{n}\sqrt{\rho_{n}}u_{nj}$ is uniformly bounded in $L^{\infty}_{T}(L^{1}(\R^{N}))$. Since the support of $\psi^{'}$ is included in the shell $C(0,1,2)$ we easily observe that the first term on the right hand side of (\ref{t0c}) is uniformly bounded in $L^{\infty}_{T}(L^{1}(\R^{N}))$. The last term is also bounded in $L^{2}_{T}(L^{1}(\R^{N}))$ since $frac{(1-\psi(\rho_{n})) \rho^{\beta+1}_{n}}{\sqrt{\mu(\rho_{n})}}$ belongs in $L^{\infty}_{T}(L^{2}(\R^{N}))$ by using (\ref{13}), $\beta\leq -\frac{1}{N}$ and the lemma \ref{tech}.\\
Finally we have seen that $\n v_{n}^{2}$ is bounded in $L^{2}_{T}(L^{1}(\R^{N}))$. And since we see easily that $v^{2}_{n}$ is also bounded in $L^{\infty}_{T}(L^{1}(\R^{N}))$, it implies that $v^{2}_{n}$ is bounded in $L^{2}_{T}(W^{1,1}(\R^{N}))$.\\
\\
Let us deal now with estimating $\p_{t} v_{n}^{2}$. It suffices to deal with the formula (\ref{moment1}) and replacing $g(\rho_{n})$ by $(1-\psi(\rho_{n}))\rho_{n}^{\beta}$. Let us start with the first term of (\ref{moment1}):
\begin{equation}
\begin{aligned}
&g(\rho_{n})\rho_{n}u_{n}\otimes  u_{n}=(1-\psi(\rho_{n}))\rho_{n}^{\beta}\sqrt{\rho_{n}}u_{n}\otimes  \sqrt{\rho_{n}}u_{n}.
\end{aligned}
\label{t1c}
\end{equation}
%\begin{equation}
%(1-\psi(\rho_{n}))\sqrt{\frac{\mu(\rho_{n})}{\rho_{n}}}\rho_{n}u_{n}\cdot\n u_{n}=\sqrt{\rho_{n}}u_{n}\cdot(\sqrt{\mu(\rho_{n})}\n u_{n})(1-\psi(\rho_{n})).
%\label{t1}
%\end{equation}
Since $\beta\leq-\frac{1}{N}$ it implies that $(1-\psi(\rho_{n}))\rho_{n}^{\beta}$ is $L^{\infty}$ bounded and then by H\"older's inequality we obtain that  $(1-\psi(\rho_{n}))\rho_{n}^{\beta}\sqrt{\rho_{n}}u_{n}\otimes  \sqrt{\rho_{n}}u_{n}$ is bounded in $L^{\infty}_{T}(L^{1}(\R^{N}))$. Next we have:
\begin{equation}
\begin{aligned}
&g^{'}(\rho_{n})\rho_{n}^{2}u_{n}{\rm div}u_{n}=-\psi^{'}(\rho_{n})\frac{\rho_{n}^{\frac{3}{2}+\beta} }{\sqrt{\mu(\rho_{n})}}    \sqrt{\rho_{n}}u_{n}\sqrt{\mu(\rho_{n})}{\rm div}u_{n},\\
&\hspace{5cm}+\beta(1-\psi(\rho_{n}) )   \frac{\rho_{n}^{\frac{1}{2}+\beta} }{\sqrt{\mu(\rho_{n})}}   \sqrt{\rho_{n}}u_{n}\sqrt{\mu(\rho_{n})}{\rm div}u_{n}.
\end{aligned}
\label{t2c}
\end{equation}
%\begin{equation}
%\begin{aligned}
%&( (1-\psi(\rho_{n}))\sqrt{\frac{\mu(\rho_{n})}{\rho_{n}}} )^{'}\rho_{n}^{2}u_{n}{\rm div}u_{n}=\\
%&=\frac{1}{2}(1-\psi(\rho_{n}))\big(\frac{\mu(\rho_{n})}{\rho_{n}} \big)^{-\frac{1}{2}}(\mu^{'}(\rho_{n})\rho_{n}-\mu(\rho_{n}))u_{n}\,{\rm div}u_{n}-\psi^{'}(\rho_{n})\sqrt{\frac{\mu(\rho_{n})}{\rho_{n}}}\rho_{n}^{2}u_{n}\cdot\n u_{n},\\
%&=\frac{1}{2}(1-\psi(\rho_{n}))\big(\frac{\mu^{'}(\rho_{n})\rho_{n}}{\mu(\rho_{n})}-1)\sqrt{\mu(\rho_{n})}{\rm div}u_{n}\sqrt{\rho_{n}}u_{n}-\psi^{'}(\rho_{n})\rho_{n}\sqrt{\rho_{n}}u_{n}\cdot \sqrt{\mu(\rho_{n})}\n u_{n},
%\end{aligned}
%\label{t2}
%\end{equation}
The first term is bounded in $L^{2}_{T}(L^{1}(\R^{N}))$ by H\"older's inequality since $\psi^{'}$ is supported in the shell $C(0,1,2)$ which implies that $\psi^{'}(\rho_{n})\frac{\rho_{n}^{\frac{3}{2}+\beta} }{\sqrt{\mu(\rho_{n})}} $ is bounded in $L^{\infty}_{T}(L^{\infty}(\R^{N}))$. The second term is also bounded in $L^{2}_{T}(L^{1}(\R^{N}))$ since $(1-\psi(\rho_{n}) )   \frac{\rho_{n}^{\frac{1}{2}+\beta} }{\sqrt{\mu(\rho_{n})}}$ is bounded in $L^{\infty}_{T}(L^{\infty}(\R^{N}))$ by using (\ref{13}) and the fact that $\beta\leq\frac{-1}{N}$.\\
Similarly we have:
\begin{equation}
\begin{aligned}
&g(\rho_{n})\mu(\rho_{n})D u_{n}=(1-\psi(\rho_{n}))\sqrt{\mu(\rho_{n})}\rho_{n}^{\beta}\sqrt{\mu(\rho_{n})}D u_{n},\\
&=(1-\psi(\rho_{n}))\sqrt{\frac{\mu(\rho_{n})}{\sqrt{\rho_{n}}}} \rho_{n}^{\frac{1}{4}+\beta}\sqrt{\mu(\rho_{n})}D u_{n}
\end{aligned}
\label{t3c}
\end{equation}
According lemma \ref{tech} we deduce easily that $(1-\psi(\rho_{n}))\sqrt{\mu(\rho_{n})}\sqrt{\mu(\rho_{n})}D u_{n}$ is bounded in $L^{2}_{T}(L^{1}(K))$ for any compact $K$. The last term gives:
\begin{equation}
\begin{aligned}
&2\n g(\rho_{n})\cdot \mu(\rho_{n}) D u_{n}=-2\psi^{'}(\rho_{n})\rho_{n}^{\beta}\n\rho_{n}\cdot \mu(\rho_{n}) D u_{n}\\
&\hspace{3cm}+\beta(1-\psi(\rho_{n}))\frac{\rho_{n}^{\beta-\frac{1}{2}}\sqrt{\mu(\rho_{n})} }{\mu^{'}(\rho_{n})}  \frac{\mu^{'}(\rho_{n})}{\sqrt{\rho_{n}}}  \n\rho_{n}\cdot \sqrt{\mu(\rho_{n})} D u_{n}. 
\end{aligned}
\label{t4c}
\end{equation}
The first term is bounded in $L^{2}_{T}(L^{1}(\R^{N}))$ by H\"older's inequality since $\psi^{'}$ is supported in the shell $C(0,1,2)$ and we proceed similarly for the second by observing that  $(1-\psi(\rho_{n}))\frac{\rho_{n}^{\beta-\frac{1}{2}}\sqrt{\mu(\rho_{n})}}{\mu^{'}(\rho_{n})}  $ is $L^{\infty}$ bounded via (\ref{BD}), (\ref{11}) and (\ref{13}). The two last term in (\ref{moment1}) are similar to treat.
(\ref{t1c}), (\ref{t2c}), (\ref{t3c}) and (\ref{t4c}) implies that $\p_{t}v^{2}_{n}$ is uniformly bounded in $L^{2}_{T}(W^{-1,1}(K))$ for any compact $K$.
\\
\\
Let us deal now with the term $v_{n}^{1}$ . We are going to distinguish two cases.\\
\\
$\bullet$ Case $\nu_{1}\geq 2$.\\
\\
We have:%start with the case $N=3$. We have:
$$\psi(\rho_{n})\mu(\rho_{n})u_{n}=\frac{\psi(\rho_{n})\mu(\rho_{n})}{\sqrt{\rho_{n}}}\sqrt{\rho_{n}}u_{n},$$
where via the lemma \ref{tech} $\frac{\psi(\rho_{n})\mu(\rho_{n})}{\sqrt{\rho_{n}}}$ is uniformly bounded in $L^{\infty}_{T}(L^{\infty}(\R^{N}))$ and $\sqrt{\rho_{n}}u_{n}$ is uniformly bounded in $L^{\infty}(0,T;L^{2}(\R^{N}))$ which implies that $\psi(\rho_{n})\mu(\rho_{n})u_{n}$ is bounded in  $L^{\infty}(0,T;L^{2}(\R^{N}))$.\\ % for any compact $K$. ON PEUT FAIRE MIEUX PAR INTERPOLATION EN DESSOUSÉ..
%for $q\in[2,+\infty[$ if $N=2$ and $q\in[2,6]$ if $N=3$.
%Since $\sqrt{\rho_{n}}u_{n}$ is bounded in $L^{\infty}(0,T;L^{2}(\R^{N}))$, we deduce that:
%$$\rho_{n}u_{n}\;\;\mbox{is bounded in}\;\;L^{\infty}(0,T;L^{q}(\R^{N}))\;\;\mbox{for all}\;\;q\in[1,\frac{3}{2}].$$
Next we have:
$$
\begin{aligned}
\p_{i}(\psi(\rho_{n})\mu(\rho_{n})u_{nj})&=\psi(\rho_{n})\sqrt{\mu(\rho_{n})}\sqrt{\mu(\rho_{n})}\p_{i}u_{nj}+\psi(\rho_{n})\mu^{'}\rho_{n})\p_{i}\rho_{n}u_{nj},\\
&=\psi(\rho_{n})\sqrt{\mu(\rho_{n})}\sqrt{\mu(\rho_{n})}\p_{i}u_{nj}+\frac{\psi(\rho_{n})\mu^{'}(\rho_{n})}{\sqrt{\rho_{n}}}\p_{i}\rho_{n}\sqrt{\rho_{n}}u_{nj}.
%\rho_{n}\p_{i}u_{nj}+u_{nj}\p_{i}\rho_{n}\\
%&=\sqrt{\rho_{n}}\sqrt{\rho_{n}}\p_{i}u_{nj}+2\sqrt{\rho_{n}}u_{nj}\p_{i}\sqrt{\rho_{n}}.\\
\end{aligned}
$$
By entropy inequality (\ref{22}) we know that $\frac{\psi(\rho_{n})\mu^{'}(\rho_{n})}{\sqrt{\rho_{n}}}\p_{i}\rho_{n}$ is uniformly bounded in $L^{\infty}_{T}(L^{2}(\R^{N}))$ which implies that $\frac{\psi(\rho_{n})\mu^{'}(\rho_{n})}{\sqrt{\rho_{n}}}\p_{i}\rho_{n}\sqrt{\rho_{n}}u_{nj}$ is uniformly bounded in $L^{\infty}_{T}(L^{1}(\R^{N}))$.\\
Let us deal now with the term $\psi(\rho_{n})\sqrt{\mu(\rho_{n})}\sqrt{\mu(\rho_{n})}\p_{i}u_{nj}$, we know that $\sqrt{\mu(\rho_{n})}\p_{i}u_{nj}$ is uniformly bounded in $L^{2}_{T}(L^{2}(\R^{N}))$. Next we know that $\psi(\rho_{n})\sqrt{\mu(\rho_{n})}$ is bounded in $L^{\infty}_{T}(L^{\infty}(\R^{N}))$ which provides on $\psi(\rho_{n})\sqrt{\mu(\rho_{n})}\sqrt{\mu(\rho_{n})}\p_{i}u_{nj}$ a $L^{2}_{T}(L^{2}(\R^{N}))$ bound.
%ave:
%$$\sqrt{\mu(\rho_{n})}=\sqrt{\frac{\mu(\rho_{n})}{\sqrt{\rho_{n}}}}\rho_{n}^{\frac{1}{4}}.$$
%By lemma \ref{tech} it implies that $\sqrt{\frac{\mu(\rho_{n})}{\sqrt{\rho_{n}}}}$  is uniformly bounded in $L^{\infty}_{T}(L^{12}(K))$ for any compact $K$ and $\rho_{n}^{\frac{1}{4}}$ in $L^{\infty}_{T}(L^{4}(\R^{N}))$. We deduce finally that  $\sqrt{\mu(\rho_{n})}\p_{i}u_{nj}$ is uniformly bounded in $L^{2}_{T}(L^{\frac{6}{5}}(K))$ for any compact $K$.
%Using lemma \ref{} and \ref{}, it is readily seen that the second term is bounded in $L^{\infty}(0,T;L^{1}(\R^{N}))$, while the first
%term is bounded in $L^{2}(0,T;L^{q}(\R^{N}))$ for all $q\in[1,\frac{3}{2}]$. 
Hence for any compact $K$:
$$\n(\psi(\rho_{n})\mu(\rho_{n})u_{n})\;\;\mbox{is bounded in}\;\;L^{2}(0,T;L^{1}(K))$$
In particular we have obtained that for all compact $K$:
$$v_{n}^{1}=\psi(\rho_{n})\mu(\rho_{n})u_{n}\;\;\mbox{is bounded in}\;\;L^{2}(0,T;W^{1,1}(K)).$$
We are now going to estimate $\p_{t}v_{n}^{1}$, it suffices to estimate each term on the right hand side of (\ref{moment1}) by replacing $g(\rho_{n})$ by $\psi(\rho_{n})\frac{\mu(\rho_{n})}{\rho_{n}}$.\\
We have then:
\begin{equation}
\begin{aligned}
&g(\rho_{n})\rho_{n}u_{n}\otimes  u_{n}=\psi(\rho_{n})\frac{\mu(\rho_{n})}{\rho_{n}}  \sqrt{\rho_{n}}u_{n}\otimes  \sqrt{\rho_{n}}u_{n},
\end{aligned}
\label{t11}
\end{equation}
%\begin{equation}
%(1-\psi(\rho_{n}))\sqrt{\frac{\mu(\rho_{n})}{\rho_{n}}}\rho_{n}u_{n}\cdot\n u_{n}=\sqrt{\rho_{n}}u_{n}\cdot(\sqrt{\mu(\rho_{n})}\n u_{n})(1-\psi(\rho_{n})).
%\label{t1}
%\end{equation}
By H\"older's inequality we obtain that  $\psi(\rho_{n})\mu(\rho_{n})u_{n}\otimes  u_{n}$ is bounded in $L^{\infty}_{T}(L^{1}(\R^{N}))$. Indeed we have used the fact that by (\ref{13}) and since $\nu_{2}\geq 2$ then $\psi(\rho_{n})\frac{\mu(\rho_{n})}{\rho_{n}} $ is bounded in $L^{\infty}_{T}(L^{\infty}(\R^{N}))$.\\
Next we have:
\begin{equation}
\begin{aligned}
g^{'}(\rho_{n})\rho_{n}^{2}u_{n}{\rm div}u_{n}&=(\psi^{'}(\rho_{n})\mu(\rho_{n})+\psi(\rho_{n})\mu^{'}(\rho_{n}))\frac{\rho_{n}^{\frac{3}{2}} }{\sqrt{\mu(\rho_{n})}}    \sqrt{\rho_{n}}u_{n}\sqrt{\mu(\rho_{n})}{\rm div}u_{n},\\
&=(\psi^{'}(\rho_{n})\sqrt{\mu(\rho_{n})}\rho_{n}^{\frac{3}{2}}+\psi(\rho_{n})\frac{\mu^{'}(\rho_{n})\rho_{n}^{\frac{3}{2}} }{\sqrt{\mu(\rho_{n})}})    \sqrt{\rho_{n}}u_{n}\sqrt{\mu(\rho_{n})}{\rm div}u_{n}
\end{aligned}
\label{t22}
\end{equation}
%\begin{equation}
%\begin{aligned}
%&( (1-\psi(\rho_{n}))\sqrt{\frac{\mu(\rho_{n})}{\rho_{n}}} )^{'}\rho_{n}^{2}u_{n}{\rm div}u_{n}=\\
%&=\frac{1}{2}(1-\psi(\rho_{n}))\big(\frac{\mu(\rho_{n})}{\rho_{n}} \big)^{-\frac{1}{2}}(\mu^{'}(\rho_{n})\rho_{n}-\mu(\rho_{n}))u_{n}\,{\rm div}u_{n}-\psi^{'}(\rho_{n})\sqrt{\frac{\mu(\rho_{n})}{\rho_{n}}}\rho_{n}^{2}u_{n}\cdot\n u_{n},\\
%&=\frac{1}{2}(1-\psi(\rho_{n}))\big(\frac{\mu^{'}(\rho_{n})\rho_{n}}{\mu(\rho_{n})}-1)\sqrt{\mu(\rho_{n})}{\rm div}u_{n}\sqrt{\rho_{n}}u_{n}-\psi^{'}(\rho_{n})\rho_{n}\sqrt{\rho_{n}}u_{n}\cdot \sqrt{\mu(\rho_{n})}\n u_{n},
%\end{aligned}
%\label{t2}
%\end{equation}
The first term is bounded in $L^{2}_{T}(L^{1}(\R^{N}))$ by H\"older's inequality since $\psi^{'}$ is supported in the shell $C(0,1,2)$ which implies that $\psi^{'}(\rho_{n})\sqrt{\mu(\rho_{n})}\rho_{n}^{\frac{3}{2}}$ is bounded in $L^{\infty}_{T}(L^{\infty}(\R^{N}))$. The second term is also bounded in $L^{2}_{T}(L^{1}(\R^{N}))$ because using (\ref{13}) we observe that $\psi(\rho_{n})\frac{\mu^{'}(\rho_{n})\rho_{n}^{\frac{3}{2}} }{\sqrt{\mu(\rho_{n})}}$ is bounded in $L^{\infty}_{T}(L^{\infty}(\R^{N}))$.\\
The third term gives:
\begin{equation}
\begin{aligned}
&g(\rho_{n})\mu(\rho_{n})D u_{n}=\psi(\rho_{n})\frac{\mu^{\frac{3}{2}}(\rho_{n})}{\rho_{n}}\sqrt{\mu(\rho_{n})}D u_{n},\\
&=\psi(\rho_{n}))\frac{\mu(\rho_{n})}{\sqrt{\rho_{n}}}\sqrt{ \frac{\mu(\rho_{n})}{\rho_{n}}} \sqrt{\mu(\rho_{n})}D u_{n}
\end{aligned}
\label{t31}
\end{equation}
By (\ref{13}) and $\nu_{1}\geq 2$  we know that $\sqrt{ \frac{\mu(\rho_{n})}{\rho_{n}}}$ is bounded in  $L^{\infty}_{T}(L^{\infty}(\R^{N}))$, since via lemma \ref{tech} $\frac{\mu(\rho_{n})}{\sqrt{\rho_{n}}}$ is bounded in $L^{\infty}_{T}(L^{6}(K))$ for any compact $K$ we deduce that $\psi(\rho_{n})\frac{\mu(\rho_{n})}{\rho_{n}}\mu(\rho_{n})D u_{n}$ is uniformly bounded in $L^{2}_{T}(L^{1}(K))$ for any compact $K$.\\
%According lemma \ref{tech} we deduce easily that $(1-\psi(\rho_{n}))\sqrt{\mu(\rho_{n})}\sqrt{\mu(\rho_{n})}D u_{n}$ is bounded in $L^{2}_{T}(L^{1}(K))$ for any compact $K$.
Next we have:
\begin{equation}
\begin{aligned}
&2\n g(\rho_{n})\cdot \mu(\rho_{n}) D u_{n}=2\psi^{'}(\rho_{n})\n\rho_{n}\cdot \mu(\rho_{n}) D u_{n}\\
&+2\psi(\rho_{n})\sqrt{\mu(\rho_{n})\rho_{n}}  \frac{\mu^{'}(\rho_{n})}{\sqrt{\rho_{n}}}\n\rho_{n}\cdot \sqrt{\mu(\rho_{n})} D u_{n}.
\end{aligned}
\label{t41}
\end{equation}
The first term is easily bounded in $L^{2}_{T}(L^{1}(\R^{N}))$ since the support of $\psi^{'}$ is included in $C(0,1,2)$. The second term is also bounded in $L^{2}_{T}(L^{1}(\R^{N}))$ by H\"older's inequality because via (\ref{13}) we deduce that $\psi(\rho_{n})\sqrt{\mu(\rho_{n})\rho_{n}}$ is bounded in $L^{\infty}_{T}(L^{\infty}(\R^{N}))$. The two last term in (\ref{moment1}) are similar to treat.
(\ref{t11}), (\ref{t22}), (\ref{t31}) and (\ref{t41}) implies that $\p_{t}v^{1}_{n}$ is uniformly bounded in $L^{2}_{T}(W^{-1,1}(K))$ for any compact $K$.\\
\\
$\bullet$ Case $0<\nu_{1}<2$.\\
\\
We are going to work with:
$$v_{n}^{1}=\psi(\rho_{n})\rho_{n}^{\beta_{1}}\rho_{n}u_{n}.$$
We have then:
$$\psi(\rho_{n})\rho_{n}^{\beta_{1}}\rho_{n}u_{n}=\psi(\rho_{n})\rho_{n}^{\beta_{1}+\frac{1}{2}}\,    \sqrt{\rho_{n}}u_{n}.$$
Since $\beta_{1}\geq \frac{-1}{2}$, it implies that $\psi(\rho_{n})\rho_{n}^{\beta_{1}}\rho_{n}u_{n}$ is bounded in $L^{\infty}_{T}(L^{2}(\R^{N}))$.\\
%where via the lemma \ref{tech} $\frac{\psi(\rho_{n})\mu(\rho_{n})}{\sqrt{\rho_{n}}}$ is uniformly bounded in $L^{\infty}_{T}(L^{\infty}(\R^{N}))$ and $\sqrt{\rho_{n}}u_{n}$ is uniformly bounded in $L^{\infty}(0,T;L^{2}(\R^{N}))$ which implies that $\psi(\rho_{n})\mu(\rho_{n})u_{n}$ is bounded in  $L^{\infty}(0,T;L^{2}(\R^{N}))$.\\ % for any compact $K$. ON PEUT FAIRE MIEUX PAR INTERPOLATION EN DESSOUSÉ..
%for $q\in[2,+\infty[$ if $N=2$ and $q\in[2,6]$ if $N=3$.
%Since $\sqrt{\rho_{n}}u_{n}$ is bounded in $L^{\infty}(0,T;L^{2}(\R^{N}))$, we deduce that:
%$$\rho_{n}u_{n}\;\;\mbox{is bounded in}\;\;L^{\infty}(0,T;L^{q}(\R^{N}))\;\;\mbox{for all}\;\;q\in[1,\frac{3}{2}].$$
Next we have:
\begin{equation}
\begin{aligned}
&\p_{i}(\psi(\rho_{n})\rho_{n}^{\beta_{1}+1}u_{nj})=\frac{\psi(\rho_{n})\rho_{n}^{\beta_{1}+1}}{  \sqrt{\mu(\rho_{n})} } \sqrt{\mu(\rho_{n})}\p_{i}u_{nj}
,\\
&\hspace{2cm}+\beta_{1}\psi(\rho_{n})\frac{\rho_{n}^{\beta_{1}}}{ \mu^{'}(\rho_{n})}\frac{\mu^{'}(\rho_{n})}{\sqrt{\rho_{n}}}\p_{i}\rho_{n} \sqrt{\rho_{n}}u_{nj}+\psi^{'}(\rho_{n})\rho_{n}^{\beta_{1}+1}\p_{i}\rho_{n} u_{nj}.
%&=\psi(\rho_{n})\sqrt{\mu(\rho_{n})}\sqrt{\mu(\rho_{n})}\p_{i}u_{nj}+\frac{\psi(\rho_{n})\mu^{'}(\rho_{n})}{\sqrt{\rho_{n}}}\p_{i}\rho_{n}\sqrt{\rho_{n}}u_{nj}.
%\rho_{n}\p_{i}u_{nj}+u_{nj}\p_{i}\rho_{n}\\
%&=\sqrt{\rho_{n}}\sqrt{\rho_{n}}\p_{i}u_{nj}+2\sqrt{\rho_{n}}u_{nj}\p_{i}\sqrt{\rho_{n}}.\\
\end{aligned}
\label{imipo}
\end{equation}
Using (\ref{13}), (\ref{BD}) and (\ref{11}) we show that:
$$|\psi(\rho_{n})\frac{\rho_{n}^{\beta_{1}}}{ \mu^{'}(\rho_{n})}|\leq \psi(\rho_{n})\rho_{n}^{\beta_{1}+\frac{1}{N}-\frac{\nu_{2}}{2N}}.$$
Since  $\beta_{1}+\frac{1}{N}-\frac{\nu_{2}}{2N}\geq 0$ it implies that $\psi(\rho_{n})\frac{\rho_{n}^{\beta_{1}}}{ \mu^{'}(\rho_{n})}$ is bounded in 
$L^{\infty}_{T}(L^{\infty}(\R^{N}))$. We deduce that the second term on the right hand side of (\ref{imipo}) is uniformly bounded in $L^{\infty}_{T}(L^{1}(\R^{N}))$ . The third term is easy to treat by using that the support of $\psi^{'}$ is a shell $C(0,1,2)$.\\
Let us deal with the first term of (\ref{imipo}). By using (\ref{13}) we get:
$$|\frac{\psi(\rho_{n})\rho_{n}^{\beta_{1}+1}}{  \sqrt{\mu(\rho_{n})} }|\leq \psi(\rho_{n})\rho_{n}^{\beta_{1}+\frac{1}{2}+\frac{1}{2N}-\frac{\nu_{2}}{4N}}.$$
Since $\beta_{1}+\frac{1}{2}+\frac{1}{2N}- \frac{\nu_{2}}{4N}\geq 0$ it provides a $L^{\infty}$ bound on $\frac{\psi(\rho_{n})\rho_{n}^{\beta_{1}+1}}{  \sqrt{\mu(\rho_{n})} }$ which insures the $L^{2}_{T}(L^{2}(\R^{N})$ bound of $\frac{\psi(\rho_{n})\rho_{n}^{\beta_{1}+1}}{  \sqrt{\mu(\rho_{n})} } \sqrt{\mu(\rho_{n})}\p_{i}u_{nj}$.
%By entropy inequality (\ref{22}) we know that $\frac{\psi(\rho_{n})\mu^{'}(\rho_{n})}{\sqrt{\rho_{n}}}\p_{i}\rho_{n}$ is uniformly bounded in $L^{\infty}_{T}(L^{2}(\R^{N}))$ which implies that $\frac{\psi(\rho_{n})\mu^{'}(\rho_{n})}{\sqrt{\rho_{n}}}\p_{i}\rho_{n}\sqrt{\rho_{n}}u_{nj}$ is uniformly bounded in $L^{\infty}_{T}(L^{1}(\R^{N}))$.\\
%Let us deal now with the term $\psi(\rho_{n})\sqrt{\mu(\rho_{n})}\sqrt{\mu(\rho_{n})}\p_{i}u_{nj}$, we know that $\sqrt{\mu(\rho_{n})}\p_{i}u_{nj}$ is uniformly bounded in $L^{2}_{T}(L^{2}(\R^{N}))$. Next we know that $\psi(\rho_{n})\sqrt{\mu(\rho_{n})}$ is bounded in $L^{\infty}_{T}(L^{\infty}(\R^{N}))$ which provides on $\psi(\rho_{n})\sqrt{\mu(\rho_{n})}\sqrt{\mu(\rho_{n})}\p_{i}u_{nj}$ a $L^{2}_{T}(L^{2}(\R^{N}))$ bound.
%ave:
%$$\sqrt{\mu(\rho_{n})}=\sqrt{\frac{\mu(\rho_{n})}{\sqrt{\rho_{n}}}}\rho_{n}^{\frac{1}{4}}.$$
%By lemma \ref{tech} it implies that $\sqrt{\frac{\mu(\rho_{n})}{\sqrt{\rho_{n}}}}$  is uniformly bounded in $L^{\infty}_{T}(L^{12}(K))$ for any compact $K$ and $\rho_{n}^{\frac{1}{4}}$ in $L^{\infty}_{T}(L^{4}(\R^{N}))$. We deduce finally that  $\sqrt{\mu(\rho_{n})}\p_{i}u_{nj}$ is uniformly bounded in $L^{2}_{T}(L^{\frac{6}{5}}(K))$ for any compact $K$.
%Using lemma \ref{} and \ref{}, it is readily seen that the second term is bounded in $L^{\infty}(0,T;L^{1}(\R^{N}))$, while the first
%term is bounded in $L^{2}(0,T;L^{q}(\R^{N}))$ for all $q\in[1,\frac{3}{2}]$. 
Hence for any compact $K$:
$$\n(\psi(\rho_{n})\rho_{n}^{\beta_{1}+1}u_{n})\;\;\mbox{is bounded in}\;\;L^{2}(0,T;L^{1}(K))$$
In particular we have obtained that for all compact $K$:
$$v_{n}^{1}=\psi(\rho_{n})\rho_{n}^{\beta_{1}+1}u_{n}\;\;\mbox{is bounded in}\;\;L^{2}(0,T;W^{1,1}(K)).$$
We are now going to estimate $\p_{t}v_{n}^{1}$, it suffices to estimate each term on the right hand side of (\ref{moment1}) by replacing $g(\rho_{n})$ by $\psi(\rho_{n})\rho_{n}^{\beta_{1}}$.\\
We have then:
\begin{equation}
\begin{aligned}
&g(\rho_{n})\rho_{n}u_{n}\otimes  u_{n}=\psi(\rho_{n})\rho_{n}^{\beta_{1}-1}  \sqrt{\rho_{n}}u_{n}\otimes  \sqrt{\rho_{n}}u_{n},
\end{aligned}
\label{t11d}
\end{equation}
%\begin{equation}
%(1-\psi(\rho_{n}))\sqrt{\frac{\mu(\rho_{n})}{\rho_{n}}}\rho_{n}u_{n}\cdot\n u_{n}=\sqrt{\rho_{n}}u_{n}\cdot(\sqrt{\mu(\rho_{n})}\n u_{n})(1-\psi(\rho_{n})).
%\label{t1}
%\end{equation}
By H\"older's inequality and the fact that $\beta_{1}\geq 1$ we obtain that  $\psi(\rho_{n})\rho_{n}^{\beta_{1}}u_{n}\otimes  u_{n}$ is bounded in $L^{\infty}_{T}(L^{1}(\R^{N}))$. Next we have:
\begin{equation}
\begin{aligned}
g^{'}(\rho_{n})\rho_{n}^{2}u_{n}{\rm div}u_{n}&=(\psi^{'}(\rho_{n})\rho_{n}^{\beta_{1}}+\beta_{1}\psi(\rho_{n})\rho_{n}^{\beta_{1}-1})\frac{\rho_{n}^{\frac{3}{2}} }{\sqrt{\mu(\rho_{n})}}    \sqrt{\rho_{n}}u_{n}\sqrt{\mu(\rho_{n})}{\rm div}u_{n},\\
&=(\psi^{'}(\rho_{n})\frac{\rho_{n}^{\beta_{1}+\frac{3}{2}} }{\sqrt{\mu(\rho_{n})}}+\beta_{1}\psi(\rho_{n})\frac{\rho_{n}^{\beta_{1}+\frac{1}{2}} }{\sqrt{\mu(\rho_{n})}})    \sqrt{\rho_{n}}u_{n}\sqrt{\mu(\rho_{n})}{\rm div}u_{n}
\end{aligned}
\label{t22d}
\end{equation}
%\begin{equation}
%\begin{aligned}
%&( (1-\psi(\rho_{n}))\sqrt{\frac{\mu(\rho_{n})}{\rho_{n}}} )^{'}\rho_{n}^{2}u_{n}{\rm div}u_{n}=\\
%&=\frac{1}{2}(1-\psi(\rho_{n}))\big(\frac{\mu(\rho_{n})}{\rho_{n}} \big)^{-\frac{1}{2}}(\mu^{'}(\rho_{n})\rho_{n}-\mu(\rho_{n}))u_{n}\,{\rm div}u_{n}-\psi^{'}(\rho_{n})\sqrt{\frac{\mu(\rho_{n})}{\rho_{n}}}\rho_{n}^{2}u_{n}\cdot\n u_{n},\\
%&=\frac{1}{2}(1-\psi(\rho_{n}))\big(\frac{\mu^{'}(\rho_{n})\rho_{n}}{\mu(\rho_{n})}-1)\sqrt{\mu(\rho_{n})}{\rm div}u_{n}\sqrt{\rho_{n}}u_{n}-\psi^{'}(\rho_{n})\rho_{n}\sqrt{\rho_{n}}u_{n}\cdot \sqrt{\mu(\rho_{n})}\n u_{n},
%\end{aligned}
%\label{t2}
%\end{equation}
The first term is bounded in $L^{2}_{T}(L^{1}(\R^{N}))$ by H\"older's inequality since $\psi^{'}$ is supported in the shell $C(0,1,2)$ which show that $\psi^{'}(\rho_{n})\frac{\rho_{n}^{\beta_{1}+\frac{3}{2}} }{\sqrt{\mu(\rho_{n})}}$ is bounded in $L^{\infty}_{T}(L^{\infty}(\R^{N}))$. The second term is also bounded in $L^{2}_{T}(L^{1}(\R^{N}))$ because using (\ref{13}) we observe that:
$$|\psi(\rho_{n})\frac{\rho_{n}^{\beta_{1}+\frac{1}{2}} }{\sqrt{\mu(\rho_{n})}}|\leq C \psi(\rho_{n})\rho_{n}^{\beta_{1}+\frac{1}{2N}-\frac{\nu_{2}}{4N}}.$$
It provides a $L^{\infty}$ bounds on $\psi(\rho_{n})\frac{\rho_{n}^{\beta_{1}+\frac{1}{2}} }{\sqrt{\mu(\rho_{n})}}$ since $\beta_{1}+\frac{1}{2N}-\frac{\nu_{2}}{4N}\geq 0$.
\\
The third term gives:
\begin{equation}
\begin{aligned}
&g(\rho_{n})\mu(\rho_{n})D u_{n}=\psi(\rho_{n}) \rho_{n}^{\beta_{1}}\sqrt{\mu(\rho_{n})}\sqrt{\mu(\rho_{n})}D u_{n},\\
%&=\psi(\rho_{n}))\frac{\mu(\rho_{n})}{\sqrt{\rho_{n}}}\sqrt{ \frac{\mu(\rho_{n})}{\rho_{n}}} \sqrt{\mu(\rho_{n})}D u_{n}
\end{aligned}
\label{t31d}
\end{equation}
By (\ref{13}) we prove that $\psi(\rho_{n}) \rho_{n}^{\beta_{1}}\sqrt{\mu(\rho_{n})}$ is bounded in $L^{\infty}_{T}(L^{\infty}(\R^{N})$ and it yields that $\psi(\rho_{n}) \rho_{n}^{\beta_{1}}\sqrt{\mu(\rho_{n})}\sqrt{\mu(\rho_{n})}D u_{n}$ is uniformly bounded in $L^{2}_{T}(L^{2}(\R^{N}))$. Next we have:
\begin{equation}
\begin{aligned}
&2\n g(\rho_{n})\cdot \mu(\rho_{n}) D u_{n}=2\psi^{'}(\rho_{n}) \rho_{n}^{\beta_{1}}\n\rho_{n}\cdot \mu(\rho_{n}) D u_{n}\\
&+2\psi(\rho_{n}) \frac{\rho_{n}^{\beta_{1}-\frac{1}{2}}    \sqrt{\mu(\rho_{n})}}{\mu^{'}(\rho_{n})}  \frac{\mu^{'}(\rho_{n})}{\sqrt{\rho_{n}}}\n\rho_{n}\cdot \sqrt{\mu(\rho_{n})} D u_{n}.
\end{aligned}
\label{t41d}
\end{equation}
The first term is easily bounded in $L^{2}_{T}(L^{1}(\R^{N}))$ since the support of $\psi^{'}$ is included in $C(0,1,2)$. In order to deal with the second term we observe that via (\ref{BD}), (\ref{11}) and (\ref{13}):
$$
\begin{aligned}
|\psi(\rho_{n}) \frac{\rho_{n}^{\beta_{1}-\frac{1}{2}}    \sqrt{\mu(\rho_{n})}}{\mu^{'}(\rho_{n})}|&=|\psi(\rho_{n}) \frac{\rho_{n}^{\beta_{1}-\frac{1}{2}}}{\sqrt{\mu^{'}(\rho_{n})}}    \sqrt{\frac{\mu(\rho_{n})} {\mu^{'}(\rho_{n})} }|,\\
&\leq C \psi(\rho_{n}) \rho_{n}^{\beta_{1}-\frac{1}{2}+\frac{1}{2N}-\frac{\nu_{2}}{4N}}.
\end{aligned}
$$
It implies that $\psi(\rho_{n}) \frac{\rho_{n}^{\beta_{1}-\frac{1}{2}}    \sqrt{\mu(\rho_{n})}}{\mu^{'}(\rho_{n})}$ is bounded in $L^{\infty}_{T}(L^{\infty}(\R^{N}))$ since
$\beta_{1}-\frac{1}{2}+\frac{1}{2N}-\frac{\nu_{2}}{4N}\geq 0$ and we deduce that the second term is bounded in  $L^{2}_{T}(L^{1}(\R^{N}))$. The two last term in (\ref{moment1}) are similar to treat.\\
(\ref{t11d}), (\ref{t22d}), (\ref{t31d}) and (\ref{t41d}) implies that $\p_{t}v^{2}_{n}$ is uniformly bounded in $L^{2}_{T}(W^{-1,1}(K))$ for any compact $K$.
\null{\hfill $\Box$}
\subsubsection*{Step 2: Convergence of $\sqrt{\rho_{n}}u_{n}$ and $\rho_{n}u_{n}$}
In the sequel we shall define $h(\rho)$ as follows:
\begin{equation}
\begin{aligned}
&h(\rho)=\psi(\rho)\frac{\mu(\rho)}{\sqrt{\rho}}+(1-\psi(\rho))\sqrt{\rho}\;\;\;\;\mbox{if}\;\;\nu_{1}\geq 2,\\
&h(\rho)=\psi(\rho)\rho^{\beta_{1}+\frac{1}{2}}+(1-\psi(\rho))\rho^{\beta+\frac{1}{2}}\;\;\;\;\mbox{if}\;\;0<\nu_{1}< 2.
\end{aligned}
\label{defh}
\end{equation}
Here $\psi$, $\beta$ and $\beta_{1}$ verify the hypothesis of lemma \ref{tech2}.
\begin{lemme}
The quantity $\sqrt{\rho_{n}}u_{n}$ strongly  converges in $L^{2}_{loc}((0,T)\times\Omega)$ to $\frac{v}{h(\rho)}$ (which is null when $v=0$).\\
In particular, we have $v(t,x)=0$ a.e on $\{\rho(t,x)=0\}$ and there exists a function $u(t,x)$ such that $v(t,x)=\sqrt{\rho}(t,x)h(\rho)(t,x)u(t,x)$ and:
$$
\begin{aligned}
&\sqrt{\rho_{n}}u_{n}\rightarrow\sqrt{\rho}u\;\;\;\mbox{strongly in}\;\;L^{2}_{loc}((0,T)\times\R^{N}),\\
&\rho_{n}u_{n}\rightarrow\rho u\;\;\;\mbox{strongly in}\;\;L^{1}(0,T;L^{1}_{loc}(\R^{N})).%,\;p\in[1,\frac{3}{2}).
\end{aligned}
$$
\label{ruse2}
\end{lemme}
\begin{remarka}
Here $u$ is not uniquely defined on the vacuum set $\{\rho(t,x)=0\}$. We will set $u=0$ on the vacuum set $\{\rho(t,x)=0\}$.
\end{remarka}
{\bf Proof:} Since $\frac{v_{n}}{h(\rho_{n})}$ is uniformly bounded in $L^{\infty}(0,T;L^{2}(\R^{N}))$, Fatou's lemma implies that:
$$\int\lim\inf\frac{v_{n}^{2}}{h(\rho_{n})^{2}}dx<+\infty.$$
We deduce that $v(t,x)=0$ a.e. in $\{\rho(t,x)=0\}$ since $h(\rho)=0$ when $\rho=0$. We can define the limit velocity by $u(t,x)$ with $u(t,x)=\frac{v(t,x)}{\sqrt{\rho}(t,x)h(\rho(t,x))}$ when $\rho(t,x)\ne0$ and $u(t,x)=0$ on  $\{\rho(t,x)=0\}$. In particular this last point implies that there is no concentration effect of $\rho_{n}u_{n}\otimes u_{n}$ on the set $\{\rho=0\}$. And for all $t>0$:
$$\int_{\R^{N}}\frac{v^{2}(t,x)}{h^{2}(\rho)(t,x)}dx=\int_{\R^{N}}\rho(t,x)|u(t,x)|^{2}dx<+\infty.$$
Furthermore applying the Fatou's lemma once more, we obtain:
$$
\begin{aligned}
\int\rho|u|^{2}\ln(1+|u|^{2})dx&\leq\int\lim\inf\rho_{n}|u_{n}|^{2}\ln(1+|u_{n}|^{2})dx\\
&\leq\lim\inf\int\rho_{n}|u_{n}|^{2}\ln(1+|u_{n}|^{2})dx,
\end{aligned}
$$
which yields $\rho|u|^{2}\ln(1+|u|^{2})\in L^{\infty}(0,T;L^{1}(\R^{N}))$.\\
Next, since $v_{n}$ and $\rho_{n}$ converge almost everywhere, we know that in $\{\rho(t,x)\ne0\}$, $\sqrt{\rho_{n}}u_{n}=\frac{v_{n}}
{h(\rho_{n})}$ converges almost everywhere to $\sqrt{\rho}u=\frac{v}{h(\rho)}$. In particular it implies that:
\begin{equation}
\begin{aligned}
&\sqrt{\rho_{n}}u_{n}1_{\{\{|u_{n}|\leq M\}\cap\{\rho>0\}}\rightarrow\sqrt{\rho}u1_{\{|u|\leq M\}}\;\;\;\mbox{almost everywhere.}\hspace{3cm}\\
&\sqrt{\rho_{n}}u_{n}1_{\{\{|u_{n}|\leq M\}\cap  \{\rho=0\}}\leq M\sqrt{\rho_{n}}\rightarrow0\;\;\;\mbox{almost everywhere.}
\end{aligned}
\end{equation}
%In $\{\rho(t,x)=0\}$, we have
%$\sqrt{\rho_{n}}u_{n}1_{\{|u_{n}|\leq M\}}\leq M\sqrt{\rho_{n}}\rightarrow0$.\\
%Moreover, we have:
%$$\sqrt{\rho_{n}}u_{n}1_{\{|u_{n}|\leq M\}}\rightarrow\sqrt{\rho}u1_{\{|u|\leq M\}}\;\;\;\mbox{almost everywhere}.$$
%As a matter of fact, the convergence holds almost everywhere in $\{\rho(t,x)\ne0\}$, and i
Following the argument of the proof of the lemma \ref{lemmeimp} for any compact $K$ we have:
\begin{equation}
\begin{aligned}
&\int_{(0,T)\times K}|\sqrt{\rho_{n}}u_{n}-\sqrt{\rho}u|^{2}dxdt\leq\int_{(0,T)\times K}|\sqrt{\rho_{n}}u_{n}1_{\{|u_{n}|\leq M\}}-\sqrt{\rho}u1_{\{|u|\leq M\}}|^{2}dxdt\\[2mm]
&\hspace{3cm}+2\int_{(0,T)\times K}|\sqrt{\rho_{n}}u_{n}1_{\{|u_{n}|\geq M\}}|^{2}dxdt+2\int_{(0,T)\times K}|\sqrt{\rho}u1_{\{|u|\geq M\}}|^{2}dxdt,\\
\end{aligned}
\label{cle}
\end{equation}
Let us deal with the first term of (\ref{cle}) we have then:
\begin{equation}
\begin{aligned}
&\int_{(0,T)\times K}|\sqrt{\rho_{n}}u_{n}1_{\{|u_{n}|\leq M\}}-\sqrt{\rho}u1_{\{|u|\leq M\}}|^{2}dxdt\leq \\
&\hspace{2cm}\int_{\{(0,T)\times K\}\cap\{\sqrt{\rho_{n}}\leq c\}}|\sqrt{\rho_{n}}u_{n}1_{\{|u_{n}|\leq M\}}-\sqrt{\rho}u1_{\{|u|\leq M\}}|^{2}dxdt\\
&\hspace{2cm}+\int_{\{(0,T)\times K\}\cap\{\sqrt{\rho_{n}}> c\}}|\sqrt{\rho_{n}}u_{n}1_{\{|u_{n}|\leq M\}}-\sqrt{\rho}u1_{\{|u|\leq M\}}|^{2}dxdt
\end{aligned}
\label{cle1}
\end{equation}
The first term of (\ref{cle1}) converges to $0$ when $n$ goes to $+\infty$ via the theorem of dominated convergence.\\
Now let us recall that via the inequality (\ref{13}) we have:
$$\rho_{n}^{1-\frac{1}{N}-\frac{1}{2}+\frac{\nu}{N}}\leq \frac{\mu(\rho_{n})}{\sqrt{\rho_{n}}}\;\;\mbox{when}\;\;\rho_{n}\geq1.$$
Since $\frac{\mu(\rho_{n})}{\sqrt{\rho_{n}}}$ is uniformly bounded in $L^{\infty}_{T}(L^{6}(\R^{N}))$ for $N=3$ (see the lemma \ref{tech}) we deduce that $\sqrt{\rho_{n}}$ is uniformly bounded in $L^{\infty}_{T}(L^{2+4\nu}(\R^{N}))$. It allows to deal with the second term of (\ref{cle1}) for $N=3$ by H\"older's inequality and Tchebytchev lemma with $c$ going to infinity.\\
For $N=2$ using the inequality (\ref{13}) and the fact that $\frac{\mu(\rho_{n})}{\sqrt{\rho_{n}}}$ is uniformly bounded in $L^{\infty}_{T}(L^{p}_{loc}(\R^{N}))$ for any $p>1$ when $N=2$ (see the lemma \ref{tech}), we deal similarly with the second term of (\ref{cle}) via H\"older's inequality and Tchebytchev lemma.\\ 
% $\rho$ Since $\sqrt{\rho_{n}}$
%It is obvious that $\sqrt{\rho_{n}}u_{n}1_{\{|u_{n}|\leq M\}}$ is bounded uniformly in $L^{\infty}(0,T;L^{3}(\Omega))$,
%so \ref{} gives the convergence of the first integral:
%$$\int|\sqrt{\rho_{n}}u_{n}1_{\{|u_{n}|\leq M\}}-\sqrt{\rho}u1_{\{|u|\leq M\}}|^{2}dxdt\rightarrow0.$$
Finally, we take advantage of the gain of integrability on the velocity furnished by the entropy (\ref{gain}):
$$\int|\sqrt{\rho_{n}}u_{n}1_{\{|u_{n}|\geq M\}}|^{2}dxdt\leq\frac{1}{\ln(1+M^{2})}\int\rho_{n}|u_{n}|^{2}\ln(1+|u_{n}|^{2})dxdt.$$
Similarly we have:
$$\int|\sqrt{\rho}u1_{\{|u|\geq M\}}|^{2}dxdt\leq\frac{1}{\ln(1+M^{2})}\int\rho|u|^{2}\ln(1+|u|^{2})dxdt.$$
Combining all the previous estimate, it yields:
$$\lim\sup_{n\rightarrow+\infty}\int|\sqrt{\rho_{n}}u_{n}-\sqrt{\rho}u|^{2}dxdt\leq\frac{C}{\ln(1+M^{2})}$$
for all $M>0$, and the lemma follows by taking $M\rightarrow+\infty$.\\
\\
$\bullet$ Let us prove now the strong convergence of $\rho_{n}u_{n}$ to $\rho u$. Since $\sqrt{\rho_{n}}u_{n}$ converges strongly in $L^{2}_{loc}((0,T)\times \R^{N})$ to $\sqrt{\rho}u$ and that via the lemma \ref{rho} $\sqrt{\rho_{n}}$ converges also in  $L^{2}_{loc}((0,T)\times \R^{N})$ to $\sqrt{\rho}$ it implies that $\rho_{n}u_{n}$ converges strongly in  $L^{1}_{loc}((0,T)\times \R^{N})$ to $\rho u$. {\hfill $\Box$}
\subsubsection*{ Step 3: Convergence of the diffusion terms}
\begin{lemme}
\label{ruse3}
We have the convergence in distribution sense up to subsequence for any $T>0$:
$$
\begin{aligned}
&\mu(\rho_{n})\n u_{n}\rightarrow \mu(\rho)\n u\;\;\;\mbox{in}\;\;{\cal D}^{'}((0,T)\times\R^{N}),\\
&\mu(\rho_{n})^{t}\n u_{n}\rightarrow \mu(\rho)^{t}\n u\;\;\;\mbox{in}\;\;{\cal D}^{'}((0,T)\times\R^{N}),\\
\end{aligned}
$$
and:
$$
\begin{aligned}
&\lambda(\rho_{n}){\rm div}u_{n}\rightarrow \lambda(\rho){\rm div}u\;\;\;\mbox{in}\;\;{\cal D}^{'}((0,T)\times\R^{N}).
\end{aligned}
$$
\end{lemme}
{\bf Proof:} Let $\phi$ be a test function, then:
\begin{equation}
\begin{aligned}
&\int \mu(\rho_{n})\n u_{n}\phi dxdt=-\int \mu(\rho_{n})u_{n}\n\phi dxdt+\int u_{n}\n \mu(\rho_{n})\phi \,dxdt\\
&\hspace{2cm}=-\int\frac{\mu(\rho_{n})}{\sqrt{\rho_{n}}}\sqrt{\rho_{n}}u_{n}\n\phi\, dxdt+\int \sqrt{\rho_{n}}u_{n}\frac{\mu^{'}(\rho_{n})}{\sqrt{\rho_{n}}}\n\rho_{n}\phi \,dxdt.\\
&\hspace{2cm}=-\int\frac{\mu(\rho_{n})}{\sqrt{\rho_{n}}}\sqrt{\rho_{n}}u_{n}\n\phi\, dxdt+\frac{1}{2}\int \sqrt{\rho_{n}}u_{n}\n f(\rho_{n})\phi \,dxdt.
\end{aligned}
\label{aimp}
\end{equation}
From lemma \ref{tech} in the appendix we know that $\frac{\mu(\rho_{n})}{\sqrt{\rho_{n}}}$ is uniformly bounded in $L^{\infty}(0,T;L^{6}_{loc}(\R^{N}))$.
Furthermore via the inequality (\ref{13}) and the convergence almost everywhere from $\rho_{n}$ to $\rho$ we know that $\frac{\mu(\rho_{n})}{\sqrt{\rho_{n}}}$ converges almost everywhere to $\frac{\mu(\rho)}{\sqrt{\rho}}$
(defined to be zero on the vacuum set). Therefore by the lemma \ref{lemmeimp}, it converges strongly in $L^{2}_{loc}((0,T)\times\R^{N})$ to  $\frac{\mu(\rho)}{\sqrt{\rho}}$. This point is enough to prove
the convergence of the first term as $\sqrt{\rho_{n}}u_{n}$ converges strongly.\\
Next since $\n f(\rho_{n})$ is bounded in $L^{\infty}(0,T;L^{2}(\R^{N}))$, up to a subsequence 
$\n f(\rho_{n})$ converges weakly to $v$ in $L^{2}_{loc}((0,T)\times\R^{N})$. In addition by Sobolev embedding we know that $f(\rho_{n})$
is bounded in $L^{\infty}(0,T;L^{6}_{loc}(\R^{N})))$. Since $f(\rho_{n})$  converges almost everywhere ($f$ is a continuous function) to $f(\rho)$, it converges strongly
in $L^{2}_{loc}((0,T)\times\R^{N})$ by using the lemma \ref{lemmeimp}. It follows that:
$$\n f(\rho_{n})\rightarrow\n f(\rho)\;\;\;L^{2}_{loc}((0,T)\times\R^{N})-\mbox{weak}.$$
It concludes the proof for the second term of (\ref{aimp}).\\
A similar argument holds for $\mu(\rho_{n})^{t}\n u_{n}$ and $\lambda(\rho_{n}){\rm div}u_{n}$ inasmuch as we have $|\lambda(\rho)|\leq C \mu(\rho)$
and $|\lambda^{'}(\rho)|\leq C \mu^{'}(\rho)$ via (\ref{11}) and the remark \ref{r2.2}.
\subsubsection*{Global existence when $u_{0}=-\n\va(\rho_{0})$ when $\mu(\rho)=\mu\rho^{\alpha}$ with $\alpha>1-\frac{1}{N}$}
The first thing to observe is that $\mu(\rho)=\mu\rho^{\alpha}$ with $\alpha>1-\frac{1}{N}$ verifies the hypothesis of theorem \ref{theo2} for the stability of the global weak solution and in particular (\ref{11}).Now it suffices to construct a sequence of global regular solution $(\rho_{n},u_{n})$ verifying the hypothesis of the theorem \ref{theo2}, in particular the uniform bound via the entropy (\ref{21}), (\ref{22}) and (\ref{gain}) and the following properties:
\begin{itemize}
\item $\rho_{0}^{n}$ converges strongly to $\rho_{0}$ in $L^{1}(\R^{N})$.
\item  $\rho_{0}^{n}u_{0}^{n}$ converges strongly to $\rho_{0}u_{0}$ in $L^{1}(\R^{N})$.
\end{itemize}
Let $\kappa$ a function belonging in the Schwarz space ${\cal S}(\R^{N})$ with $\kappa>0$ and $\int_{\R^{N}}\kappa dx=1$. We define $\kappa_{n}$ by:
$$\kappa_{n}=n^{N}\kappa(n\cdot).$$
Let us take for example $\kappa(x)=(2\pi)^{\frac{-N}{2}}e^{\frac{-|x|^{2}}{2}}$. By using the theorem \ref{theo1} and by setting $\rho_{0}^{n}=\rho_{0}+\frac{f_{0}}{n}$  with $f_{0}$ continuous and belonging in $L^{1}(\R^{N})\cap L^{\infty}(\R^{N})\cap W^{1,3}(\R^{N})\cap  W^{1,1}(\R^{N})$ and strictly positive and $\rho_{0}$ continuous we know that there exists a global regular solution of the system (\ref{3systemev}) with initial data $(\rho_{0}^{n},-\n\va(\rho_{0}^{n}))$. Indeed $ \rho_{0}^{n}$ verifies the hypothesis of theorem \ref{theo2} since  $ \rho_{0}^{n}$ is continuous in $L^{1}(\R^{N})$ and strictly positive. We also observe that $\rho_{0}^{n}$ converges strongly to $\rho_{0}$ in $L^{1}(\R^{N})$. Let us deal now with a more simple case when we assume that $\rho_{0}$ belongs also in $L^{\infty}(\R^{N})\cap W^{1,3}(\R^{N})\cap W^{1,1}(\R^{N})$.\\
Next we have:
\begin{equation}
\begin{aligned}
\sqrt{\rho_{0}^{n}}\va^{'}(\rho_{0}^{n})\n \rho_{0}^{n}&=\alpha (\rho_{0}^{n})^{\alpha-\frac{3}{2}}\n \rho_{0}^{n},\\
&=\alpha (\rho_{0}+\frac{f_{0}}{n})^{\alpha-\frac{3}{2}}\n \rho_{0}+\frac{\alpha}{n} (\rho_{0}+\frac{f_{0}}{n})^{\alpha-\frac{3}{2}}\n f_{0}.
\end{aligned}
\label{initialtech}
\end{equation}
Let us distinguish two cases $\alpha\geq \frac{3}{2}$ and $1-\frac{1}{N}\alpha\leq \frac{3}{2}$.\\
\\
$\bullet$ $\alpha\geq \frac{3}{2}$. In this case since $\rho_{0}$ belongs in $H^{1}(\R^{N})\cap L^{\infty}(\R^{N})$ we deduce that $\n f(\rho_{0}^{n})$ and $\sqrt{\rho_{0}^{n}}u_{0}^{n}$ are uniformly bounded in $L^{\infty}(L^{2}(\R^{N}))$. Indeed we have:
$$
\begin{aligned}
&|\sqrt{\rho_{0}^{n}}\va^{'}(\rho_{0}^{n})\n \rho_{0}^{n}|\leq (\alpha +1)(\|\rho_{0}\|_{L^{\infty}(\R^{N})}+\|f_{0}\|_{L^{\infty}(\R^{N})})^{\alpha-\frac{3}{2}}\big(|\n \rho_{0}|+|\n f_{0}|\big),
\end{aligned}
$$
which implies the previous statement.\\
%Furthermore we have if $\alpha\geq 2$:
%\begin{equation}
%\begin{aligned}
%&|\sqrt{\rho_{0}^{n}}u_{0}^{n}\sqrt{\ln(1+|u_{0}^{n}|^{2})}|=| \sqrt{\rho_{0}^{n}}u_{0}^{n}\sqrt{\ln(1+4( \rho_{0}^{n})^{2(\alpha-2)}|\n \rho_{0}^{n}|^{2})}|,\\[2mm]
%&\leq |\sqrt{\rho_{0}^{n}}u_{0}^{n}|\sqrt{\ln(1+4 (\|\rho_{0}\|_{L^{\infty}(\R^{N})}+\|f\|_{L^{\infty}(\R^{N})})^{2(\alpha-2)}(2|\n \rho_{0}|_{L^{\infty}(\R^{N})}^{2}+2|\n f_{0}|_{L^{\infty}(\R^{N})}^{2}))}.\\[2mm]
%\end{aligned}
%\label{gainini1}
%\end{equation}
%\\
%Since we have assumed that $f_{0}$ and $\rho_{0}$ are in $W^{1,\infty}(\R^{N})$, it proves that $\sqrt{\rho_{0}^{n}}u_{0}^{n}\sqrt{\ln(1+|u_{0}^{n}|^{2})}$ is uniformly bounded in $L^{2}(\R^{N})$.\\
Let us  prove that%tackle the problem when $\frac{3}{2}\leq\alpha<2$, it suffices in this case to check that
 $\sqrt{\rho_{0}^{n}}(u_{0}^{n})^{1+\frac{1}{p}}$ is uniformly bounded in $L^{2}(\R^{N})$ for a $p>1$ large enough. It will be then sufficient in order to show that $\sqrt{\rho_{0}^{n}}u_{0}^{n}\sqrt{\ln(1+|u_{0}^{n}|^{2})}$ is uniformly bounded in $L^{2}(\R^{N})$. We have then for $p$ large enough $(1+\frac{1}{p})(\alpha-2)+\frac{1}{2}\geq 0$ and:
\begin{equation}
\begin{aligned}
&|\sqrt{\rho_{0}^{n}}(u_{0}^{n})^{1+\frac{1}{p}}|\leq C (\rho_{0}^{n})^{(1+\frac{1}{p})(\alpha-2)+\frac{1}{2}}|\n \rho_{0}^{n}|^{1+\frac{1}{p}},\\[2mm]
&\leq  C (\|\rho_{0}\|_{L^{\infty}(\R^{N})}+\|f_{0}\|_{L^{\infty}(\R^{N})})^{(1+\frac{1}{p})(\alpha-2)+\frac{1}{2}}(|\n \rho_{0}|^{1+\frac{1}{p}}+|\n f_{0}|^{1+\frac{1}{p}}).
%&\leq |\sqrt{\rho_{0}^{n}}u_{0}^{n}|\sqrt{\ln(1+4 (\|\rho_{0}\|_{L^{\infty}(\R^{N})}+\|f\|_{L^{\infty}(\R^{N})})^{2(\alpha-2)}(2|\n \rho_{0}|_{L^{\infty}(\R^{N})}^{2}+2|\n f_{0}|_{L^{\infty}(\R^{N})}^{2}))}.\\[2mm]
\end{aligned}
\label{gainini1}
\end{equation}
Since we have assumed that $\rho_{0}$ and $f_{0}$ are in $W^{1,1}(\R^{N})\cap W^{1,3}(\R^{N})$ it shows that for $p$ large enough, $\sqrt{\rho_{0}^{n}}(u_{0}^{n})^{1+\frac{1}{p}}$ is bounded in $L^{2}(\R^{N})$ which implies the desired result.\\
%\begin{equation}
%\begin{aligned}
%&|\rho_{0}^{n}|u_{0}^{n}|^{2}\ln(1+|u_{0}^{n}|^{2})|= \sqrt{\rho_{0}^{n}}u_{0}^{n}\sqrt{\ln(1+4( \rho_{0}^{n})^{2(\alpha-2)}(|\n \rho_{0}|^{2}+\frac{1}{n^{2}}|\n f_{0}|^{2}+\frac{2}{n} \n \rho_{0}\cdot \n f_{0}  )},\\
%&\leq  \sqrt{\rho_{0}^{n}}u_{0}^{n}\sqrt{\ln(1+4[ \rho_{0}^{2(\alpha-2)}|\n \rho_{0}|^{2}+n^{2(1-\alpha)}f_{0}^{2(\alpha-2)}|\n f_{0}|^{2}+2n^{1-\alpha}\rho_{0}^{\alpha-2}f_{0}^{\alpha-2}| \n \rho_{0}\cdot \n f_{0}|  )}
%\end{aligned}
%\label{gainini}
%\end{equation}
Let us finished by proving the strong convergence of $\rho_{0}^{n}u_{0}^{n}$ to $\rho_{0}u_{0}$. We have then:
\begin{equation}
|\rho_{0}^{n}u_{0}^{n}-\rho_{0}u_{0}|=|\alpha[ (\rho_{0}+\frac{f_{0}}{n})^{\alpha-1}-\rho_{0}^{\alpha-1}]\n \rho_{0}+\frac{\alpha}{n} (\rho_{0}+\frac{f_{0}}{n})^{\alpha-1}\n f_{0}|.
\label{initialappro}
\end{equation}
Since $\rho_{0}$ and $f_{0}$ are bounded in $L^{\infty}(\R^{N})$ we deduce that:
$$\|(\rho_{0}+\frac{f_{0}}{n})^{\alpha-1}-\rho_{0}^{\alpha-1}\|_{L^{\infty}(\R^{N})}\rightarrow _{n\rightarrow+\infty}0.$$
This last inequality and H\"older's inequality show in particular that $[ (\rho_{0}+\frac{f_{0}}{n})^{\alpha-1}-(\rho_{0})^{\alpha-1}]\n \rho_{0}$ converges strongly in $L^{1}(\R^{N})$ to $0$. The second term of (\ref{initialappro}) is easy to treat. It gives that $\rho_{0}^{n}u_{0}^{n}$ converges strongly to $\rho_{0}u_{0}$.\\
In particular via the first part of the theorem \ref{theo2}, it implies that $(\rho_{n},u_{n})$ converge in distribution sense to a global weak solution $(\rho,u)$.\\
\\
$\bullet$ $1-\frac{1}{N}<\alpha<\frac{3}{2}$. In this case we shall assume moreover that $\sqrt{\rho_{0}}|\n\va(\rho_{0})|^{1+\frac{1}{p}}$ and  $\sqrt{f_{0}}|\n\va(f_{0})|^{1+\frac{1}{p}}$ are bounded in $L^{2}(\R^{N})$ for $p$ large enough. We have via (\ref{initialtech}):
\begin{equation}
\begin{aligned}
|\sqrt{\rho_{0}^{n}}\va^{'}(\rho_{0}^{n})\n \rho_{0}^{n}|\leq \alpha \rho_{0}^{\alpha-\frac{3}{2}}|\n \rho_{0}|+\alpha n^{\frac{1}{2}-\alpha} f_{0}^{\alpha-\frac{3}{2}}|\n f_{0}|.
\end{aligned}
\label{initialtech1}
\end{equation}
Since we have $\alpha>\frac{1}{2}$ it implies that $\sqrt{\rho_{0}^{n}}\va^{'}(\rho_{0}^{n})\n \rho_{0}^{n}$ is uniformly bounded in $L^{2}(\R^{N})$ by using the fact that $\sqrt{\rho_{0}}\n\va(\rho_{0})$ and  $\sqrt{f_{0}}\n\va(f_{0})$ are bounded in $L^{2}(\R^{N})$.
It implies that $\n f(\rho_{0}^{n})$ and $\sqrt{\rho_{0}^{n}}u_{0}^{n}$ are uniformly bounded in $L^{\infty}(L^{2}(\R^{N}))$. Next by (\ref{gainini1}) we have:
\begin{equation}
\begin{aligned}
|\sqrt{\rho_{0}^{n}}(&u_{0}^{n})^{1+\frac{1}{p}}|\leq C (\rho_{0}^{n})^{(1+\frac{1}{p})(\alpha-2)+\frac{1}{2}}|\n \rho_{0}^{n}|^{1+\frac{1}{p}},\\[2mm]
&\leq C (\rho_{0}^{n})^{(1+\frac{1}{p})(\alpha-2)+\frac{1}{2}}(|\n \rho_{0}|^{1+\frac{1}{p}}+\frac{1}{n^{1+\frac{1}{p}}}|\n f_{0}|^{1+\frac{1}{p}}),\\
&\leq C (\rho_{0}^{(1+\frac{1}{p})(\alpha-2)+\frac{1}{2}}|\n \rho_{0}|^{1+\frac{1}{p}}+n^{(1+\frac{1}{p})(2-\alpha)-\frac{3}{2}-\frac{1}{p}}       f_{0}^{\frac{1}{2}+(1+\frac{1}{p})(\alpha-2)}|\n f_{0}|^{1+\frac{1}{p}},\\
&\leq C (\rho_{0}^{(1+\frac{1}{p})(\alpha-2)+\frac{1}{2}}|\n \rho_{0}|^{1+\frac{1}{p}}+n^{(\frac{1}{2}-\alpha+\frac{1}{p}(1-\alpha)}       f_{0}^{\frac{1}{2}+(1+\frac{1}{p})(\alpha-2)}|\n f_{0}|^{1+\frac{1}{p}}.
%&\leq |\sqrt{\rho_{0}^{n}}u_{0}^{n}|\sqrt{\ln(1+4 (\|\rho_{0}\|_{L^{\infty}(\R^{N})}+\|f\|_{L^{\infty}(\R^{N})})^{2(\alpha-2)}(2|\n \rho_{0}|_{L^{\infty}(\R^{N})}^{2}+2|\n f_{0}|_{L^{\infty}(\R^{N})}^{2}))}.\\[2mm]
\end{aligned}
\label{gainini2}
\end{equation}
Since $\alpha>\frac{1}{2}$ by choosing $p$ large enough we obtain that $n^{(\frac{1}{2}-\alpha+\frac{1}{p}(1-\alpha)}$ which is uniformly bounded. By the fact that  $\sqrt{\rho_{0}}|\n\va(\rho_{0})|^{1+\frac{1}{p}}$ and  $\sqrt{f_{0}}|\n\va(f_{0})|^{1+\frac{1}{p}}$ are bounded in $L^{2}(\R^{N})$, it implies that $\sqrt{\rho_{0}^{n}}(u_{0}^{n})^{1+\frac{1}{p}}$ is bounded in $L^{2}(\R^{N})$  and so that $\rho_{0}^{n}|u_{0}^{n}|^{2}\ln(1+|u_{0}^{n}|^{2})$ is uniformly bounded in $L^{2}(\R^{N})$.\\
Finally by (\ref{initialappro}) we are going to prove that $\rho_{0}^{n}u_{0}^{n}$ converges strongly to $\rho_{0}u_{0}$. We start by remarking that when $\alpha\leq 1$:
$$|[ (\rho_{0}+\frac{f_{0}}{n})^{\alpha-1}-\rho_{0}^{\alpha-1}]\n \rho_{0}|\leq 2 \rho_{0}^{\alpha-1}|\n \rho_{0}|.$$
We deduce by the theorem of dominated convergence that $ [(\rho_{0}+\frac{f_{0}}{n})^{\alpha-1}-\rho_{0}^{\alpha-1}]\n \rho_{0}$ converges strongly to $0$ in $L^{1}(\R^{N})$. The second term on the right hand side of (\ref{initialtech}) goes also trivially to $0$. It the case where $1\leq\alpha<\frac{3}{2}$ it suffices to observe that:
$$\||[ (\rho_{0}+\frac{f_{0}}{n})^{\alpha-1}-\rho_{0}^{\alpha-1}]\|_{L^{\infty}(\R^{N})}\rightarrow_{n\rightarrow+\infty}0.$$
It implies in particular since $\n\rho_{0}$ belongs in $W^{1,1}(\R^{N})$ that $[ (\rho_{0}+\frac{f_{0}}{n})^{\alpha-1}-\rho_{0}^{\alpha-1}]\n \rho_{0}$ converges strongly to $0$ in $L^{1}(\R^{N})$. It achieves the proof of  the strong convergence of $\rho_{0}^{n}u_{0}^{n}$ to $\rho_{0}u_{0}$. Finally it implies that $(\rho_{n},u_{n})$ converge also in the case $1-\frac{1}{N}<\alpha<\frac{3}{2}$ to a global weak solution $(\rho,u)$ of the system (\ref{3systemev}).\\
\\
We have previously proved the existence of a global weak solution $(\rho,u)$  of the system (\ref{3systemev}) by assuming extra conditions on the initial density $\rho_{0}$, typically $\rho_{0}$ belonging in  $L^{\infty}(\R^{N})\cap H^{1}(\R^{N})\cap W^{1,1}(\R^{N})$, $\rho_{0}$ continuous and $\sqrt{\rho_{0}}|\n\va(\rho_{0})|^{1+\frac{1}{p}}$ bounded in $L^{2}(\R^{N})$ for $p$ large enough. Let us deal with the general case, it suffices by a second approximation on the initial data $(\rho_{0}^{n},\sqrt{\rho_{n}}(u_{0}^{n})^{1+\frac{1}{p}})$ with $p$ large enough to pass to the limit by using the first part of the theorem \ref{theo2}. More precisely by a convolution approximation we choose $\rho_{0}^{n}=\rho_{0}*\kappa_{n}$ which belongs in $L^{\infty}(\R^{N})\cap H^{1}(\R^{N})\cap W^{1,1}(\R^{N})$ and is continuous and we set $\sqrt{\rho_{n}}(u_{0}^{n})^{1+\frac{1}{p}}=(\sqrt{\rho}u_{0}^{1+\frac{1}{p}})*\kappa_{n}$.\\
\\
Let us now describe  the form of the solution $(\rho,u)$, in particular we are interested in proving that $\rho$ is also the unique global strong solution  of the system \ref{P}) (for a such result, we refer to the theorem \ref{theo2.4}). To do this we recall that our first approximation $(\rho_{n},u_{n})$ is solution of (\ref{P}) with a initial data strictly positive in $L^{1}(\R^{N})$ and continuous. Let us recall that from theorem \ref{theo2.4} the porous media equation verifies a crucial property which ensures the uniqueness the $L^{1}$ contraction principle. Let us apply this property to the sequence $(\rho_{n})_{n\in\mathbb{N}}$, we have then for all $n,m\in \mathbb{N}$ :
$$\|\rho_{n}(t,\cdot)-\rho_{m}(t,\cdot)\|_{L^{1}(\R^{N})}\leq \|\rho_{0}^{n}-\rho_{0}^{m}\|_{L^{1}(\R^{N})}.$$
Since we know that $\rho_{0}^{n}$ converge to $\rho_{0}$ in $L^{1}(\R^{N})$ it implies that $(\rho_{n})_{n\in\mathbb{N}}$ is a Cauchy sequence in $C([0,T],L^{1}(\R^{N}))$ for any $T>0$ which implies that $(\rho_{n})_{n\in\mathbb{N}}$ converges strongly to $\rho_{1}$ in $C([0,T],L^{1}(\R^{N}))$. But since via the first part of the theorem \ref{theo2} we know that $\rho_{n}$ converges also strongly in $C([0,T],L^{1}_{loc}(\R^{N}))$ to $\rho$, it implies that $\rho=\rho_{1}$. Furthermore by the definition of the $L^{1}$ solution of the porous media equation (indeed these last one are defined as limit of energy global weak solution after a regularization of the initial data $\rho_{0}$ where we use the fundamental $L^{1}$ contraction principle which ensures the uniqueness of a such process, we called this the limit solutions in the $L^{1}$ setting. We refer to the proof of the theorem \ref{theo2.4} where we explained precisely how are defined the $L^{1}$ solutions. For more details on the theory the reader can also consult the chapter 6 and 9 of the excellent book of V\'azquez \cite{Vaz}), we know that $\rho_{1}$ is  the unique solution of porous media equation with $\rho^{1}_{0}=\rho_{0}$. It proves that $\rho=\rho_{1}$ is the global unique solution of (\ref{P}) with initial data $\rho_{0}$ which belongs in $L^{1}(\R^{N})$. We proceed similarly for the second approximation $(\rho_{n},u_{n})$ since each time the approximated sequel verify the $L^{1}$ contraction principle of the porous media equation. It concludes the proof of the theorem \ref{theo2}.
%is bounded in a space better than $L^{\infty}(0,T;L^{2}(\R))$, it will allow to give the convergence of the momentum (step 4) and finally the strong convergence of $\sqrt{\rho_{n}}v_{n}\otimes \sqrt{\rho_{n}}u_{n}$ in $L_{w}^{2}((0,T)\times\R)$ (step 5). In fact by considering the effective velocity $v_{n}$ we do not need extra information on the density $\rho_{n}$ to deal with the capillarity terms (as it is the case in \cite{5BDL,fH2}). Furthermore at the difference with \cite{5BDL} we are able to treat the momentum term without using test functions depending of the density (or in an other way to add a cold pressure).
%The last step shows the convergence of the momentum term $\sqrt{\rho_{n}}v_{n}$ to
%$\sqrt{\rho}u+\frac{2\kappa}{\mu}\n\rho$.
 {\hfill $\Box$}
\subsection{Proof of Theorem \ref{theo3} and corollary \ref{cor}}
\label{section6}
\subsection{Proof of theorem \ref{theo3}}
We are now going to prove that if we have some global weak solution $(\rho_{\ep},u_{\ep})$ for the system (\ref{0.2}) in the sense of the definition \cite{defMV} (or see \cite{MV}), then these global weak solution converge in distribution sense to a quasi-solution $(\rho,u)$ with initial data  $(\rho_{0},u_{0})$. To prove this, it will suffice  to use the same compactness argument  than in the previous section; except that we shall deal with the pressure term and that the entropy (\ref{gaina}) is quite more complicated since there is a reminder term to deal with. Via the entropy (\ref{21a}), (\ref{22a}) we have:
\begin{equation}
\begin{aligned}
&\int_{\R^{N}}\big[\rho_{\ep}|u_{\ep}(t,x)|^{2}(t,x)+\frac{\ep}{\gamma-1}\rho^{\gamma}_{\ep}\big]\,dx+\int^{t}_{0}\int_{\R^{N}}\mu(\rho_{\ep})|D u_{\ep}|^{2}dxdt\\
&\hspace{2cm}+\int^{t}_{0}\int_{\R^{N}}\lambda(\rho_{\ep})|{\rm div} u_{\ep}|^{2}dxdt\leq  \int_{\R^{N}}\big[\rho_{0}|u_{0}|^{2}(x)+\frac{\ep}{\gamma-1}\rho^{\gamma}_{0}\big]dx.\\[2mm]
&\int_{\R^{N}}\big[\rho_{\ep}|u_{\ep}(t,x)|^{2}+\rho_{\ep}|\n\va(\rho_{\ep})|^{2}(t,x)\big]\,dx%+\int^{t}_{0}\int_{\R^{N}}\mu(\rho_{\ep})|\n u_{\ep}|^{2}dxdt+\int^{t}_{0}\int_{\R^{N}}\lambda(\rho_{\ep})|{\rm div} u_{\ep}|^{2}dxdt\\
%&\hspace{1cm}
+\ep\int^{t}_{0}\int_{\R^{N}}\n\va(\rho_{\ep})\cdot\n\rho_{\ep}^{\gamma}dxdt\\
&\hspace{3cm}\leq C( \int_{\R^{N}}\big(\rho_{0}|u_{0}|^{2}(x)+\rho_{0}|\n\va(\rho_{0})|^{2}(x)+\frac{\ep}{\gamma-1}\rho^{\gamma}_{0}(x)\big)\,dx).
\end{aligned}
\label{23}
\end{equation}
\begin{lem}
We are going to distinguish two cases.
\begin{itemize}
\item $\nu_{1}\geq 2$ $\sqrt{\e}\rho_{\e}^{\frac{\gamma}{2}}$ is uniformly bounded in $L^{2}_{T}(L^{6}(\R^{N}))$ for any $T>0$ when $N=3$ and in
$L^{2}_{T}(L^{q}(\R^{N}))$ for any $T>0$  and $q\geq 2$ when $N=2$.
\item $0<\nu_{1}<2$ 
\begin{enumerate}
\item $\e\rho_{\ep}^{\gamma-\frac{1}{N}+\frac{\nu_{1}}{2N}}$ is uniformly bounded in $L^{1}_{T}(L^{3}(\R^{N}))$ for any $T>0$ when $N=3$.
\item  $\e  \rho_{\ep}^{\gamma-\frac{1}{N}+\frac{\nu_{1}}{2N}}$ is bounded in any $L^{1}_{T}(L^{p}(\R^{N}))$ with $p\in[2,+\infty[$ when $N=2$.
\end{enumerate}
\end{itemize}
\label{tech3}
\end{lem}
{\bf Proof:} The inequality (\ref{23}) insures a uniform bound of $\sqrt{\ep}\sqrt{\va^{'}(\rho_{\ep})\rho_{\ep}^{\gamma-1}}\n\rho_{\ep}$ in $L^{2}((0,T),L^{2}(\R^{N}))$ for any $T>0$. Let us evaluate $\sqrt{\va^{'}(\rho_{\ep})\rho_{\ep}^{\gamma-1}}\n\rho_{\ep}$, we have using (\ref{11}) and (\ref{BD}):
$$
\begin{aligned}
\sqrt{\va^{'}(\rho_{\ep})\rho_{\ep}^{\gamma-1}}\n\rho_{\ep}&=\sqrt{\frac{N\lambda(\rho_{\ep})+2\mu(\rho_{\ep})}{N}\rho_{\ep}^{\gamma-3}+2(1-\frac{1}{N})\mu(\rho_{\ep})\rho_{\ep}^{\gamma-3}}\;\;\n\rho_{\ep},\\
&\geq \sqrt{2(1-\frac{1}{N}+\frac{\nu_{1}}{2})\mu(\rho_{\ep})\rho_{\ep}^{\gamma-3}}\;\;\n\rho_{\ep},
\end{aligned}
$$
Via (\ref{13}) we deduce that:
$$
\begin{aligned}
&\rho_{\ep}^{\frac{\gamma}{2}-1-\frac{1}{2N}+\frac{\nu_{1}}{4N}}|\n\rho_{\ep}|\leq C|\sqrt{\va^{'}(\rho_{\ep})\rho_{\ep}^{\gamma-1}}\n\rho_{\ep}|\;\;\forall \rho_{\ep}>1,\\
&\rho_{\ep}^{\frac{\gamma}{2}-1-\frac{1}{2N}+\frac{\nu_{2}}{4N}}|\n\rho_{\ep}|\leq C|\sqrt{\va^{'}(\rho_{\ep})\rho_{\ep}^{\gamma-1}}\n\rho_{\ep}|\;\;\forall \rho_{\ep}\leq 1.
\end{aligned}
$$
Let set $\psi$ a $C^{\infty}_{0}$ function such that $\psi=1$ on $B(0,1)$ and $\psi=0$ on $^{c}B(0,2)$. It implies that since $\mbox{supp}\,\psi^{'}$ is included in the shell $C(0,1,2)$:
\begin{itemize}
\item $\sqrt{\e} \psi(\rho_{\e})\n \rho_{\ep}^{\frac{\gamma}{2}-\frac{1}{2N}+\frac{\nu_{2}}{4N}}$ is uniformly bounded in $L^{2}_{T}(L^{2}(\R^{N}))$ for any $T>0$.
\item $\sqrt{\e}(1-\psi(\rho_{\e}) )\n \rho_{\ep}^{\frac{\gamma}{2}-\frac{1}{2N}+\frac{\nu_{1}}{4N}}$ is uniformly bounded in $L^{2}_{T}(L^{2}(\R^{N}))$ for any $T>0$.
\end{itemize}
Next by (\ref{23}) we know that $\ep^{\frac{1}{\gamma}}\rho_{\ep}$ is uniformly bounded in $L^{\infty}_{T}(L^{\gamma}(\R^{N}))$ for any $T>0$ which implies that $\sqrt{\e}\rho^{\frac{\gamma}{2}}_{\ep}$  is uniformly bounded in $L^{\infty}_{T}(L^{2}(\R^{N}))$ for any $T>0$. Let us deal with two different cases.
\subsubsection*{$\nu_{1}\geq 2$} %ON PEUT FAIRE MIEUX ON A DE LA PERTE
In this case we have : $-\frac{1}{2N}+\frac{\nu_{1}}{4N}\geq 0$ it implies that  $\sqrt{\e}(1-\psi(\rho_{\e}) )\n \rho_{\ep}^{\frac{\gamma}{2}}$ is uniformly bounded in $L^{2}_{T}(L^{2}(\R^{N}))$ for any $T>0$. And since  $\ep^{\frac{1}{2}}\rho_{\ep}^{\frac{\gamma}{2}}$ is uniformly bounded in $L^{\infty}_{T}(L^{2}(\R^{N}))$ for any $T>0$, we deduce that $\ep^{\frac{1}{2}} (1-\psi(\rho_{\e}) )\rho_{\ep}^{\frac{\gamma}{2}}$ is uniformly bounded in $L^{\infty}_{T}(L^{2}(\R^{N}))$ for any $T>0$. We deduce that  $\ep^{\frac{1}{2}} (1-\psi(\rho_{\e}) )\rho_{\ep}^{\frac{\gamma}{2}}$ is uniformly bounded in $L^{2}_{T}(H^{1}(\R^{N}))$ for any $T>0$.  Easily since $\ep^{\frac{1}{2}} \psi(\rho_{\e}) \rho_{\ep}^{\frac{\gamma}{2}}$ is uniformly bounded in $L^{\infty}_{T}(L^{1}(\R^{N}))\cap L^{\infty}_{T}(L^{\infty}(\R^{N}))$, we deduce by Sobolev embedding that $\ep^{\frac{1}{2}}  \rho_{\ep}^{\frac{\gamma}{2}}$
is uniformly bounded in $L^{2}_{T}(L^{6}(\R^{N}))$ when $N=3$ and in $L^{2}_{T}(L^{q}(\R^{N}))$ for any $q\geq 2$.
\subsubsection*{$0<\nu_{1}<2$}
Let us distinguish the case $N=2$ and $N=3$.\\
\\
$\bullet  N=3$\\
We know that $\sqrt{\e}(1-\psi(\rho_{\e}) )\n \rho_{\ep}^{\frac{\gamma}{2}-\frac{1}{2N}+\frac{\nu_{1}}{4N}}$ is uniformly bounded in $L^{2}_{T}(L^{2}(\R^{N}))$ for any $T>0$. Furthermore since $\sqrt{\e}\rho_{\e}^{\frac{\gamma}{2}}$ is uniformly bounded in $L^{\infty}_{T}(L^{1}(\R^{N}))$ we deduce that $\e^{\frac{1}{p}}(1-\psi(\rho_{\e}) ) \rho_{\ep}^{\frac{\gamma}{2}-\frac{1}{2N}+\frac{\nu_{1}}{4N}}$ is uniformly bounded in $L^{\infty}_{T}(L^{p}(\R^{N}))$ with $p(\frac{\gamma}{2}-\frac{1}{2N}+\frac{\nu_{1}}{4N})=2$ with $p>2$ (indeed it is possible because $\frac{\gamma}{2}-\frac{1}{2N}+\frac{\nu_{1}}{4N}>0$) . Indeed we have:
$$\|\e^{\frac{1}{p}}\rho_{\e}^{\frac{\gamma}{2}-\frac{1}{2N}+\frac{\nu_{1}}{4N}}\|_{L^{\infty}_{T}(L^{p}(\R^{N}))}=\|\e^{\frac{1}{\gamma}}\rho_{\e}\|^{\frac{\gamma}{p}}_{L^{\infty}_{T}(L^{\gamma}(\R^{N}))}.$$
Since $p>2$ it implies that $\sqrt{\e}(1-\psi(\rho_{\e}) ) \rho_{\ep}^{\frac{\gamma}{2}-\frac{1}{2N}+\frac{\nu_{1}}{4N}}$ is uniformly bounded in $L^{\infty}_{T}(L^{p}(\R^{N}))$ and also in $L^{\infty}_{T}(L^{2}(\R^{N}))$ because $\sqrt{\e}(1-\psi(\rho_{\e}) ) \rho_{\ep}^{\frac{\gamma}{2}-\frac{1}{2N}+\frac{\nu_{1}}{4N}}$ is strictly positive only a set of finite measure (it is a direct consequence of the Tchebytchev lemma). We have shown that $\sqrt{\e}(1-\psi(\rho_{\e}) ) \rho_{\ep}^{\frac{\gamma}{2}-\frac{1}{2N}+\frac{\nu_{1}}{4N}}$ is uniformly bounded in $L^{2}_{T}(H^{1}(\R^{N}))$ for any $T>0$. By Sobolev embedding we deduce that $\e (1-\psi(\rho_{\e}) )\rho_{\ep}^{\gamma-\frac{1}{N}+\frac{\nu_{1}}{2N}}$ is uniformly bounded in $L^{1}_{T}(L^{3}(\R^{N}))$ for any $T>0$.\\
Since we know that $\psi(\rho_{\e})\rho_{\e}$ is uniformly bounded in $L^{\infty}_{T}(L^{1}(\R^{N}))\cap L^{\infty}_{T}(L^{\infty}(\R^{N}))$ we deduce that $\e\psi(\rho_{\e})\rho_{\ep}^{\gamma-\frac{1}{N}+\frac{\nu_{1}}{2N}}$ is uniformly bounded in $L^{1}_{T}(L^{3}(\R^{N}))$ for any $T>0$ because if $\gamma-\frac{1}{N}+\frac{\nu_{1}}{2N}\geq 1$ this is obvious by interpolation. In the other case $\e\psi(\rho_{\e})\rho_{\ep}^{\gamma-\frac{1}{N}+\frac{\nu_{1}}{2N}}$ is uniformly bounded in $L^{\infty}_{T}(L^{p}(\R^{N}))$ with $p(\gamma-\frac{1}{N}+\frac{\nu_{1}}{2N})=1$. And we have $p=\frac{1}{\gamma-\frac{1}{N}+\frac{\nu_{1}}{2N}}\leq 2$, we conclude also by interpolation in order to prove that $\e\psi(\rho_{\e})\rho_{\ep}^{\gamma-\frac{1}{N}+\frac{\nu_{1}}{2N}}$ is uniformly bounded in $L^{\infty}_{T}(L^{3}(\R^{N}))$.\\ 
Finally we obtain that $\e\rho_{\ep}^{\gamma-\frac{1}{N}+\frac{\nu_{1}}{2N}}$ is uniformly bounded in $L^{1}_{T}(L^{3}(\R^{N}))$ for any $T>0$.\\
\\
$\bullet$  $N=2$\\
% since $\sqrt{\e}\n\big( (1-\psi(\rho_{\e}) )\rho_{\ep}^{\frac{\gamma}{2}-\frac{1}{2N}+\frac{\nu_{1}}{4N}}\big)$ is uniformly bounded in $L^{2}_{T}(L^{2}(\R^{N}))$ for any $T>0$
%and $\sqrt{\ep}\rho_{\ep}^{\frac{\gamma}{2}}$  is uniformly bounded in $L^{\infty}_{T}(L^{2}(\R^{N}))$, we deduce that  
Similarly we obtain that $\sqrt{\e}(1-\psi(\rho_{\e}) )\rho_{\ep}^{\frac{\gamma}{2}-\frac{1}{2N}+\frac{\nu_{1}}{4N}}$ is uniformly bounded in $L^{2}_{T}(H^{1}(\R^{N}))$ for any $T>0$. It provides a uniform bound for $\sqrt{\e}(1-\psi(\rho_{\e}) )\rho_{\ep}^{\frac{\gamma}{2}-\frac{1}{2N}+\frac{\nu_{1}}{4N}}$ in any $L^{2}_{T}(L^{q}(\R^{N}))$ with $q\in[2,+\infty[$. It means that $\e(1-\psi(\rho_{\e}) )\rho_{\ep}^{\gamma-\frac{1}{N}+\frac{\nu_{1}}{2N}}$ is in any $L^{1}_{T}(L^{p}(\R^{N}))$ with $p\in[1,+\infty[$. Since $\psi(\rho_{\e})\rho_{\e}$ is uniformly bounded in $L^{\infty}_{T}(L^{1}(\R^{N}))\cap L^{\infty}_{T}(L^{\infty}(\R^{N}))$ we deduce that $\e \psi(\rho_{\e}) \rho_{\ep}^{\gamma-\frac{1}{N}+\frac{\nu_{1}}{2N}}$ is bounded in any $L^{1}_{T}(L^{p}(\R^{N}))$ with $p\in[1,+\infty[$ if $\gamma-\frac{1}{N}+\frac{\nu_{1}}{2N}\geq 1$ and in any 
$L^{1}_{T}(L^{p}(\R^{N}))$ with $[p_{1},+\infty[$ where $p_{1}(\gamma-\frac{1}{N}+\frac{\nu_{1}}{2N})=1$ in the other case where $p_{1}\leq 2$. It implies that:
\begin{itemize}
\item  $\e  \rho_{\ep}^{\gamma-\frac{1}{N}+\frac{\nu_{1}}{2N}}$ is bounded in any $L^{1}_{T}(L^{p}(\R^{N}))$ with $p\in[1,+\infty[$ if $\gamma-\frac{1}{N}+\frac{\nu_{1}}{2N}\geq 1$.
\item  $\e  \rho_{\ep}^{\gamma-\frac{1}{N}+\frac{\nu_{1}}{2N}}$ is bounded in any $L^{1}_{T}(L^{p}(\R^{N}))$ with $p\in[p_{1},+\infty[$ with $p_{1}(\gamma-\frac{1}{N}+\frac{\nu_{1}}{2N})= 1$ in other case.
\end{itemize}
It achieves the proof of the lemma. \null{\hfill $\Box$}\\
\\
In the following lemma we are going to get uniform estimate on the pressure $\rho_{\ep}^{\gamma}$.
\begin{lemme}
Let us distinguish two cases:
\begin{itemize}
\item  $\nu_{1}\geq 2$. The pressure $\ep\rho_{\ep}^{\gamma}$ is bounded in $L^{\frac{5}{3}}((0,T)\times\R^{N})$ when $N=3$ and $L^{r}((0,T)\times\R^{N})$
for all $r\in[1,2[$ when $N=2$.% In particular, $\rho_{\ep}^{\gamma}$ converges to $\rho^{\gamma}$ strongly in $L^{1}_{loc}((0,T)\times\R^{N})$.
\item $0<\nu_{1}<2$. The pressure $\e\rho_{\e}^{\gamma}$ is bounded in $L^{r_{1}}_{T}(L^{r_{1}}(\R^{N}))$
with $r_{1}=2-\frac{2-\nu_{1}}{6(1+\nu_{1})}$ when $N=3$.
\end{itemize}
\label{lempression}
\end{lemme}
{\bf Proof:}\\
$\bullet$ $\nu_{1}\geq 2$:\\ 
We have seen in the lemma \ref{tech3} that when $N=2$, $\sqrt{\ep} \rho_{\ep}^{\frac{\gamma}{2}}$ is bounded in $L^{2}(0,T;L^{q}(\R^{N}))$ for any $q\geq 2$. We deduce that $\ep \rho_{\ep}^{\gamma}$ is bounded in $L^{1}(0,T;L^{p}(\R^{N}))\cap L^{\infty}(L^{1}(\R^{N}))$ for all $p\in[1,+\infty[$, hence by interpolation $\ep\rho_{\ep}^{\gamma}$
is bounded in $L^{r}((0,T)\times\R^{N})$ for all $r\in[1,2[$.\\
When $N=3$, by Sobolev embedding we only obtain that $\sqrt{\ep}\rho_{\ep}^{\frac{\gamma}{2}}$ bounded in $L^{2}(0,T;L^{6}(\R^{N}))$ which gives that  $\ep \rho_{\ep}^{\gamma}$ is uniformly bounded in $L^{1}(0,T;L^{3}(\R^{N}))$. By H\"older inequality we have:
$$\|\ep\rho_{\ep}^{\gamma}\|_{L^{\frac{5}{3}}((0,T)\times\R^{N}}\leq\|\ep\rho_{\ep}^{\gamma}\|_{L^{\infty}(0,T;L^{1}(\R^{N}))}^{\frac{2}{5}}
\|\ep\rho_{\ep}^{\gamma}\|_{L^{1}(0,T;L^{3}(\R^{N}))}^{\frac{3}{5}}$$
Hence $\ep\rho_{\ep}^{\gamma}$ is bounded in $L^{\frac{5}{3}}((0,T)\times\R^{N})$.\\
\\
$\bullet$ $0<\nu_{1}< 2$:\\
\\
When $N=3$ we know via the lemma \ref{tech3} that $\e\rho_{\e}^{\gamma-\frac{1}{N}+\frac{\nu_{1}}{2N}}$ is bounded in $L^{1}_{T}(L^{3}(\R^{N}))$. We have in particular that:
$$\e\rho_{\e}^{\gamma}=(\e\rho_{\e}^{\gamma-\frac{1}{N}+\frac{\nu_{1}}{2N}})\rho_{\e}^{\frac{1}{N}-\frac{\nu_{1}}{2N}}.$$
Via lemma \ref{tech} we have seen that $\rho_{\e}$ is uniformly bounded in $L^{\infty}_{T}(L^{1+\nu_{1}}(\R^{N}))$ when $N=3$. We define $p$ such that:
$$p(\frac{1}{N}-\frac{\nu_{1}}{2N})=p(\frac{1}{3}-\frac{\nu_{1}}{6})=1+\nu_{1}\;\;\Leftrightarrow\;\; p=\frac{6(1+\nu_{1})}{2-\nu_{1}}\;\;\mbox{with}\;\;\nu_{1}<2.$$
By H\" older's inequality we obtain with $\frac{1}{p_{1}}=\frac{1}{3}+\frac{1}{p}$ that:
$$\|\e\rho_{\e}^{\gamma}\|_{ L^{1}_{T}(L^{p_{1}}(\R^{N}))}    =\|\e\rho_{\e}^{\gamma-\frac{1}{N}+\frac{\nu_{1}}{2N}}\|_{L^{1}_{T}(L^{3}(\R^{N}))}
\|\rho_{\e}^{\frac{1}{N}-\frac{\nu_{1}}{2N}}\|_{ L^{\infty}_{T}(L^{p}(\R^{N}))} .$$
Now since $\e\rho_{\e}^{\gamma}$ is bounded in $L^{1}_{T}(L^{p_{1}}(\R^{N}))\cap L^{\infty}_{T}(L^{1}(\R^{N}))$ we have by interpolation that $\e\rho_{\e}^{\gamma}$ is bounded in $L^{r}_{T}(L^{q}(\R^{N}))$ with:
$$
\begin{cases}
\begin{aligned}
&\frac{1}{r}=\theta,\\
&\frac{1}{q}=\frac{\theta}{p_{1}}+1-\theta,
\end{aligned}
\end{cases}
$$
It implies the following relation $\frac{1}{q}+\frac{1}{r}(1-\frac{1}{p_{1}})=1$. In particular we obtain that $\e\rho_{\e}^{\gamma}$ is bounded in $L^{r_{1}}_{T}(L^{r_{1}}(\R^{N}))$
with $r_{1}=2-\frac{2-\nu_{1}}{6(1+\nu_{1})}$.
{\hfill $\Box$}\\
\\
%Since we already now that $\rho_{n}^{\gamma}$ converges almost everywhere to $\rho^{\gamma}$, those bounds yields the strong convergence of $\rho_{n}^{\gamma}$ in $L^{1}_{loc}((0,T)\times\R^{N})$. 
We are going finally to prove that $\sqrt{\rho_{\ep}}|u_{\ep}|(\ln(1+|u_{\ep}|^{2})^{\frac{1}{2}}$ is uniformly bounded in $L^{\infty}_{T}(L^{2}(\R^{N}))$.
\begin{lemme}
 $\sqrt{\rho_{\ep}}|u_{\ep}|(\ln(1+|u_{\ep}|^{2})^{\frac{1}{2}}$ is uniformly bounded in $L^{\infty}_{T}(L^{2}(\R^{N}))$ for the following situations:
 \begin{itemize}
 \item $\nu_{1}\geq 2$ and:
 $$
 \begin{aligned}
 &\frac{5}{6}+\frac{\nu_{2}}{12}<\gamma<2+\frac{\nu_{1}}{2}\;\;\;\mbox{if}\;\;N=3,\\
 &\frac{5}{6}+\frac{\nu_{2}}{12}<\gamma<\frac{5}{6}+\frac{7}{12}\nu_{1}\;\;\;\mbox{if}\;\;N=3,\\
 &\frac{1}{4}+\frac{\nu_{2}}{8}<\gamma\;\;\;\mbox{if}\;\;N=2.
 \end{aligned}
 $$
 
\item $0<\nu_{1}< 2$ and:
$$
 \begin{aligned}
 &\frac{5}{6}+\frac{\nu_{2}}{12}<\gamma<\frac{(4-\nu_{1})(1+\nu_{1})}{2-\nu_{1}}\;\;\;\mbox{if}\;\;N=3,\\
 &\frac{5}{6}+\frac{\nu_{2}}{12}<\gamma<\frac{5}{6}+\frac{7}{12}\nu_{1}\;\;\;\mbox{if}\;\;N=3,\\
 &\frac{1}{4}+\frac{\nu_{2}}{8}<\gamma\;\;\;\mbox{if}\;\;N=2.
 \end{aligned}
 $$
%\begin{equation}
%\begin{aligned}
%&2\gamma-\frac{2}{3}-\frac{\nu_{1}}{6}< (2-\frac{2-\nu_{1}}{6(1+\nu_{1})})\gamma\;\;\Leftrightarrow\;\;\gamma<2+\frac{\nu_{1}}{2}.\\
%&2\gamma-1+\frac{1}{N}-\frac{\nu_{2}}{2N}>1
%\end{aligned}
 \end{itemize}
 \label{lemcondi}
\end{lemme}
{\bf Proof:} Let us come back to the inequality (\ref{gaina}), we have  $\forall \delta\in(0,2)$:
\begin{equation}
\begin{aligned}
&\int_{\R^{N}}\rho_{\ep}\frac{1+|u_{\ep}|^{2}}{2}\ln(1+|u_{\ep}|^{2})(t,x)dx+\nu\int^{t}_{0}\int_{\R^{N}}\mu(\rho_{\ep})(1+\ln(1+|u_{\ep}|^{2}))|D u_{\ep}|^{2}(t,x)dxdt\\
&\leq C\int^{t}_{0}\int_{\R^{N}}\mu(\rho_{\ep})|\n u_{\ep}|^{2}(t,x)dxdt+C_{\delta}\ep^{2}\int^{t}_{0}\big(\int_{\R^{N}}(\frac{\rho_{\ep}^{2\gamma-\frac{\delta}{2}}}{\mu(\rho_{\ep})})^{\frac{2}{2-\delta}}dx\big)dt.
\end{aligned}
\label{24}
\end{equation}
Let us deal with this inequality in order to prove the uniform bound of $\rho_{\ep}\frac{1+|u_{\ep}|^{2}}{2}\ln(1+|u_{\ep}|^{2})$ in $L^{\infty}_{T}(L^{1}(\R^{N}))$ for any $T>0$. To do this we shall estimate the right hand side of (\ref{24}). Using the inequality (\ref{13}) we obtain that:
\begin{equation}
\begin{aligned}
&\frac{\rho^{2\gamma-\frac{\delta}{2}}}{\mu(\rho)}\leq \frac{1}{C}\rho^{2\gamma-\frac{\delta}{2}-1+\frac{1}{N}-\frac{\nu_{1}}{2N}}\;\;\;\forall\rho>1,\\
&\frac{\rho^{2\gamma-\frac{\delta}{2}}}{\mu(\rho)}\leq \frac{1}{C}\rho^{2\gamma-\frac{\delta}{2}-1+\frac{1}{N}-\frac{\nu_{2}}{2N}}\;\;\;\forall\rho\leq 1.
\end{aligned}
\label{imptech}
\end{equation}
Let us distinguish two cases when $\nu_{1}\geq 2$ and when $0<\nu_{1}<2$.\\
\\
$\bullet$ $\nu_{1}\geq 2$ and $N=3$\\
From the lemma \ref{lempression}, we know that $\e\rho_{\e}^{\gamma}$ is uniformly bounded in $L^{\frac{5}{3}}((0,T)\times\R^{N})\cap L^{\infty}_{T}(L^{1}(\R^{N}))$ for $N=3$ which implies that it exists $C>0$ such that:
\begin{equation}
\e^{\frac{5}{3}}\int^{t}_{0}\int_{\R^{N}}\rho_{\e}^{\frac{5}{3}\gamma}dxdt\leq C.
\label{crucialconv}
\end{equation}
In particular it implies that for $\delta$ small enough $\e^{2}(1_{\{\rho_{\ep}\geq 1\}}\frac{\rho_{\ep}^{2\gamma-\frac{\delta}{2}}}{\mu(\rho_{\ep})})^{\frac{2}{2-\delta}}$ is uniformly bounded in $L^{1}((0,T)\times\R^{N})$ under the conditions that:
\begin{equation}
\begin{aligned}
&2\gamma-1+\frac{1}{N}-\frac{\nu_{1}}{2N}< \frac{5}{3}\gamma\;\;\Leftrightarrow\;\;\gamma<2+\frac{\nu_{1}}{2}.
\end{aligned}
\label{p1condi}
\end{equation}
Indeed by Tchebytchev lemma we have:
$$|\{x,\;|\rho_{\e}(t,x)|\geq 1\}|\leq \|\rho_{0}\|_{L^{1}(\R^{N})}.$$
We choose $p$ such that $p(2\gamma-\frac{\delta}{2}-1+\frac{1}{N}-\frac{\nu_{1}}{2N})\frac{2}{2-\delta}= \frac{5}{3}\gamma$ with $\delta$ small enough. By H\"older's inequality we have:
$$
\begin{aligned}
&\e^{2}\int^{T}_{0}\int_{\R^{N}} 1_{\{\rho_{\e}\geq 1\}}\rho^{(2\gamma-\frac{\delta}{2}-1+\frac{1}{N}-\frac{\nu_{1}}{2N})\frac{2}{2-\delta}}dxdt\leq \e^{2}\int^{T}_{0}(\int_{\R^{N}} \rho^{\frac{5}{3}\gamma}dx)^{\frac{1}{p}}
\|\rho_{0}\|_{L^{1}(\R^{N})}^{\frac{1}{p^{'}}}dt,\\
&\hspace{6cm}\leq \e^{2-\frac{5}{3p}}   (T\|\rho_{0}\|_{L^{1}(\R^{N})})^{\frac{1}{p^{'}}}  (\int^{T}_{0}\int_{\R^{N}}\e^{\frac{5}{3}} \rho^{\frac{5}{3}\gamma}dx dt)^{\frac{1}{p}}.
\end{aligned}
$$
with $\frac{1}{p^{'}}=1-\frac{1}{p}$ and $p>1$. It implies by (\ref{crucialconv}) that $\e^{2}1_{\{\rho_{\e}\geq 1\}}\rho^{(2\gamma-\frac{\delta}{2}-1+\frac{1}{N}-\frac{\nu_{1}}{2N})\frac{2}{2-\delta}}$ is uniformly bounded in $L^{1}((0,T)\times\R^{N})$ under the hypothesis (\ref{p1condi}).\\
Let us deal now with the term $\e^{2}1_{\{\rho_{\e}\leq 1\}}\rho^{(2\gamma-\frac{\delta}{2}-1+\frac{1}{N}-\frac{\nu_{1}}{2N})\frac{2}{2-\delta}}$, it suffices to assume that:
\begin{equation}
\begin{aligned}
&2\gamma-1+\frac{1}{N}-\frac{\nu_{2}}{2N}> 1\;\;\Leftrightarrow\;\;\gamma>\frac{5}{6}+\frac{\nu_{2}}{12}.
\end{aligned}
\label{p2condi}
\end{equation}
Indeed we know that $1_{\{\rho_{\e}\leq 1\}}\rho_{\e}$ is bounded in $L^{\infty}_{T}(L^{1}(\R^{N}))\cap L^{\infty}_{T}(L^{\infty}(\R^{N}))$ which insures the uniform bound of $1_{\{\rho_{\e}\leq 1\}}\rho^{(2\gamma-\frac{\delta}{2}-1+\frac{1}{N}-\frac{\nu_{1}}{2N})\frac{2}{2-\delta}}$ in $L^{\infty}_{T}(L^{1}(\R^{N}))$ by interpolation. This achieves the proof of the case $N=3$.\\
\\
The last situation consists when $N=3$ in using the lemma \ref{tech} when $2\gamma-1+\frac{1}{N}-\frac{\nu_{1}}{2N}\leq 1+\frac{3\nu_{1}}{N}$ which is equivalent to $\gamma\leq \frac{5}{6}+\frac{7}{12}\nu_{1}$. Indeed via the lemma \ref{tech} we know that $\rho_{\e}^{1+\frac{\nu_{1}}{2N}}$ is uniformly bounded in $L^{\infty}_{T}(L^{6}(\R^{N}))$. It implies that
$\e^{2}1_{\{\rho_{\e}\geq 1\}}\rho^{(2\gamma-\frac{\delta}{2}-1+\frac{1}{N}-\frac{\nu_{1}}{2N})\frac{2}{2-\delta}}$ is uniformly bounded in $L^{\infty}_{T}(L^{1}(\R^{N}))$ when
$2\gamma-1+\frac{1}{N}-\frac{\nu_{1}}{2N}\leq 1+\frac{3\nu_{1}}{N}$ . In order to bounde $\e^{2}1_{\{\rho_{\e}\leq 1\}}\rho^{(2\gamma-\frac{\delta}{2}-1+\frac{1}{N}-\frac{\nu_{1}}{2N})\frac{2}{2-\delta}}$ we use the same hypothesis than in the previous case.
%Using the hypothesis  of the theorem \ref{theo3} $\gamma<3-\frac{3}{N}+\frac{3\nu_{1}}{2N}=2+\frac{\nu_{1}}{2}$ when $N=3$ and since $\gamma>\frac{1}{2}-\frac{1}{2N}+\frac{\nu_{2}}{4N}=\frac{1}{3}+\frac{\nu_{2}}{12}$, it implies by interpolation that $\frac{\rho_{\ep}^{2\gamma-\frac{\delta}{2}}}{\mu(\rho_{\ep})}1_{\{\rho_{\ep}\leq 1\}}$ is uniformly bounded in $L^{\infty}_{T}(L^{1}(R^{N}))$ for $\delta$ small enough since $\rho_{\e}$ is uniformly bounded in $L^{\infty}_{T}(L^{1}(R^{N}))$. And since $\frac{\rho_{\ep}^{2\gamma-\frac{\delta}{2}}}{\mu(\rho_{\ep})}1_{\{\rho_{\ep}\leq 1\}}$ is uniformly bounded in $L^{\infty}((0,T)\times\R^{N})$ we conclude that $\e^{2}(1_{\{\rho_{\ep}\leq 1\}}\frac{\rho_{\ep}^{2\gamma-\frac{\delta}{2}}}{\mu(\rho_{\ep})})^{\frac{2}{2-\delta}}$ is uniformly bounded in $L^{1}((0,T)\times\R^{N})$.\\
\\
\\
$\bullet$ $\nu_{1}\geq 2$, $N=2$\\
In this case, the situation is quite simple, indeed we know via the lemma \ref{tech} that $\rho_{\e}$ is bounded in $L^{\infty}_{T}(L^{q}(\R^{N}))$ for any $q\geq 1$. In particular it implies that $1_{\{\rho_{\e}\geq 1\}}\rho^{2\gamma-\frac{\delta}{2}-1+\frac{1}{N}-\frac{\nu_{1}}{2N}}$ is bounded in $L^{\infty}_{T}(L^{1}(\R^{N}))$ without any specific condition. However we shall require a hypothesis for dealing with the term $1_{\{\rho_{\e}\leq 1\}}\rho^{2\gamma-\frac{\delta}{2}-1+\frac{1}{N}-\frac{\nu_{1}}{2N}}$ which is similar to the previous section:
\begin{equation}
\begin{aligned}
&2\gamma-1+\frac{1}{N}-\frac{\nu_{2}}{2N}> 1\;\;\Leftrightarrow\;\;\gamma>\frac{1}{4}+\frac{\nu_{2}}{8}.
\end{aligned}
\label{p3condi}
\end{equation}
$\bullet$ $0<\nu_{1}< 2$, $N=2$\\
The proof in this situation is exactly the same than in the previous case by using the lemma \ref{tech}. We need only:
\begin{equation}
\begin{aligned}
&2\gamma-1+\frac{1}{N}-\frac{\nu_{2}}{2N}> 1\;\;\Leftrightarrow\;\;\gamma>\frac{1}{4}+\frac{\nu_{2}}{8}.
\end{aligned}
\label{p4condi}
\end{equation}
$\bullet$ $0<\nu_{1}< 2$, $N=3$\\
%\subsubsection*{Step 2: Convergence of the pressure}
\\
Via the lemma \ref{lempression} we have seen that $\e\rho_{\e}^{\gamma}$ is bounded in $L^{r_{1}}_{T}(L^{r_{1}}(\R^{N}))$
with $r_{1}=2-\frac{2-\nu_{1}}{6(1+\nu_{1})}$. By using exactly the same argument than in the previous case we have two possibility to bound the last term on the right hand side of (\ref{24}):
\begin{equation}
\begin{aligned}
&2\gamma-\frac{2}{3}-\frac{\nu_{1}}{6}< (2-\frac{2-\nu_{1}}{6(1+\nu_{1})})\gamma\;\;\Leftrightarrow\;\;\gamma<2+\frac{\nu_{1}}{2}.\\
&2\gamma-1+\frac{1}{N}-\frac{\nu_{2}}{2N}>1
\end{aligned}
\label{p5condi}
\end{equation}
or
\begin{equation}
\begin{aligned}
&\frac{5}{6}+\frac{\nu_{2}}{12}<\gamma \frac{5}{6}+\frac{7\nu_{1}}{12}.
\end{aligned}
\label{p5condi}
\end{equation}
It achieves the proof of the lemma \ref{lemcondi}.
{\hfill $\Box$}\\
\\
We have now proved that $\rho_{\ep}\frac{1+|u_{\ep}|^{2}}{2}\ln(1+|u_{\ep}|^{2})$ is uniformly bounded in $L^{\infty}_{loc}(L^{1}(\R^{N}))$. We can then pass to the limit when $\e$ goes to $0$, more precisely by using lemmas \ref{rho}, \ref{ruse2}  we show that $\rho_{\ep}$, $\rho_{\ep}u_{\ep}$ and $\sqrt{\rho_{ep}}u_{\ep}\otimes \sqrt{\rho_{ep}}u_{\ep}$ converges in distribution sense to $\rho$, $\rho u$ and  $\sqrt{\rho} u\otimes \sqrt{\rho} u$ and lemma \ref{ruse3} give us the convergence in distribution sense of the diffusion term. 
Furthermore the lemmas \ref{rho} and \ref{ruse2} give us the following desired strong convergence:
\begin{itemize}
\item $\rho_{\ep}$ converges strongly to $\rho$ in $C([0,T],L^{1+\alpha}_{loc}(\R^{N}))$ for $\alpha$ small enough when $N=3$.
\item $\rho_{\ep}$ converges strongly to $\rho$ in $C([0,T],L^{p}_{loc}(\R^{N}))$ for any $p\geq 1$  when $N=2$.
\item $\sqrt{\rho_{\ep}}u_{\e}$ converges strongly to $\sqrt{\rho}u$ in $L^{2}_{loc}((0,T)\times\R^{N}))$ for any $T>0$ .
\end{itemize}
It remains only to deal with the term $\e\n\rho_{\e}^{\gamma}$ and to prove that it converges in distribution sense to $0$. 
\begin{lemme}
Let us distinguish two cases:
\begin{itemize}
\item When $N=3$\\
 $\e^{\alpha}\rho_{\e}^{\gamma}$ converges strongly to $0$ in $L^{\infty}_{T}(L^{1}_{loc}(\R^{N}))$ for any $\alpha>0$ when $\e$ goes to $0$ for:
$$
\begin{aligned}
&\gamma<1+\nu_{1}.
\end{aligned}
$$ 
 $\e^{\alpha}\rho_{\e}^{\gamma}$ converges strongly to $0$ in $L^{\frac{5}{3}-\alpha}((0,T)\times\R^{N}))$ for $\nu_{1}\geq 2$ and for any $\alpha>0$ small enough when $\e$ goes to $0$.\\
  $\e^{\alpha}\rho_{\e}^{\gamma}$ converges strongly to $0$ in $L^{1+\alpha}_{T}(L^{r_{1}-\alpha}(\R^{N}))$ for $0<\nu_{1}<2$ and for any $\alpha>0$ small enough when $\e$ goes to $0$.
  \item When $N=2$\\
 $\e^{\alpha}\rho_{\e}^{\gamma}$ converges strongly to $0$ in $L^{\infty}_{T}(L^{1}_{loc}(\R^{N}))$ for any $\alpha>0$ when $\e$ goes to $0$.
\end{itemize}
\label{lepress}
\end{lemme}
{\bf Proof:} When $N=2$ we know via the lemma \ref{tech} that $\rho_{\e}$ is bounded in $L^{\infty}_{T}(L^{q}(\R^{N}))$ for any $q\geq 1$. It implies trivially that  $\e^{\alpha}\rho_{\e}^{\gamma}$ converges strongly to $0$ in $L^{\infty}_{T}(L^{1}_{loc}(\R^{N}))$ for any $\alpha>0$ when $\e$ goes to $0$.\\
When $N=3$ we are in a similar situation when $\gamma<1+\nu_{1}$ via the lemma \ref{tech}. If $\nu_{1}\geq 2$ we have seen in the lemma \ref{lempression} that $\e\rho_{\e}^{\gamma}$ is uniformly bounded in $L^{\frac{5}{3}}((0,T)\times\R^{N}))$, combining this result with the fact that $\e\rho_{\e}^{\gamma}$ is uniformly bounded in $L^{\infty}_{T}(L^{1}(\R^{N}))$ and an interpolation argument we obtain the result that we wish. In the case $0<\nu_{1}<2$ we apply a similar argument with $r_{1}$ by using the lemma \ref{lempression}. It concludes the proof of the lemma. {\hfill $\Box$}
\\
\\
Furthermore in the same spirit by using exactly the same arguments than in the proof of theorem \ref{theo2} and the previous estimate on the pressure we prove also the stability of the global weak solutions for the system (\ref{a01}).
\subsection{Proof of the corollary \ref{cor}}
From the previous theorem we know that $(\rho_{\e},u_{\e})$ converges to a quasi solution $(\rho,u)$ and that $\rho_{\e}$ converges strongly to $\rho$ in $C([0,T],L^{1}_{loc}(\R^{N}))$. In particular when $\mu(\rho)=\mu\rho^{\alpha}$ if we assume that there exists a unique global quasi solution we know via the theorem \ref{theo2} that this quasi solution verifies the porous media or the fast diffusion equation in function of $\alpha$ inasmuch as $\rho$ is solution of (\ref{P}).\\
In particular when we assume that the initial density $\rho_{0}$ has a compact support, it implies that when $\alpha>1$ the support of the density $\rho$ remains bounded along the time. Indeed it consists merely of not ing that we can find a delayed Barrenblatt solution centered for instance at $0$ that lies on the top of $\rho_{0}$, it means:
$$0\leq \rho_{0}(x)\leq U_{m}(\tau,x)\;\;\;\;\forall x\in\R^{N},$$
with $m$ large enough and $\tau>0$. By theorem \ref{theo2.4} and the maximum principle we know that:
$$0\leq \rho(t,x)\leq (t+\tau)^{-\gamma_{1}}F(\frac{x}{(t+\tau)^{\beta}})\;\;\;\;\forall x\in\R^{N},$$
with $F(x)=(C-\frac{(\alpha-1)\gamma_{1}}{2\alpha}|x|^{2})_{+}^{\frac{1}{\alpha-1}}$. In particular it implies an information on the expansion of the support of the solution since we observe that the support of $\rho(t,\cdot)$ is included in a set $E(t)=C B(0,M(t+\tau)^{\frac{\beta}{2}}$ with $C>0$, $M>0$ independent of $t$.\\
Let us prove now that  $\rho_{\e}$ converges strongly to $\rho$ in $C([0,T],L^{1}(\R^{N}))$. We know for the moment that $\rho_{\e}$ converges strongly to $\rho$ in $C([0,T],L^{1}_{loc}(\R^{N}))$, it suffices to consider $K$ a compact set large enough such that for any $t\in[0,T]$ we have $\mbox{supp}\rho(t,\cdot)\subset K$, we have then:
\begin{equation}
\|\rho(t,\cdot)-\rho_{\e}(t,\cdot)\|_{L^{1}(K)}\rightarrow_{\e\rightarrow 0}0.
\label{fort2}
\end{equation}
Now we have by conservation of the mass and the fact that $\rho(t,\cdot)=0$ in $K^{c}$ for any $t\in[0,T]$:
\begin{equation}
\begin{aligned}
\|\rho(t,\cdot)-\rho_{\e}(t,\cdot)\|_{L^{1}(K^{c})}&= \|\rho_{\e}(t,\cdot)\|_{L^{1}(K^{c})},\\
&=\|\rho_{0}\|_{L^{1}(\R^{N})}-\|\rho_{\e}(t,\cdot)\|_{L^{1}(K)}.
\end{aligned}
\label{fort1}
\end{equation}
In particular since  $\rho_{\e}$ converges strongly to $\rho$ in $C([0,T],L^{1}(K))$, (\ref{fort1}) implies that $\|\rho_{\e}(t,\cdot)\|_{L^{1}(K)}$ converges uniformly on $[0,T]$ to $\|\rho(t,\cdot)\|_{L^{1}(K)}$ when $\e$ goes to $0$.  But since $\|\rho(t,\cdot)\|_{L^{1}(K)}=\|\rho_{0}\|_{L^{1}(\R^{N})}$ for any $t\in[0,T]$ (indeed the support of $\rho(t,\cdot)$ is completely included in $K$), it induces that $\|\rho(t,\cdot)-\rho_{\e}(t,\cdot)\|_{L^{1}(K^{c})}$ converges uniformly on $[0,T]$ to $0$ when $\e$ goes to $0$.
\\
Finally we have by using (\ref{fort1}):
$$\|\rho(t,\cdot)-\rho_{\e}(t,\cdot)\|_{L^{1}(\R^{N})}\leq \|\rho(t,\cdot)-\rho_{\e}(t,\cdot)\|_{L^{1}(K)}+\|\rho_{0}\|_{L^{1}(\R^{N})}-\|\rho_{\e}(t,\cdot)\|_{L^{1}(K)}.$$
It implies that $\|\rho(t,\cdot)-\rho_{\e}(t,\cdot)\|_{L^{1}(\R^{N})}$ converges uniformly on $[0,T]$ to $0$ when $\e$ goes to $0$ and we have shown that $\rho_{\e}$ converges strongly to $\rho$ in $C([0,T],L^{1}(\R^{N}))$. In particular it implies that for $\e$ small enough $\rho_{\e}$ is the sum of a solution with compact support on $[0,T]$ and of a term of small $L^{1}$ norm. In this sense we can claim that the propagation speed of the free boundary of $\rho_{\e}$ is not so far to be finite at a small $L^{1}$ perturbation.\\
This implies in particular by interpolation that $\rho_{\e}$  converges strongly to $\rho$ in $C([0,T],L^{p}(\R^{N}))$ for any $T>$, any $p\leq 1+\alpha$ with $\alpha$ small enough if $N=3$ and with $p\geq 1$ if $N=2$.\\
In particular since via the theorem \ref{theo2.6}, we have:
$$\|\rho(t)\|_{L^{p}(\R^{N})}\leq Ct^{-\sigma_{p }}\|\rho_{0}\|^{\alpha_{p}}_{L^{1}(\R^{N})},$$
with $\sigma_{p}=\frac{N(\alpha-1)+2p}{(N(\alpha-1)+2)p}$ and $\alpha_{p}=\frac{N(p-1)}{(N(\alpha-1)+2)p}$. It shows that up a remainder term of small norm in $L^{p}$, the $L^{p}$ norm of  $\rho_{\e}$ decrease in time for small $\e$. It means that in some sense the density is subjected to a damping effect in time for the $L^{p}$ norm which is very surprising since this effect seems purely non linear.\\
%{\bf REGARDER !!!!}Let us prove in fact that for $p>1$ $\rho_{\e}$  converges strongly to $\rho$ in $C([0,+\infty),L^{p}(\R^{N}))$. We have then:
%$$\|\rho(t,\cdot)-\rho_{\e}(t,\cdot)\|_{L^{p}(K)}\leq \|\rho\|_{L^{p}(K)}+\|\rho_{0}\|_{L^{p}(\R^{N})}-.$$
\\
Let us deal now with the time asymptotic behavior of $\rho_{\e}$. We expect that $\rho_{\e}(t,\cdot)$ goes asymptotically in time to the Barrenblatt solution $U_{m}$ of (\ref{P}) of mass $\|\rho_{0}\|_{L^{1}(\R^{N})}=m$. We have then:
\begin{equation}
\begin{aligned}
&\|U_{m}(t,\cdot)-\rho_{\e}(t,\cdot)\|_{L^{1}(\R^{N})}\leq \|U_{m}(t,\cdot)-\rho(t,\cdot)\|_{L^{1}(\R^{N})}+\|\rho(t,\cdot)-\rho_{\e}(t,\cdot)\|_{L^{1}(\R^{N})}.
\end{aligned}
\label{fort3}
\end{equation}
Via the theorem \ref{theo2.7} we know that $\|U_{m}(t,\cdot)-\rho(t,\cdot)\|_{L^{1}(\R^{N})}$ converges asymptotically to $0$ when $t$ goes to $+\infty$. The second term converges also to $0$ when $\e$ goes to $0$ since we have shown that $\rho_{\e}$ converges strongly to $\rho$ in $C([0,T],L^{1}(\R^{N}))$ for any $T>0$. In particular it implies that for any $\alpha>0$ it exists $T>0$ such that:
$$\|U_{m}(t,\cdot)-\rho(t,\cdot)\|_{L^{1}(\R^{N})}\leq\alpha\;\;\;\forall t>T.$$
Furthermore it exists $\e_{0}>0$ such that for all $\e\leq\e_{0}$ we have:
$$\|\rho(t,\cdot)-\rho_{\e}(t,\cdot)\|_{L^{1}(\R^{N})}\leq\alpha\;\;\;\forall t\in[0,nT],\;\;\mbox{with}\;n\in\mathbb{N}.$$
It implies that for all $\alpha>0$ it exits $T>0$  such that for all $n\in\mathbb{N}$ it exits $\e_{0}>0$ such that for all $0<\e\leq\e_{0}$ we have:
$$\|U_{m}(t,\cdot)-\rho_{\e}(t,\cdot)\|_{L^{1}(\R^{N})}\leq 2\alpha\;\;\;\forall t\in[T,nT].$$
In this sense we observe that for $\e$ small enough the solution $\rho_{\e}$ tends to converge asymptotically to a Barrenblatt solution of mass $\|\rho_{0}\|_{L^{1}(\R^{N})}$.
%{\bf UTILISER L'INVARIANCE PAR SCALING POUR LE GLOBAL ET SE RAMENER AU CAS TEMPS FINI!!!!}
%{\bf QUASI PRINCIPE DU MAXIMUM DONC QUASI UNICITE!!!}
%ON DOIT AVOIR UNE CONVERGENCE FORTE GLOBALE DANS $L^{p}$ PUISEQUE CA DECROIT!!!!
%\begin{lem}
%\end{lem}
%We can easily observe that via the energy estimates the right hand side of (\ref{24}) is uniformly bounded in $\ep$. The last step corresponds to use the same compactness argument than in \cite{MV} to show that $(\rho_{\ep},u_{\ep})$ converges in distribution sense to a quasi-solution $(\rho,u)$ when $\ep$ goes to $0$ with initial data $(\rho_{0},u_{0})$. Let us point out that $\ep\rho^{\gamma}_{\ep}$ goes to $0$ in distribution sense. Indeed it suffice to observe that $\rho_{\ep}$ and $\sqrt{\rho_{\ep}}\n\va(\rho_{\ep})$ are uniformly bounded in $\ep$ respectively in $L^{\infty}((0,T),L^{1}(\R^{N})$ and $L^{\infty}((0,T),L^{2}(\R^{N})$ for any $T>0$. In particular by Sobolev embedding and interpolation (for $\gamma$ not so large as in \cite{MV}) we obtain that $\rho_{\ep}^{\gamma}$ with $\alpha>0$ is uniformly bounded in $L^{1}_{loc}$ what means that $\ep\rho^{\gamma}_{\ep}$ goes to $0$.
\null{\hfill $\Box$}
%\section{Vanishing viscosity process for the quasi-solutions}
%Let us start with the case of shallow water viscosity coefficient $\mu(\rho)=\mu\rho$ and $\lambda(\rho)=0$. We have then:
%$$
%\begin{cases}
%\begin{aligned}
%&\p_{t}\rho-\mu\D\rho=0,\\
%&\rho(0,\cdot)=\rho_{0}.
%\end{aligned}
%\end{cases}
%$$
%It implies that:
%$$\rho(t,x)=\frac{1}{(4\pi\mu t)^{\N}}\int_{\R^{N}}\rho_{0}(y)\exp (-\frac{|x-y|^{2}}{4\mu t})dy.$$
%Next we have:
%$$
%\begin{aligned}
%u(t,x)&=-\mu\n\ln\rho(t,x),\\
%&=\frac{\int_{\R^{N}}(\frac{x-y}{2 t})\rho_{0}(y)\exp (-\frac{|x-y|^{2}}{4\mu t})dy}{\int_{\R^{N}}\rho_{0}(y)\exp (-\frac{|x-y|^{2}}{4\mu t})dy}.
%\end{aligned}
%$$
%And we recall that:
%$$u_{0}(y)=-\mu\n\ln\rho_{0}(y).$$
%We obtain then:
%$$\rho_{0}(y)$$
%\subsection*{Acknowledgement}
\subsection*{Acknowledgement} I would like to thank Alexis Vasseur for a very interesting question this problem in Padova.

\end{document}